\documentclass[12pt,reqno]{amsart}
\usepackage[LGR,T1]{fontenc}
\usepackage[utf8]{inputenc}

\title[]{Feynman-Kac formula under a finite entropy condition}
\date{June, 2022}


\author[]{Christian Léonard}
\address{Modal’X, UMR CNRS 9023. UPL, Univ Paris Nanterre, F92000 Nanterre, France}
\email{christian.leonard@math.cnrs.fr}
 \thanks{This research is partially  granted by Labex MME-DII (ANR-11-LBX-0023).\\
 Data sharing not applicable to this article as no datasets were generated or analysed during the current study.}

\usepackage{amssymb, amsmath, amsfonts, latexsym, enumerate}
\usepackage{comment}
\usepackage{mathabx} 
\usepackage[usenames,dvipsnames]{pstricks}
\usepackage{epsfig}
\RequirePackage[colorlinks,linkcolor=blue,citecolor=blue,urlcolor=blue]{hyperref}
\usepackage[english]{babel}

\usepackage{xcolor}

\newtheorem{theorem}[equation]{Theorem}
\newtheorem{lemma}[equation]{Lemma}
\newtheorem{proposition}[equation]{Proposition}
\newtheorem{corollary}[equation]{Corollary}

\newtheorem{definition}[equation]{Definition}

\newtheorem{hypotheses}[equation]{Hypotheses}

\theoremstyle{remark}
\newtheorem{remark}[equation]{Remark}
\newtheorem{remarks}[equation]{Remarks}

\numberwithin{equation}{section}

\newcommand{\RR}{\mathbb{R}}
\newcommand{\Rn}{\mathbb{R}^n}
\newcommand{\PP}{\mathbb{P}}

\newcommand{\1}{\mathbf{1}}

\newcommand{\ttimes}{\!\times\!}

\newcommand\pf{_{\#}}

\newcommand{\Leb}{\mathrm{Leb}}

\renewcommand{\ae}{\textrm{-}{a.e.}}
\newcommand{\Id}{\mathrm{Id}}
\newcommand{\scal}{\!\cdot\!}

\DeclareMathOperator{\dom}{dom}

\DeclareMathOperator{\supp}{supp}

\DeclareMathOperator{\tr}{tr}
\DeclareMathOperator{\MP}{MP}

\DeclareMathOperator{\range}{range}

\newcommand{\boulette}[1]{$\bullet$\ Proof of #1.}
\newcommand{\Boulette}[1]{\par\medskip\noindent $\bullet$\ Proof of #1.}
\newcommand{\sbt}{\,\begin{picture}(-1,1)(-1,-3)\circle*{3}\end{picture}\ }

\newcommand\seq[2]{(#1_#2)_{#2\ge1}}

\newcommand\Lim[1]{\lim_{#1\rightarrow\infty}}

\newcommand\Limh{\lim_{h\to 0^+}}
\newcommand\II[2]{\int_{[#1,#2]}}

\newcommand{\eqlaw}{\overset{\textrm{Law}}=}

\newcommand{\ud}{\frac{1}{2}}

\newcommand{\ii}{[0,T]}
\newcommand\XX{ \mathcal{X}}

\newcommand\XXb{{\overline{\XX}}}
\newcommand\iX{{\ii\ttimes\XX}}

\newcommand{\ZZ}{\Rn}
\newcommand\ZZZ{\ZZ\times\ZZ}

\newcommand{\iZ}{\ii\times\ZZ}
\newcommand\OO{\Omega}

\newcommand\Ob{\overline{\OO}}

\newcommand\PX{\mathrm{P}(\XX)}

\newcommand\MX{\mathrm{M}(\XX)}

\newcommand\PO{\mathrm{P}(\OO)}
\newcommand\MO{\mathrm{M}(\OO)}

\newcommand\Iii{\int_{\ii}}

\newcommand\IZ{\int_{\ZZ}}

\newcommand\IiZ{\int_{\iZ}}

\newcommand\Xb{\overline{X}}
\newcommand\Pb{\overline{P}}

\newcommand\Qb{\overline{Q}}

\newcommand{\gradt}{\widetilde{\nabla}}
\newcommand{\Gt}[2]{ \gradt ^{ #1, {\scriptscriptstyle #2}}}

\newcommand\LQ{\mathcal{L}^Q}

\newcommand{\Lf}{\overrightarrow{L}}
\newcommand{\Lb}{\overleftarrow{L}}
\newcommand{\LLf}{\overrightarrow{ \mathcal{L}}}
\newcommand{\LLF}{ \mathcal{L}} 
\newcommand{\LLb}{\overleftarrow{ \mathcal{L}}}
\newcommand{\LLL}[2]{ \mathcal{L} ^{ #1, {\scriptscriptstyle #2}}}
\newcommand{\LLLf}[2]{ \LLf ^{ #1, {\scriptscriptstyle #2}}}
\newcommand{\LLLb}[2]{ \LLb ^{ #1, {\scriptscriptstyle #2}}}

\newcommand{\GaF}{\Gamma}

\newcommand{\vf}{ \vv}

\renewcommand{\bf}{ \beta}

\newcommand{\pb}{\bar{\pp}}
\newcommand{\qb}{\bar{\qq}}
\newcommand{\rb}{\bar{\rr}}

\newcommand{\mm}{ \mathsf{m}}

\newcommand{\pp}{ \mathsf{p}}
\newcommand{\qq}{ \mathsf{q}}
\newcommand{\rr}{ \mathsf{r}}
\newcommand{\vv}{ \mathsf{v}}

\newcommand{\ww}{ \mathsf{W}}

\newcommand{\yy}{ \mathsf{Y}}

\newcommand{\bbf}[2]{ \bf ^{ #1|#2}}

\renewcommand{\aa}{ \mathsf{a}}
\newcommand{\BB}{ \mathsf{b}}
\newcommand{\BBf}{ \BB}

\newcommand{\cc}{ \mathsf{c}}
\renewcommand{\gg}{ \mathsf{g}}

\newcommand{\mU}{\mm^U}
\newcommand{\ms}{\mm _{ \ast}}
\newcommand{\vU}{\vv^U}
\newcommand{\vs}{\vv _{ \ast}}
\newcommand{\RU}{{R^U}}
\newcommand{\Us}{{U _{ \ast}}}
\newcommand{\Uo}{U^\diamond }
\newcommand{\UUo}{ \mathcal{U}^\diamond }

\begin{document}

\begin{abstract} 
Motivated by entropic optimal transport, we are interested in the Feynman-Kac formula associated to     the  parabolic equation $	( \mathsf{L}+V)g
	=0$ with a final nonnegative boundary condition and a Markov generator $ \mathsf{L}:= \partial_t + \BB\scal \nabla + \Delta _{ \aa}/2$. 
	It is well-known that  when the drift $\BB$, the diffusion matrix $\aa$ and the scalar potential $V$ are regular enough and not growing too fast, the classical  solution $g$ of this PDE,
 is represented by the Feynman-Kac formula
 $
 g_t(x)=E_R[\exp \left( \int _{[t,T]} V(s,X_s)\,ds\right) g(X_T)\mid X_t=x]
 $ where $R$ is the Markov measure with  generator $ \mathsf{L}$. 

We do not assume  that  $g$, $\BB$ and $V$ are  regular, and only require  that  their growth is controlled by a finite entropy condition. These hypotheses are less restrictive than  the standard assumptions of the theory of viscosity solutions,  and allow for instance $V$ to belong to some Kato class.  We prove that $g$ defined by  the Feynman-Kac formula belongs to the domain of the \emph{extended} generator $ \mathcal{L}$ of the Markov measure $R$ and  satisfies the \emph{trajectorial} identity: $
[(\mathcal{L}   +V)g] (t,X_t)=0,\ dtdP\ae
$
where the path measure $P$ is defined by   $	
 P:= f(X_0)\exp \left( \Iii V(t,X_t)\,dt\right) g(X_T)\ R,
$	
with    $ f:\ZZ\to [0, \infty)$  another nonnegative function.  We also show that the forward drift $\BB^P$ of $P$ satisfies $\BB^P(t,X_t)=[\BB+ \aa \gradt \log g](t,X_t),$ $dtdP\ae,$ where $\gradt$  is some extension of the standard derivative.

Our  probabilistic approach relies on stochastic derivatives, semimartingales, Girsanov's theorem  and the Hamilton-Jacobi-Bellman equation satisfied by $\log g$. 
\end{abstract}

\maketitle 
\tableofcontents

\section{Introduction}

Let us call for practical use in this article,  Feynman-Kac equation   the linear parabolic equation
\begin{align}\label{eq-FK}
\left\{
\begin{array}{ll}
( \partial_t +\mathsf{A}	+V)g
	=0,
	&\quad  0\le t<T,\\
 g(T,\sbt)= g_T, &\quad t=T,
\end{array}
\right.
\tag{FK}
\end{align}
where the numerical function $g:\iZ\to \RR$ is the unknown, $V:\iZ\to\RR$ is a scalar potential seen as a multiplicative operator, $ g_T:\ZZ\to[0, \infty)$ is a given \emph{nonnegative} function, and $ \mathsf{A}$ is the Markov diffusion generator
\begin{align*}
 \mathsf{A}:=\BB\scal \nabla + \Delta _{ \aa}/2
\end{align*}
whose coefficients are a velocity field
$\BB:\iZ\to\ZZ$  and a diffusion matrix field  $\aa:\iZ\to S_+$  taking its values in the set $S_+$ of \emph{nonnegative} symmetric  $n\times n$-matrices. We denote for simplicity $ \Delta_\aa:= \sum _{ 1\le i,j\le n} a _{ ij} \partial_i \partial_j$ where $\aa=(a _{ ij}) _{ 1\le i,j\le n}$. 
\\
When the fields $\aa, \BB$ and $V$ are regular enough   and not growing too fast, it is a  consequence of Itô formula that a solution to this equation is given by the Feynman-Kac formula
\begin{align}\label{eq-74}
g(t,x)=E_R \left[\left.\exp \left( \int _{ [t,T]} V(s,X_s)\, ds\right)\ g_T(X_T) \,\right|\, X_t=x \right] , 
\qquad x\in\ZZ, 0\le t\le T,
\end{align}
where $X$ is the canonical process, $R$ is the law of  a Markov process with generator $ \mathsf{A},$ and we denote by $E_R$ the expectation with respect to the measure $R$ and $E_R(\sbt\mid \sbt)$ the corresponding conditional expectation. This formula is named after R. Feynman and M. Kac for their contributions \cite{Fey48, Kac49,Kac51} in the late 1940's.  Defining
\begin{align}\label{eq-FKsg}
S_t^r u(x):=E_R \left[\left.\exp \left( \int _{ [r,t]} V(s,X_s)\, ds\right)\ u(X_t) \,\right|\, X_r=x \right] , \qquad 0\le r\le t\le T,
\end{align}
for any function $u$ on $\ZZ$ such that this expression is meaningful, 
 we see that $g_t=S^t_T g_T,$
and the collection of  linear operators $(S_t^r) _{0\le r\le t\le T }$ is the Feynman-Kac semigroup, where $g_t(x):=g(t,x).$
The stationary version of equation \eqref{eq-FK} 
\begin{align*}
( \BB\scal \nabla + \Delta _{ \aa}/2	+V)g
	=0,
\end{align*}
when  $\BB, \aa, V$ and $g$ do not depend on $t$,
is the stationary Schrödinger equation.
\\
The logarithmic  transformation
\begin{align}\label{eq-lt}
\psi:=\log g
\end{align} 
links  \eqref{eq-FK} to the Hamilton-Jacobi-Bellman equation 
\begin{align}\label{eq-HJB}
\left\{
\begin{array}{ll}
 e ^{ - \psi} (\partial_t +\mathsf{A})e ^{  \psi}+V
	=0,
	&\quad  0\le t<T,\\
 \psi(T,\sbt)= \psi_T, &\quad t=T,
\end{array}
\right.
\tag{HJB}
\end{align}
(formally, divide \eqref{eq-FK} by $g$ and replace $g$ by $ e ^{ \psi}$, provided that $g>0$). It   was a keystone of  Schrödinger's original derivation of his eponym equation\footnote{Beware, with our notation the role of the wave function $\Psi$ is played by $g$, not $\psi,$ see \eqref{eq-Sch} below.} because it permits to primarily work with some nonlinear Hamilton-Jacobi-Bellman (HJB)  equation  which is well-suited to carry both features of particle mechanics  and wave evolution, and then to transform it  into   Schrödinger's linear equation. It is also of importance in the theory of controlled Markov processes, see \cite[Ch.\,6]{FS06}.

Typical results about  classical  -- i.e.\ $\mathcal{C} ^{ 1,2}$ --  solutions of \eqref{eq-HJB} require that $\aa$ is uniformly  positive definite and that  $\aa, \BB, V$ and $g_T$ are  $ \mathcal{C} ^{ 1,2}_b,$ see  \cite{Kry87,evans-PDE}. We also know that when $\aa, \BB, V$ and $g_T$ are continuous, but $\aa$ might not be uniformly positive definite,  and the solution $g$ of \eqref{eq-FK} is also continuous (the Feller property of $ \mathsf{A}$ implies this continuity in several cases), then   
\(\psi:=\log g
\)   where $g$ is given by the Feynman-Kac formula \eqref{eq-74} is the viscosity solution of \eqref{eq-HJB}, see \cite[thm.\,II.5.1]{FS06}.

On the other hand, Kac proved in  \cite{Kac49}   (in  one dimension) that if $V$ is an  upper  and locally bounded measurable function and $ \mathsf{A}=  \partial_x^2$ is the generator of the Brownian motion, then $g$ given by \eqref{eq-74} solves  \eqref{eq-FK} in some weak sense. It was discovered  
later with $\mathsf{A}=   \Delta$ the generator of the Brownian motion in $\ZZ$ that when $V$ belongs the Kato class   (a set of lowly regular measurable functions which might not be locally bounded but  with some integrability properties), see Definition \ref{def-kato},  that the Feynman-Kac operator $S_t$ defined at \eqref{eq-FKsg} is a continuous operator  from $L^p$ to $ L^p$ with $1\le p\le \infty$ and that $g$ given by \eqref{eq-74} is continuous, see \cite[Ch.\,3]{CZ95}.

The ``Feynman-Kac transform''  of $R$ which is  the path  measure defined by
 \begin{align}\label{eq-77}
 P:= f_0(X_0)\exp \left( \Iii V(t,X_t)\,dt\right) g_T(X_T)\ R
 \end{align}
where  $ f_0:\ZZ\to [0, \infty)$  is another nonnegative function, is a generalization of Doob's $h$-transform  of $R$  \cite{Doob57,Doob84}, which is recovered by choosing $f_0=1, V=0$ and taking $g_T=h.$ When the solution $g$ of \eqref{eq-FK} is $ \mathcal{C} ^{ 1,2},$ with standard stochastic calculus arguments one proves   \cite{RY99, FS06}  that 
$P$  is the law of  a Markov diffusion process with the same matrix field $\aa$ as $R$, and  drift field
\begin{align}\label{eq-gir-b}
\BB^P=\BB+ \aa \nabla \psi.
\end{align}
The path measure $P$   is the solution to the Schrödinger problem \eqref{eq-SP} below,  a topic  also called \emph{entropic optimal  transport} which is tightly related to optimal transport, and   is currently an active field of research.  More will be said about entropic optimal transport in a  moment  in this introductory section.

\subsection*{Main results of the article}

In  this article, we prove with   probabilistic techniques that $g$ given by  the Feynman-Kac formula \eqref{eq-74} belongs to the domain of the \emph{extended} generator of the Markov measure $R$ and satisfies the trajectorial equation \eqref{eq-78} below, extending equation \eqref{eq-FK} in a probabilistic way. The diffusion matrix field $\aa$ is supposed to be  regular (typically $\aa= \sigma \sigma^*$ with $ \sigma$ locally Lipschitz) and  invertible, but   the coefficients $\BB, V $ and the datum $g_T$ are \emph{neither  assumed to be regular, nor locally bounded}. 

The notion of 
extended generator of a Markov measures  is directly connected to the  notion of semimartingale which   plays a central role in this article. 
The main hypothesis of this article  is that the \emph{relative entropy} of $P$ with respect to the reference path measure $R$ is finite, i.e.
\begin{align}\label{eq-entropy}
H(P|R):=E_P\log \Big( \frac{dP}{dR}\Big) < \infty.
\end{align}
We prove in addition that $g$ admits some  generalized spatial derivative $\gradt g$ and extend \eqref{eq-gir-b}:  
$P$  is the law of  a Markov diffusion process with the same matrix field $\aa$ as $R$, and  drift field
\begin{align}\label{eq-gir}
\BB^P=\BB+ \aa \gradt \psi,
\end{align}
with $\gradt\psi=\gradt(g)/g$.
The rigorous statement of  this formula requires some "almost everywhere" cautions, see  Theorem \ref{res-13} for the exact result.  
The main interest of this result is that it holds even when  not much is known a priori about the regularity of $ g$. For instance,   if $V$ is in the Kato class (a natural assumption in theoretical physics),  $g$ might be continuous but one does not know that it is differentiable in general.

The rule of the game in this article is to prohibit   regularity hypotheses on $V$ and $g_T$ stronger than the finite entropy condition \eqref{eq-entropy}. This is suggested by our main motivation which is the   entropic optimal transport. But  it also appears that this finite entropy assumption is very  efficient to derive  low regularity results.

\subsubsection*{First result}

Theorem \ref{res-31} is our first main result. Its  approximate statement is as follows. \\
\emph{
Suppose that
\begin{enumerate}[(i)]
 \item
$\aa$ and $\BB$ are such that $R$ is the ``unique'' solution to the martingale problem $\MP(\aa,\BB)$  (see  definition \eqref{eq-85} and  assumption \eqref{eq-H17a}),
\item
$\aa$ is locally bounded and invertible, 
\item
 $\aa, \BB, V $ and $g_T$ are such that $H(P|R)< \infty,$ where $P$ is defined at \eqref{eq-77}.
\end{enumerate}
Then, 
 $g$  defined by \eqref{eq-74} belongs to the  domain $\dom\LLL RP$  of the extended generator $\LLL RP$ of $R$ localized by $P$  (see Definition \ref{defh-05b}) and }
\begin{align}\label{eq-78}
[(\LLL RP  +V)g] (t,X_t)=0,\quad dtdP\ae
\end{align}

\begin{remarks}\ \begin{enumerate}[(a)]
\item
The tricky part of this result is: $g\in\dom\LLL RP.$  
\item
The hypothesis:  $H(P|R)< \infty,$  is an integrability assumption on the data $\aa, \BB, g_T$ and $V$ which requires  almost no regularity from $V$ and $g_T.$
\item
The extended generator $\LLL RP$ is localized by $P$, see Definition \ref{defh-05b},  and  \eqref{eq-78} is only valid $dtdP$-almost everywhere, rather than  $dtdR\ae$ a priori. It is partly because of this self-reference to the observed path measure $P$,  that it is possible to get rid of some regularity and growth restrictions on $g_T$ and $V.$
\item
We require that the growths of  $\BB$ and $\aa$   are  such that the reference path measure $R$ exists and their regularities are sufficient  for $R$ to be  the ``unique'' solution to its martingale problem - typically Lipschitz regularity. But some entropic argument allows us to depart from the regularity of  $\BB$, see Section \ref{sec-GC}. 
\item
On the other hand,  the additional hypothesis that $\aa$ is invertible is important for our approach to work: It is there to ensure the Brownian martingale representation theorem  which implies that the domain of the extended generators of $R$ and $P$  are algebras which are  stable by $C^2$ transformations, see Lemma \ref{res-23}.
\end{enumerate}\end{remarks}

\subsubsection*{Second result}
Theorem \ref{res-13} which extends \eqref{eq-HJB} and \eqref{eq-gir-b} is our second main result.  Its approximate statement is as follows.
\\
\emph{
Under the same assumptions as before, the function $\psi=\log g$ where $g$ is given by the Feynman-Kac formula \eqref{eq-74} solves the following extended HJB equation
\begin{align}\label{eq-Lsol}
\big(\LLL RP \psi+| \Gt RP \psi|^2_\aa /2 +V\big)(t,X_t)=0,
\qquad dtdP\ae
\end{align}
with $ \psi_T=\log g_T,$  $P_T\ae$ 
\\
The Feynman-Kac measure $P$ solves the martingale problem 
$\MP(\aa,\BB^P)$ where
}
\begin{align}\label{eq-z09-b}
\BB^P=\BB+\aa \Gt RP \psi.
\end{align}
Martingale problems are defined at Definition \ref{defh-02} -- we slightly depart  from the standard definition, see Remark \ref{rem-MP}-(c). The existence of $\Gt RP \psi$ and its definition are stated  at Proposition \ref{res-14}  which is an extended Itô formula. When $ \psi$ is differentiable, we have: $\Gt RP \psi=\nabla \psi.$

In fact, our first result (Theorem \ref{res-31}) about equation \eqref{eq-FK} is a corollary of the above extended HJB  equation.

\subsubsection*{Third result}

A sufficient condition on $\aa,\BB, V, f_0$ and $g_T$ is stated at Theorem \ref{res-32} for $H(P|R)$ to be finite. As it is rather technical, we do not propose in this introduction  an approximate description of this set of assumptions.  However, let us say that  it is large enough to include potentials $V$ in the Kato class. 
Therefore we are in position to apply our previous results about the FK and HJB equations and the FK-transform  $P.$ In particular, 
\\
\emph{the representation \eqref{eq-gir-b} of the drift of $P$ which is inaccessible when $V$ is a generic element of the Kato class ($g$ is only known to be continuous in  most favorable situations) admits   the extension  \eqref{eq-z09-b}.}

\subsubsection*{Fourth set of results}

Our last results are of  different nature than previous ones. They state that \emph{stochastic derivatives and extended generators are essentially the same.}  This already well-known assertion (since Nelson's monograph \cite{Nel67} and its use by Föllmer in a series of works \cite{Foe85b,Foe86,Foe85})  is the content of Propositions \ref{resh-01} and \ref{resh-03}. The proofs of Proposition \ref{resh-01}, and also Theorem \ref{res-13}, deeply rely on the convolution Lemma \ref{resh-21}. This technical lemma permits to fill some gaps in the already published literature on the subject. We had to put on a solid ground  the use of stochastic derivatives to compute extended generators because it is central in the approach of  the present  paper. 
\\
In addition, the convolution Lemma \ref{resh-21} is used in an essential manner when deriving time reversal formulas in the recent article
\cite{CCGL20}.

\subsection*{Some comments and remaining questions}

 It is important to stress that this paper primarily states that $g$ satisfies the identity \eqref{eq-78} and  $ \psi$  satisfies the identity \eqref{eq-Lsol}. These statements should not be  interpreted as $ g$ solves equation \eqref{eq-78} and   $  \psi$ solves equation \eqref{eq-Lsol}. We only look at $g$ defined by the Feynman-Kac formula, and $ \psi:=\log g,$ and do not study any new notion of solution to the \emph{equations} (FK) and (HJB): the existence of such a solution would automatically be given by the Feynman-Kac formula. On the other hand, the uniqueness problem is not even clear to state: unique in which regularity class? while we do not even know the minimal regularity of $g$ under our general hypotheses.

We present some problems  which are not treated in this article,  and make a couple of comments.
\begin{enumerate}[(i)]
\item
As already noticed, if the Feynman-Kac formula defines a  function $g$ which  is continuous,  
then  $ \psi$ is the viscosity solution of \eqref{eq-HJB}.  We do not know whether it  remains a viscosity solutions when it is discontinuous, under the only hypothesis that $H(P|R)$ is finite.  

\item
Because \eqref{eq-78} and \eqref{eq-Lsol} are trajectorial statements, they are more precise than any  solution to a PDE. 
The specificity of a pathwise representation in a probabilistic context is twofold: one can play with \emph{stopping times} or with \emph{couplings}, and sometimes both.   Looking at
\begin{align*}
g(t,x)=E_R \left[\left.\1 _{ \{t\le \tau\}}\exp \left( \int _{ [t, \tau]} V(s,X_s)\, ds\right)\  g _{ \mathrm{\,final}}(X_ \tau) \,\right|\, X_t=x \right] , 
\end{align*}
where $ \tau$ is a stopping time, is tempting. Similarly, one may ask whether couplings are of some use when looking at comparison principles or functional inequalities.
\item

On the other hand, this advantage is balanced by some drawbacks. In particular, the powerful stability properties  of viscosity solutions along  convergence schemes   might not be recovered via  \eqref{eq-Lsol} . Typically, in case of a vanishing viscosity convergence, because the supports of diffusion path measures with different diffusion matrices are disjoint, a trajectorial solution does not permit us to use pointwise convergence.

\item
Statements \eqref{eq-78} and \eqref{eq-Lsol}   rely on the existence of some Markov path measure $R$. This  is  restrictive in comparison to the   general definition of  viscosity solution  which only requires the existence of some semigroup   obeying the maximum principle  \cite[Ch.\,2]{FS06}.

\end{enumerate}

\subsection*{Standard approaches for computing the generator of $P$}

There are three main ways to look at the dynamics of $P.$ They rely on (i) Markov semigroups, (ii) Dirichlet forms and (iii) semimartingales.

\subsubsection*{Markov semigroups}

Let $P\in\PO$ be a Markov measure (see Definition \ref{defh-09} below) and $(T_{s,t}^P) _{ 0\le s \le t\le T}$ be its semigroup on some Banach function space $(U,\|\cdot\|_U).$ For instance $U$ may be the space of all  bounded Borel measurable functions   equipped with the topology of uniform convergence. Its infinitesimal generator is $A^P=(A^P_t)_{t\in\ii}$ with
\begin{equation}\label{eq-41}
	A^P_tu(x):=\|\cdot\|_U\hbox{-}\Limh \frac 1h E_P[\left.u(X_{t+h})-u(X_t) \,\right|\, X_t=x],
	\quad u\in\dom A ^P
\end{equation}
where the domain $\dom A^P$ of $A^P$ is precisely the set of all functions $u\in U$ such that the above strong limit exists for all $t\in [0,1)$ and $x\in\ZZ.$
One can prove rather easily (see \cite[Ch.\ VIII,\S 3]{RY99} for instance, in the diffusion case) that when $V$ is zero and $g$ is positive and regular enough,  the generator $A^P$ of the Markov semigroup associated with $P$ is given for regular enough functions $u$ on $\ZZ,$ by
\begin{equation}\label{eq-40}
A^P_t u(x)=A^R u(x)+\frac{\Gamma(g,u)}{g}(t,x),\quad (t,x)\in\ii\times\ZZ
\end{equation}
where $\Gamma$ is the carré du champ operator, defined for all functions $u,v$ such that $u,v$ and the product $uv$ belong to the domain $\dom A^R$ of $A^R,$ by
$$
	\Gamma(u,v)=A^R(uv)-u A^Rv-vA^Ru.
$$
For Eq.\,\eqref{eq-40} to be meaningful, it is necessary that for all $t\in\ii,$ $g_t$ and the product $g_tu$ belong to $\dom A^R.$ 
But with a non-regular potential $V$, $g$ might be non-regular as well. There is no reason why $g_t$ and $g_tu$ are in $\dom A^R$ in general.
Clearly, one must drop the semigroup approach and work with  Dirichlet forms or semimartingales. 

\subsubsection*{Dirichlet forms}

The Dirichlet form theory is natural for constructing irregular processes and has been employed in similar contexts, see \cite{Alb00}.    It is made-to-measure for reversible processes. Though there exists a theory for non-symmetric Dirichlet forms \cite{MR92,O13,St99}, it is not fully efficient for our purpose.

\subsubsection*{Semimartingales}

Working with semimartingales means that instead of the infinitesimal semigroup generators $A^R$ and $A^P$, we consider \emph{extended generators} in the sense of the Strasbourg school \cite{DM4}, see Definition \ref{defh-05} below. This natural idea has already been implemented by P.-A. Meyer and W.A. Zheng \cite{MZ84, MZ85} in the context of stochastic mechanics and also by  P. Cattiaux and the author in \cite{CL94,CL96} for solving  related entropy minimization problems. But we still had to face the remaining problem of giving some sense to $\Gamma(g_t,u)$ in \eqref{eq-40}. Consequently, restrictive assumptions were imposed: reversibility in \cite{MZ85},  and  in \cite{CL96}: the standard hypothesis that the domains of the extended generators of $R$ and $P$ contain  ``large'' sub-algebras. In practice this  requirement is uneasy to satisfy, except for standard regular processes. It is all right for the reference measure $R$, but typically when $V$ blows up, $P$ is  singular and this large sub-algebra assumption does not seem to be accessible with standard arguments. Moreover, there is no known criterion  for this property to be inherited  from $R$ by $P$ when $P\ll R.$ In contrast, extended generators and considerations about the carré du champ  allow us to overcome this obstacle in the present article,  see  Lemma \ref{res-23} which is based on Lemma \ref{res-22}.

\subsection*{Some motivations}

Let us  present our motivations for proving the main results of this paper. Although very much related to each other, we have two distinct problems in mind: (i)  the first one is  coming from theoretical physics in connection with the interpretation of dissipative evolutions and quantum mechanics, (ii) the second one is related to the entropic optimal transport problem, also called the Schrödinger problem.

\begin{enumerate}[(i)]
\item
\emph{A remarkable analogy between dissipative and quantum evolutions.}\ 
In the early thirties, Schrödinger \cite{Sch31,Sch32}  addressed the entropy minimization problem
\begin{align}\label{eq-SP}
\inf\{H(Q| R^V); Q \textrm{ path measure such that } Q_0= \mu_0, Q_T= \mu_T \}
\end{align}
where the  couple $(Q_0,Q_T)$ of initial and final marginals of $Q$ is prescribed to be equal to some fixed  $( \mu_0, \mu_T),$ and the reference measure is $R^V:= \exp \left( \Iii V(t,X_t)\,dt\right)\, R.$    Its solution $P$ is called a \emph{Schrödinger bridge} whose  time-marginal flow $(P_t) _{ 0\le t\le T}$ is called an \emph{entropic interpolation} between $ \mu_0$ and $ \mu_T.$

In the specific case where $R$ is the law of a reversible Brownian and $V=0$, Schrödinger essentially showed  that formula \eqref{eq-77} gives the general shape of the Schrödinger bridge and that the entropic interpolation  admits the Radon-Nikodym derivative 
\begin{align}\label{eq-Born-SP}
\frac{dP_t}{dx}
	=f_t(x)g_t(x) 
\end{align}
where $g$ solves \eqref{eq-FK} and 
$f_t(x)=E_R \left[\left.f_0(X_0)\exp \left( \int _{ [0,t]} V(s,X_s)\, ds\right) \,\right|\, X_t=x \right]$ 
solves a forward-time analogue of \eqref{eq-FK}
\begin{align}\label{eq-FK-f}
\left\{
\begin{array}{ll}
(-\partial_t + \widetilde{\mathsf{A}}	+V)f
	=0,
	&\quad  0< t\le T,\\
 f(0,\sbt)=f_0, &\quad t=0,
\end{array}
\right.
\end{align}
which again is a Feynman-Kac equation 
with $ \widetilde{\mathsf{A}}= \tilde\BB\scal \nabla + \Delta _{ \aa}/2$
where the vector field  $\tilde\BB$ is the drift of the time reversal of $R.$
\\
He    noticed a striking analogy (in his own words) between the solution of the thermodynamical problem \eqref{eq-SP} described by the product formula \eqref{eq-Born-SP}  and  Born's formula
\begin{align*}
 \frac{d \mu_t}{dx}=| \Psi_t|^2(x)= \Psi_t(x) \overline\Psi_t(x)
\end{align*}
where $ \mu_t$ is the probability of presence of some quantum system at time $t$ and   $\Psi$ is the wave  function describing its evolution, solution of 
 his eponym equation:
\begin{align}\label{eq-Sch}
i\hbar \partial_t \Psi= \Big(-\frac{\hbar^2}{2m} \Delta +V\Big)\,\Psi
\quad
\iff
\quad
-i\hbar \partial_t \overline\Psi= \Big(-\frac{\hbar^2}{2m} \Delta +V\Big)\,\overline\Psi
\end{align}
and $\overline\Psi$ is the complex conjugate of $\Psi.$
To see this analogy, remark that with $V$ instead of $-V$, taking $ \tau=it$, we have $-i \partial_t= \partial_\tau$. It follows that    equation \eqref{eq-FK-f}   is an analogue of the Schrödinger equation for $\Psi$, while   the time-reversed equation   \eqref{eq-FK} is an  analogue of the Schrödinger equation for $\overline\Psi$. 
\\
Not only does the Feynman-Kac transform $P$ defined at \eqref{eq-77} provide us with interesting analogies with quantum mechanics, but  its  family of bridges $P(\sbt\mid X_r=x, X_t=y)$   is the  classical thermodynamics analogue  of the propagator  appearing in Feynman's approach to quantum mechanics \cite{FH65}: the ill-defined Feynman integral is replaced by a  stochastic integral. 
\\
A couple of potential applications of the present article are \begin{enumerate}
\item[-\ ]
the EQM (Euclidean quantum mechanics) program launched a long time ago by Zambrini \cite{Zam86,ChZ03, Zam12}, and 
\item[-\ ]
the closely related theory of Bohmian mechanics \cite{Bo52,DT09}. 
\end{enumerate}
EQM describes the  thermal evolution of particle systems, tracking as much as possible the analogies with quantum mechanics, while Bohmian mechanics describes the \emph{pathwise}  evolution of quantum systems  with the same prediction as standard quantum theory.  In these settings, $V$ is a typical scalar potential generating a  force field  such as the Coulomb potential which is \emph{not locally bounded.}  A natural set of such potentials is the Kato class. This is the reason why we prove at Section 6 that scalar potentials in the Kato class satisfy our main finite entropy hypothesis \eqref{eq-entropy} under standard additional assumptions on $R$.

\item
\emph{Entropic optimal transport.}\ 
Entropic optimal transport (EOT) is another name for the Schrödinger problem \eqref{eq-SP}  emphasizing its connection with (standard) optimal transport (OT). It is an active research topic. 

EOT is often considered  taking $R$ to be the law of a diffusion process, and seen as an entropic regularization of quadratic OT. In this context, it is known that for "generic" marginal constraints $ \mu_0$ and $ \mu_T$, the solution of \eqref{eq-SP} is \eqref{eq-77}, i.e.\ the Schrödinger bridge is the FK-transform $P=f_0(X_0)g_T(X_T)\,R^V$ for some couple of functions $(f_0,g_T)$ solving the \emph{Schrödinger system}
\begin{align}\label{eq-SS}
\left\{
\begin{array}{lcl}
f_0g_0&=&d \mu_0/dR_0\\
f_Tg_T&=&d \mu_T/dR_T
\end{array}
\right.
\end{align}
where  $g_0$ is defined by the Feynman-Kac formula \eqref{eq-74} and $f_T$ is defined in a time-symmetric way, that is for all $ x\in\ZZ$ and  $0\le t\le T,$
\begin{align*}
g_0(x)&:=E_R \left[\left.\exp \left( \int _{ [0,T]} V(t,X_t)\, dt\right)\ g_T(X_T) \,\right|\, X_0=x \right],\\
f_T(x)&:=E_R \left[f_0(X_0)\ \left.\exp \left( \int _{ [0,T]} V(t,X_t)\, dt\right) \,\right|\, X_0=x \right].
\end{align*}
The Schrödinger system \eqref{eq-SS}  follows from: $dP_t/dR_t=f_tg_t,$ $0\le t\le T,$ see Proposition \ref{res-16} below, and extends \eqref{eq-Born-SP}.
The logarithms $ \varphi_0:=\log f_0$ and $ \psi_T:=\log g_T$  are called the \emph{Schrödinger potentials}, in analogy with  the Kantorovich potentials of OT.  

The main application to EOT of the results of the present article is the  identity \eqref{eq-z09-b}: $\BB^P=\BB+\aa \Gt RP \psi,$ which characterizes the evolution of the Schrödinger bridge $P$.  In this formula $ \psi$ satisfies the extended HJB identity  \eqref{eq-Lsol}  with the boundary condition $ \psi_T$ at time $t=T$ which requires to solve the Schrödinger system \eqref{eq-SS}.

This improves already known results in the literature, in particular because $\BB$ is allowed to be a singular drift.  Indeed,  our approach does not require that $\aa$, $\BB$ and $V$ meet the standard hypotheses of regularity theorems for the Feynman-Kac equation \eqref{eq-FK} which  are obtained using  PDE techniques \cite{evans-PDE} or stochastic flows \cite{Kun97}.
\end{enumerate}

\subsection*{Future works}
It is the purpose of EQM to transpose well-established results in stochastic analysis of variational processes to standard quantum mechanics, and the other way round. From the start, the EQM program sees the Feynman-Kac equation (FK) as the stochastic deformation of Schrödinger's equation. 
The associated entropic interpolation:   the marginal flow $(P_t)$ of the FK path measure $P$, is governed by some Newton equation in the Otto-Wasserstein space. This  was proved by  Conforti  in \cite{Co18},  assuming that $f_t$ and $g_t$ are regular enough and not growing too fast. Our aim in a  work in progress is to relax these assumptions and to extend this  result under the finite entropy condition \eqref{eq-entropy}.

On the other hand,  it is worth  expressing Bohmian mechanics as the evolution of a fluid also living in the Otto-Wasserstein space. This is suggested by von Renesse's work \cite{vR11} where Schrödinger's equation is  seen at a heuristic level  as  a Newton equation in the Otto-Wasserstein space, again.  In this perspective,   the Newton equation  associated to \eqref{eq-FK} is the classical rigorous analogue of Newton's equation in Bohm's view of quantum dynamics.   This is also a  work in progress.

The time-symmetry of formula \eqref{eq-77} suggests that the forward-in-time equation \eqref{eq-FK-f} is as important as its backward-in-time analogue \eqref{eq-FK}. This is crucial in many aspects of EOT and EQM. And indeed this time-symmetry is an important ingredient of these works in progress. It relies on the recent article \cite{CCGL20} about time reversal of diffusion processes under a finite entropy condition. 

\subsection*{Further developments}
Let us raise some remaining problems.

\begin{itemize}
\item
\emph{Jumps}.\ 
Replacing the diffusion measure $R$ by the law of a time-continuous  Markov process with jumps would lead us to a similar trajectorial  identity associated a nonlocal  PDE:  $( \partial_t+ \mathsf{A}+V)g=0,$ with $ \mathsf{A}$ a  Markov  generator expressed with a jump kernel, and to the associated nonlocal HJB equation. However, the  conceptual limitation one might face is the martingale representation theorem (MRT), which is crucial  to prove the extended Itô formula (Lemma \ref{res-23}, Proposition \ref{res-14}). Indeed, to our knowledge MRT results are only available for specific classes of jump processes, see \cite{Z09} for instance.

\item
\emph{Manifold}.\  What happens when replacing $\ZZ$ by a Riemannian manifold? Unlike the jump setting, there exists a well established Brownian MRT on a complete connected Riemannian manifold, see \cite{F94,Hsu97}. It remains to establish a Girsanov theory under a finite entropy condition in the spirit of  \cite{Leo11a}, but this must be straightforward.

\item
\emph{Viscosity solutions}.\  The links between   the trajectorial representation of the solution to the HJB equation encountered at Theorem \ref{res-13}, see \eqref{eq-Lsol}, and its viscosity solution remain to be clarified.

\item
\emph{Small noise limit}.\ 
Suppose that the Feynman-Kac path measure $P ^{ \epsilon}$, built on some
$f_0 ^{ \epsilon},g_T ^{ \epsilon}, \BB ^{ \epsilon},$ $ V ^{ \epsilon}=V/ \epsilon$ and $\aa ^{ \epsilon}:= \epsilon\aa,$ satisfies the hypotheses of Theorem \ref{res-13} and  converges weakly to some $P^0$ in $\PO$ as $ \epsilon$ tends to zero.  This  implies that the entropic  interpolations converge to the displacement interpolations because taking the marginal flow is a continuous mapping.  But one can ask if there is some \emph{almost sure convergence}  of the sample paths.  This is an open problem whose solution should rely on large deviation results to obtain exponentially small probability of deviations, ensuring some argument based on Borel-Cantelli lemma. In this respect, see the recent article \cite{BGN21a}  which establishes a  large deviation principle in the setting of the {static} Schrödinger problem.  Its dynamical counterpart remains to be proved.
\end{itemize}

\subsection*{Literature}

Schrödinger \cite{Sch31,Sch32} only considered the case where $R$ is the law of a reversible Brownian motion and $V=0$. His arguments for deriving \eqref{eq-Born-SP} are  profound but not rigorous (at this time the axioms of probability theory were   unsettled and  the Wiener process was unknown), but the impressive strength of the physicist's arguments is sufficient to convince the reader. The extension to a non-zero potential $V$ under the hypothesis that  the solutions $f$ and $g$ of the Feynman-Kac equations in both directions of time  are regular enough was performed by Zambrini \cite{Zam86}.
The actual writing of the {Schrödinger problem} \eqref{eq-SP} in terms of entropy, as well as its 
formulation as a large deviation problem for the empirical measure of a system of particles is due to Föllmer \cite{Foe85}.  More about this active field of research, in particular its tight connection with optimal transport, can be found in the survey paper \cite{Leo12e}.

\subsection*{Outline of the present approach}

A key feature of our approach is  the logarithmic transformation \eqref{eq-lt}
 because it enables us to  take  advantage of a  connection between \eqref{eq-HJB} and Girsanov's theory.
More precisely,  elementary stochastic calculus gives the following expression of the Radon-Nikodym derivative of $P$ defined at \eqref{eq-77} with respect to $R$:
\begin{equation}\label{eq-76}
\begin{split}
\frac{dP}{dR}&:=f_0(X_0)\exp \left( \Iii V(t,X_t)\,dt\right) g_T(X_T)\\
	&\,=f_0(X_0) g_0(X_0)\, \exp \left( \psi_T(X_T)- \psi_0(X_0)-\Iii [e ^{ - \psi} (\partial_t +\mathsf{A})e ^{  \psi}](t,X_t)\, dt \right) ,
\end{split}
\end{equation}
 where 
\begin{align}\label{eq-75}
 \psi(t,x)=\log g(t,x)=\log\, E_R \left[\left.\exp \left( \int _{ [t,T]} V(s,X_s)\, ds\right)\ g_T(X_T) \,\right|\, X_t=x \right] , 
\end{align}
provided that $\aa, \BB, g_T$ and $V$ satisfy some growth conditions   and that \emph{$ \psi$ is regular enough} to apply Itô formula and give sense to $ (\partial_t +\mathsf{A}) e ^\psi$. At first sight the identity \eqref{eq-76} is reminiscent to  \eqref{eq-HJB}, and indeed it establishes a strong link between \eqref{eq-FK} and \eqref{eq-HJB}. The main problem we have  to face is to develop  this simple idea,   when one does not know much about the a priori regularity of $ \psi$. In particular, $(\partial_t +\mathsf{A})e ^{  \psi}$ is a priori undefined.

 A good thing to do   is  to compare the above expression of $dP/dR$ with the one obtained by Girsanov's theory. Indeed, this  provides us with valuable informations on $ \psi$, and therefore  on the solution of \eqref{eq-FK}. This is possible at the price of working with \emph{extended generators} instead of standard generators of Markov semi-groups, allowing us to  extend Itô formula to the domain of the extended generator under the important requirement that \eqref{eq-entropy}: $H(P|R)<  \infty,$ holds.

This entropy estimate is in fact a finite energy condition which carries some control of the generalized derivative $ \gradt \psi$ of  $ \psi$  which takes part of an extended Itô formula (this is the place where it is needed that $\aa$ is invertible to make sure that the martingale representation theorem is valid). On the other hand, Girsanov's theory tells us that $\gradt \psi$ is precisely the additional drift which ``translates'' $R$ to $P$, see \eqref{eq-gir}.  For a better understanding of the key point of the  present strategy which takes advantage of Girsanov's formula to allow us to get rid of an a priori regularity of the solution $g$, see Remark \ref{rem-05}.

\subsection*{Outline of the paper}

 The  specific features of the present  extension of the Feynman-Kac equation \eqref{eq-FK}  are introduced at Section \ref{sec-EDM} which contains both standard results about finite entropy diffusion path measures, and a bit of new material designed for our purpose (especially the extended Itô formula at Proposition \ref{res-14}). This material  is based on extended generators, a standard notion which is revisited at Section \ref{sec-SD}, using Nelson stochastic derivatives, very much in the spirit of the seminal paper \cite{Foe86} by H. Föllmer (see also \cite{Foe85} and more recently \cite{CCGL20}), and  generalized  at Section \ref{sec-sd-ext}.  Using the preliminary material established in the first four sections, Section \ref{sec-HJB} is dedicated to the proofs of our main results. At Section \ref{sec-GC}, sufficient conditions are established  on the coefficients of equation \eqref{eq-FK}   for $H(P|R)$ to be finite.

\section{Stochastic derivatives. Main results}
\label{sec-SD}

After reviewing basic notions of  semimartingale theory in a general setting (càdlàg paths in a Polish space), we prove at Propositions \ref{resh-01} and  \ref{resh-03} that under some integrability condition,  the extended generator and the stochastic derivative of a Markov measure coincide. This is Nelson's way of looking at diffusion generators  \cite{Nel67}.  

The aim of the present section is to  provide rigorous  proofs  of these general results.  To our knowledge, although these notions, in connection with the notion  of martingale problem, were introduced in the late sixties \cite{Nel67,K69,SV69a}, such detailed proofs  do not appear in the literature. However,  the guideline they provide  and the recognition of the relevance of these notions for our purpose are fully credited to Föllmer \cite{Foe86}.  Our main technical tool is  the convolution  Lemma \ref{resh-21}.

\subsection*{Notation and setting}
The set of all probability measures on a measurable set $A$ is denoted by $ \mathrm{P}(A)$ and the set of all nonnegative $ \sigma$-finite measures on $A$ is  $ \mathrm{M}(A).$
The push-forward of a measure $ \qq \in  \mathrm{M}(A)$ by the measurable map $f:A\to B$ is   $f\pf\qq (\sbt):=\qq (f\in \sbt)\in \mathrm{M}(B).$ 
\\
The state space $\XX$ is assumed to be Polish and is equipped with its Borel $ \sigma$-field.
The path space is the set $$\OO=D(\ii,\XX)$$ of all  càdlàg $\XX$-valued trajectories $\omega=(\omega_t)_{t\in\ii}\in\OO.$ It is equipped with the canonical $\sigma$-field:  $\sigma(X_t;t\in\ii)$ which is generated by the canonical process $X=(X_t)_{t\in\ii}$ defined for each $t\in\ii$ and $\omega\in\OO$ by $$X_t(\omega)=\omega_t\in\XX.$$
We denote $\Ob:=\ii\times\OO,$ $\XXb:=\iX$, and   for any $t\in\ii,$ 
$$
\Xb_t:=(t,X_t)\in\XXb,
$$
and  any function $u:\iX\to\RR$,   $$u(\Xb): (t, \omega)\in\Ob\mapsto u(t, \omega_t)\in\RR.
$$
 We call any positive measure  $Q\in\MO$ on $\OO$ a \emph{path measure}. For any $\mathcal{T}\subset\ii,$ we denote $X_\mathcal{T}=(X_t)_{t\in \mathcal{T}}$ and the push-forward measure $Q_\mathcal{T}=(X_\mathcal{T})\pf Q.$ In particular, for any
 $0\le r\le s\le T,$ $X_{[r,s]}=(X_t)_{r\le t\le s}$, $Q_{[r,s]}=(X_{[r,s]})\pf Q$, and $Q_t=(X_t)\pf Q\in\MX$ denotes the law of the position $X_t$ at time $t$. If $Q\in\PO$ is a probability measure, then $Q_t\in\PX$.
\\
For any $0\le t\le T$,   $\Qb  :=\Leb _{ [0,T]}\otimes Q $ is the product measure 
\[\Qb  (dtd \omega):=dtQ (d \omega),\qquad dtd \omega\subset\Ob
\]
and we denote 
\[\qb (dtd x):=dtQ_t (d x),\qquad dtd x\subset\XXb.\]
For any $Q\in\PO,$ we denote 
\begin{align*}
[Q]:=(Q_t; 0\le t\le T)\in \PX ^{ \ii}
\end{align*}
its time marginal flow. 
\\
For any random time $ \tau$, we denote $Y_t^\tau:=Y_{t\wedge \tau}$ and $\Xb_t^\tau:=(t\wedge \tau,X_{t\wedge \tau}).$

 \subsection*{Filtration} \emph{The forward filtration associated with $Q\in\MO$  is the $Q$-completion of the canonical filtration.}  It  fulfills the so-called ``usual hypotheses'': it is right continuous and contains the $Q$-null sets.  
 Under this hypothesis  it is known that any  $Q$-martingale  admits a càdlàg version, see for instance \cite{Lowther-cadlag}. We shall choose this version in all cases.


\subsection*{Basic notions}

We recall the definitions of Markov measure, extended generator and stochastic derivative.

 \begin{definition}[Markov measure] \label{defh-09}
 A path measure $Q$ such that  $Q_t$ is  $ \sigma$-finite for all $t$ is called a conditionable path measure.
  A path measure $Q\in\MO$ is said to be Markov if it is conditionable and for any $0\le t\le T,$ $Q(X _{ [t,T]}\in\sbt\mid X _{ [0,t]})=Q(X _{ [t,T]}\in\sbt\mid X _t).$ 
 \end{definition}

  The reason for requiring $Q$ to be conditionable is that  it allows for defining  the conditional expectations $E_Q(\sbt\mid X_ \mathcal{T})$ for any $ \mathcal{T}\subset \ii$ even in the case where $Q$ is an unbounded measure, see \cite[Def.\,1.10]{Leo12b}.  Remark that the definition of a  $Q$-martingale remains unchanged when $Q$ is unbounded because  wether $Q$ is bounded or not any conditioning $Q(\sbt\mid \mathcal{C})$ of $Q$ is a  \emph{probability} measure.

Let $Q$ be a path measure. Recall that a process
$M$ is called a \emph{local $Q$-martingale} if there exists a
sequence $\seq\tau k$ of $\ii\cup\{\infty\}$-valued stopping times
such that $\Lim k\tau_k=\infty,$ $Q\ae$ and for each $k\ge1,$ the
stopped process $M^{\tau_k}$ is a uniformly integrable
$Q$-martingale.
 A process $Y$ is called a \emph{special
$Q$-semimartingale} if $Y=B+M,$ $Q\ae$ where $B$ is a predictable
bounded variation process and $M$ is a local $Q$-martingale.

\begin{definition}[Nice semimartingale]\label{defh-01}
A process $Y$ is called a nice\footnote{This is a ``local'' definition in the sense that this notion  probably appears somewhere else with another name.} $Q$-semimartingale if $Y=B+M$ where $M$ is a local $Q$-martingale and the bounded variation process $B$
has \emph{absolutely continuous} sample paths $Q\ae$
\end{definition}

The notion of extended generator was introduced by H. Kunita \cite{K69} and extensively used by P.\,A.\,Meyer and his collaborators, see \cite{DM4}. Here is a variant of his definition.

\begin{definition}[Extended generator of a path measure]\label{defh-05}
Let $Q\in\MO$ be a conditionable path measure. A measurable function $u$ on $\XXb$ is said to be in
the domain of the extended generator of $Q$ if there exists an
adapted process   $\big(v(t, X _{ [0,t]});\ 0\le t\le T\big)$   such that
 $\Iii| v(t, X _{ [0,t]})|\,dt<\infty,$ $Q\ae$ and the process 
 $$M^u_t:=u(\Xb_t)-u(\Xb_0)-\II0t v(s, X _{ [0,s]})\,ds,\quad 0\le t\le T,$$ is a local $Q$-martingale. 
We denote
$$
   \LQ u(t, \omega):=v(t, \omega _{ [0,t]})
$$
and call $\LQ $ the extended generator of $Q.$ The domain
of the extended generator of $Q$ is denoted by $\dom
\LQ .$
\end{definition}

\begin{remark}[Special  case where $Q$ is Markov] \ 
 It is proved at  Corollary \ref{res-18} that when $Q$ is Markov,  $\LQ u$  only depends on the current position: 
$
\LQ u(t, \omega)=\LQ u(t, \omega_t).
$
It is also shown at Corollary \ref{res-19} that under some hypotheses
$$
 \LQ =\partial_t+ \mathcal{G}_t,
$$
 if  $ (\mathcal{G}_t) _{ 0\le t\le T}$ is the generator of the semigroup associated to the Markov measure $Q$.
\end{remark}

We go on with technical considerations.

\begin{remarks}\ \begin{enumerate}[(a)]
\item
For any measurable function $u$ on $\XXb$ in
$\dom\LQ $,  the process $u(\Xb)$ is a nice
$Q$-semimartingale.
\item
The adapted process $t\mapsto\II0t
v(s,X _{ [0,s]})\,ds$ is predictable since it is continuous. 
\item
$M^u$ admits a càdlàg $Q$-version as a local $Q$-martingale (we always choose this regular version). 
\item
The notation $v=\mathcal{L}u$ almost rightly suggests that $v$ is
a function of $u.$ Indeed, when $u$ is in $\dom\LQ ,$
the Doob-Meyer decomposition of the special semimartingale
$u(\Xb_t)$  into its predictable bounded variation part $\int
v_s\,ds$ and its local martingale part is unique. But one can
modify $v=\LQ u$ on a  zero-potential set without
breaking the martingale property.   As a consequence,  $u\mapsto
\LQ u$ is a single-valued linear operator  once  $\LQ u$ is identified with a class of functions modulo zero-potential sets. In the present context, a measurable  set $A\subset \iX$ is a zero-potential set if $\int_A dQ_tdt=0.$ 

\item
  This definition of the extended generator is suggested in the comment  XV-(21) of \cite{DM4}. It slightly differs from  the definition XV-(20) of \cite{DM4} which is expressed in terms of  resolvents.
\end{enumerate}
\end{remarks}

Extended generators are connected with   martingale problems which were introduced by Stroock and Varadhan \cite{SV69a,SV69b,SV79}.

\begin{definition}[Martingale problem]\label{defh-02}
Let $\mathcal{C}$ be a class of measurable real functions $u$ on $\XXb$
and for each $u\in\mathcal{C},$ let $\mathcal{L}u:\Ob\to\RR$  be some
adapted process. Take also a positive $ \sigma$-finite measure
$\mu_0\in\MX.$ 
\\
One says that the conditionable path measure $Q\in\MO$ is a solution to the
martingale problem 
$$
\MP(\mathcal{C},\mathcal{L};\mu_0)
$$ if
$Q_0=\mu_0\in\MX$ and for all $u\in\mathcal{C},$ 
$Q(\Iii |\mathcal{L}u(t,\omega _{ [0, t]})|\,dt=\infty)=0$
and
the process
$
 M^u_t:=   u(\Xb_t)-u(\Xb_0)-\II0t \mathcal{L}u(s,X _{ [0,s]})\,ds
$
is  a local $Q$-martingale.
\end{definition}

\begin{remarks}\ \begin{enumerate}[(a)] \label{rem-MP}
\item
As in Definition \ref{defh-05}, the local martingale $M^u$ admits a  càdlàg $Q$-version. 

\item
Playing with the definitions, it is clear that any path
measure $Q\in\MO$  is a
solution to $\MP(\mathcal{C},\mathcal{L}^Q;Q_0)$ where $\mathcal{L}^Q$ is the extended generator of $Q$ and
$\mathcal{C}$ is any nonempty subset of $\dom\mathcal{L}^Q.$

\item
In any standard definition of a martingale problem, it is assumed that for any $u\in \mathcal{C}$ and \emph{all} $ \omega\in\OO,$ we have: $\Iii |\mathcal{L}u(t,\omega _{ [0,t]})|\,dt<\infty$ (and not only $Q\ae$ as above). This will not be convenient for our purpose because when looking at $Q\in\PO$ such that the relative entropy $H(Q|R)$ with respect to some reference path measure $R$ is finite, the extended generator $ \mathcal{L}^ Q$ of $Q$ is only defined $\Qb\ae$, see \eqref{eq-H04} below for instance.
\end{enumerate}\end{remarks}

Our aim is to show at Proposition \ref{resh-01} that the extended generator can be computed by means of a stochastic derivative. 
Nelson's definition \cite{Nel67} of the stochastic derivative is the following.

\begin{definition}[Stochastic derivative]\label{defh-04}
Let $Q\in\MO$ be a conditionable path measure and $u$ be a measurable real function on $\XXb$ such that $E_Q|u(\Xb_s)|<\infty$ for all $0\le s\le T.$ 
\begin{enumerate}
\item
We say
that $u$ admits a stochastic derivative  under $Q$ at time
$t\in[0,T)$ if   the following limit
\begin{align}\label{eq-F03}
 L^Qu(t,X _{ [0,t]}):=\Limh E_Q\left(\left.\frac1h [u(\Xb_{t+h})-u(\Xb_t)]
     \,\right|\,  X _{ [0,t]}  \right)
\end{align}
exists in $L^1(Q).$
\\
In this case,  $L^Qu(t,\sbt)$ is called  the stochastic derivative of $u$ at time $t$.
\item
If $u$ admits a stochastic derivative for  almost all
$t,$ we say that $u$ belongs to the domain $\dom L^Q $ of
the stochastic derivative $L^Q$ of $Q.$
\end{enumerate}
\end{definition}

\begin{remark}[Special  case where $Q$ is Markov] 
If $Q\in\MO$ is Markov, it is immediate to see that for any function $u\in\dom L^Q,$
\begin{align*}
 L^Qu(t,X _{ [0,t]})= L^Qu(t,X _t),\quad Q\ae,\quad \textrm{for all }t,
\end{align*}
with an obvious abuse of notation.
The right-hand side of this identity defines a function    on $\iX$ because for every $t$,  $L^Qu(t,X _t)$ is a conditional expectation (with respect to $X_t$)  and on a Polish space any conditional expectation admits a regular version. In the present situation, the function $(t,x)\mapsto \LQ u(t,x)$ is defined $\qq_t\ae$ for every $t.$ 
\\
While for all $t,$ $\LQ u(t,\sbt)$ is measurable, it does not follow from Definition \ref{defh-04} that $L^Q u$ is jointly measurable on $\iX.$ However, when $u$ fulfils the hypotheses of Proposition \ref{resh-01} below, that is: $u$ is in $\dom\LQ$ and satisfies
$E_Q\Iii
 \big|\LQ u(t,X _{ [0,t]})\big|^p\,dt<\infty$ for some  $p \ge 1,$ if follows from this proposition that $L^Qu=\LQ u$, implying that $L^Qu$ is jointly measurable.
\end{remark}

\subsection*{A convolution lemma}

  Next technical result will be used at several places in the rest of the article.

\begin{lemma}\label{resh-21}
For all $h>0,$ let
$k^h$ be a measurable nonnegative  convolution kernel such that $\supp
k^h\subset[-h,h]$ and $\int_{\mathbb{R}}k^h(s)\,ds=1.$
Let $Q$ be a $ \sigma$-finite positive measure on $\OO$ and  $v$ be a function in $L^p(\Qb)$ with    $1\le p< \infty$. 
\\
Define for all $h>0, t\in\ii$ and $\omega\in\OO$,\   $k^h *
v(t, \omega):=\int_{\ii}k^h(t-s)v_s( \omega)\,ds$.\\
Then, $k^h *v$ is in $L^p(\Qb)$ and
    $
    \Limh k^h *v=v\  \textrm{in }L^p(\Qb).
    $
\end{lemma}

We see that $k^h (s)\, ds$ is a probability measure on $\RR$ which
converges weakly to the Dirac measure $\delta_0$ as $h$ tends down to zero. 
\\
 We shall  use this lemma with $p=1$ or $2,$ and with $k^h=\frac 1h \1_{[-h,0]}$ or $\frac 1h \1_{[0,h]}.$ 

\begin{proof}
In this lemma, we endow $\OO$ with the Skorokhod topology  which turns it into a Polish space and has the interesting property that its Borel $\sigma$-field matches with the canonical $\sigma$-field.

We  start the proof by
showing that $k^h *v\in L^p(\Qb)$ and more precisely
\begin{equation}\label{eqh-42}
    \|k^h*v\|_{L^p(\Qb)}\le \|v\|_{L^p(\Qb)}<\infty.
\end{equation}
Since $v\in L^p(\Qb)$, for $Q$-almost all $\omega,$
$v(\sbt,\omega)\in L^p(\ii)$. By Jensen's inequality applied with the probability measure $k^h(t)dt$, and the standard $L^1$ estimate of convolution
\begin{multline*}
\Iii |k^h*v(t, \omega)|^p\,dt\le \Iii dt\Iii |v(s, \omega)|^pk^h(t-s)\,ds=\|k^h* |v(\sbt, \omega)|^p\|_1\\
 \le \|k^h\|_1\| |v(\sbt, \omega)|^p\|_1=\Iii |v(t, \omega)|^p\,dt.
\end{multline*}
Letting $h\to 0^+,$
\begin{equation}\label{eqh-50}
\|k^h*v(\sbt,\omega)\|^p_{L^p(\ii)}\le
\|v(\sbt,\omega)\|^p_{L^p(\ii)}.
\end{equation}
Integrating  with
respect to  $Q(d\omega)$ leads to \eqref{eqh-42}.

Now, we prove the convergence. We first show that   the proof can be reduced to the case where $Q$ is a bounded measure. As $Q$ is $ \sigma$-finite, there is an increasing sequence $( \Omega_m)$  of measurable subsets of $\Omega$ such that $\cup_m \Omega_m= \Omega$ and $Q( \Omega_m)< \infty$ for all $m$. With   \eqref{eqh-42} for $\1 _{ \Ob\setminus\Ob_m}\,\Qb= \overline{\1 _{ \OO\setminus\OO_m}\, Q} $ instead of $\Qb$  and the dominated convergence theorem, we see that
\begin{equation}\label{eqh-59}
\begin{split}
0\le \|\1 _{ \Ob\setminus\Ob_m}\times (k^h*v)\|_{L^p(\Qb)}
	&= \| k^h*v\|_{L^p(\1 _{ \Ob\setminus\Ob_m}\,\Qb)}\\
	&\le  \| v\|_{L^p(\1 _{ \Ob\setminus\Ob_m}\,\Qb)}=
	 \|\1 _{ \Ob\setminus\Ob_m}v\|_{L^p(\Qb)} \underset{m\to\infty}{\to} 0
\end{split}
\end{equation}
where $\Ob_m:=\ii\times\OO_m.$ Considering $\1 _{ \OO_m}Q$ for arbitrarily large $m$ instead of $Q$, one can assume without loss of generality that $Q$ is bounded. By the same token, since $Q$ is a  bounded  nonnegative measure on a Polish space, it is tight: there exists a compact subset   containing an arbitrarily close to 1 proportion of the mass of $Q.$ Hence, one can assume without loss of generality that $Q$ is a bounded nonnegative measure with a compact support $\OO_o\subset \OO.$
\\
As $\Ob_o:=\ii\times \OO_o$ is Polish, the space $C(\Ob_o)$ of all   continuous functions on the compact set  $\Ob_o$ is dense in $L^p(\Qb)$, remember that $1\le p< \infty$. Therefore  we can  approximate
$v$ in $L^p(\Qb)$ by a sequence $\seq vn$ in $C(\Ob_o).$
For all $h$ and $n$
\begin{eqnarray*}
 \|k^h *v-v\|_{L^p(\Qb)}
  &\le& \|k^h *(v-v_n)\|_{L^p(\Qb)} +\|k^h *v_n-v_n\|_{L^p(\Qb)}
    +\|v_n-v\|_{L^p(\Qb)}  \\
   &\le& \|k^h *v_n-v_n\|_{L^p(\Qb)}+2 \|v-v_n\|_{L^p(\Qb)}
\end{eqnarray*}
where we used \eqref{eqh-42}.
Take an arbitrary small $\eta>0$ and choose $n$ large enough for
$\|v-v_n\|_{L^p(\Qb)}\le\eta$ to hold. Then,
\begin{equation}\label{eqh-01}
     \|k^h *v-v\|_{L^p(\Qb)}\le
     \|k^h *v_n-v_n\|_{L^p(\Qb)}+2\eta.
\end{equation}
Fix this $n.$ Since $\Ob_o$ is compact,  $v_n\in C(\Ob_o)$  is a
uniformly continuous function. Therefore, for all $\eta>0,$ there
exists $h(\eta)>0$ such that for any $t,t',\omega,\omega'$
satisfying $|t-t'|+d_{\OO_o}(\omega,\omega')\le h(\eta)$ (with $d_{\OO_o}$ the Skorokhod distance),  we have
$|v_n(t',\omega')-v_n(t,\omega)|\le\eta.$   In particular, with $\omega=\omega',$
we see that
\begin{equation*}
    |t'-t|\le h(\eta) \Rightarrow
    \sup_{\omega\in\OO_o}|v_n(t',\omega)-v_n(t,\omega)|\le\eta.
\end{equation*}
Because of the property: $\supp k^h\subset[-h,h],$ we deduce from
this that for any $\omega\in\OO_o,$
$|k^h*v_n(t, \omega)-v_n(t, \omega)|\le\Iii |v_n(t-s, \omega)-v_n(t, \omega)|k^h(s)\,ds\le\eta$
as soon as $h\le h(\eta)/2.$ We used the convention that $v_n(t, \omega)$ vanishes when $t$ is outside $\ii.$ Consequently $ \|k^h
*v_n-v_n\|_{L^p(\Qb)}\le \eta\, Q(\OO) ^{ 1/p}.$ Finally, with \eqref{eqh-01}
this leads us to $\|k^h *v-v\|_{L^p(\Qb)}\le [2+Q(\OO) ^{ 1/p}]\,\eta.$
Since $\eta$ is arbitrary and $Q(\OO)$ is finite, this shows that
    $
    \lim_{h\rightarrow0}\|k^h
    *v-v\|_{L^p(\Qb)}=0,
    $
which is the desired result.
\end{proof}

\begin{corollary}\label{res-14b}
Assume that in addition to the hypotheses of Lemma \ref{resh-21}, for any $0\le t\le T,$ the random variable $v(t,\sbt)$ is $ \sigma(X _{ [0,t]})$-measurable, resp.\ $ \sigma(X_t)$-measurable.  Then,  the process $v^h$  defined by $v^h(t, \omega):= E_Q[k^h*v(t)\mid X _{ [0,t]}= \omega _{ [0,t]}],$ resp.\ $v^h(t, \omega):= E_Q[k^h*v(t)\mid X_t= \omega_t],$ is in  $L^p(\Qb)$ and $\Limh v^h=v$ in $L^p(\Qb).$
\end{corollary}

\begin{proof}
By Jensen's inequality
\begin{align*}
\|v^h-v\| _{ L^p(\Qb)}^p
	&=\int _{ \XXb} |E_Q\left.[k^h*v(t)  \,\right|\, X_{[0,t]}]-v(t)|^p\,d\Qb
	=\int _{ \XXb} |E_Q\left.[k^h*v(t)-v(t)  \,\right|\, X_{[0,t]}]|^p\,d\Qb\\
	&\le \int _{ \XXb} E_Q[\left.|k^h*v(t)-v(t)|^p  \,\right|\, X_{[0,t]}]\,d\Qb
	=\int _{ \XXb} |k^h*v(t)-v(t)|^p  \,d\Qb\\
	&=\|k^h*v-v\| _{L^p(\Qb)}^p \underset{h\to 0^+}\longrightarrow 0,
\end{align*}
where the vanishing limit is the content of Lemma \ref{resh-21}. Replace  $X _{ [0,t]}$ by $X_t$ for the other statement.
\end{proof}

\subsection*{Extended generators and stochastic derivatives coincide}
The main result of this section is the following Proposition \ref{resh-01} which states that whenever $u$ is in the domain of the extended generator $\LQ,$ one can compute  $\LQ u$  using the stochastic derivative:
\begin{align}\label{eq-F04}
\LQ u=L^Qu,\quad \Qb\ae
\end{align}
On the other hand, it will be proved later at Proposition \ref{resh-03} that whenever  the limit \eqref{eq-F03} defining the stochastic derivative $L^Qu$ exists in $L^1(\Qb),$  $u$ is also in the domain of the extended generator and \eqref{eq-F04} is satisfied.

\begin{proposition}\label{resh-01} For any $Q\in\MO,$ if  $u$ is in $\dom\LQ$ and satisfies
$E_Q\Iii
 \big|\LQ u(t,X _{ [0,t]})\big|^p\,dt<\infty$ for some  $p \ge 1,$ then
\begin{equation}\label{eqh-43}
    \Limh E_Q\II0{T-h} \left|
    \frac 1hE_Q\left.\Big[ u(\Xb_{t+h})-u(\Xb_t) \,\right|\, X _{ [0,t]}\Big]-\LQ u(t,X _{ [0,t]})\right|^p\,dt=0.
\end{equation}
In particular, this implies that
$u\in\dom L^Q,$ and
the limit 
$$
 \LQ u(t, X _{ [0,t]})=L^Q u(t, X _{ [0,t]}):=\Limh \frac 1hE_Q\left.\Big[u(\Xb_{t+h})-u(\Xb_t)\,\right|\,  X _{ [0,t]}  \Big]$$ takes place in $L ^p(\Qb).$
\end{proposition}


\begin{proof}  
The specific convolution
kernel $k^h=\frac 1h \1_{[-h,0]}$ gives 
\[
k^h*v_t=\frac 1h \int _{ [t,(t+h)\wedge T]} v_s\,ds.
\]
Denoting
$v_t:= \LQ u(t, X _{ [0,t]})$, with the
definition of the extended generator,  we see that $h\mapsto[u(\Xb_{(t+h)\wedge T})-u(\Xb_t)]-h k^h*v_t$ is a local
martingale with zero expectation.
It follows that there exists a sequence $\seq{\tau}k$ of $[0,T]\cup \left\{ \infty\right\} $-valued stopping times such that $\Lim k \tau_k= \infty,$ $Q\ae,$ and for any $0\le t<T$, $0<h\le T$ and $k\ge1,$
\begin{equation*}
\frac 1hE_Q\left.[u(\Xb_{(t+h)\wedge \tau_k\wedge T})-u(\Xb _{ t\wedge\tau_k})\,\right|\,  X_{[0,t]}]
=E_Q[k^h*v _{ t\wedge \tau_k}\mid X_{[0,t]}].
\end{equation*}
Since  it is assumed that $E_Q\Iii
| v_t|^p\,dt<\infty,$ by  Jensen's inequality, \eqref{eqh-50} and  dominated convergence, we obtain that the right-hand side  tends to $E_Q[k^h*v _{ t}\mid X_{[0,t]}]$ as $k$ tends to infinity (along a subsequence), for almost all $t$, leading to 
 \begin{equation*}
\frac 1hE_Q[u(\Xb_{(t+h)\wedge T})-u(\Xb _{ t})\mid X_{[0,t]}]
=E_Q[k^h*v _{ t}\mid X_{[0,t]}].
\end{equation*} 
On the other hand, as   $v_t$  is $\sigma(X_{[0,t]})$-measurable, we see with Corollary \ref{res-14b} that \\
$\Limh E_Q\Iii|k^h*v_t-v_t|^p\,dt=0.$ Gathering these considerations leads to \eqref{eqh-43}.
\end{proof}

We provide a  result which is complementary to Proposition \ref{resh-01}, below at Proposition \ref{resh-03}. Its  proof  relies on the following easy analytic result.

\begin{lemma}\label{resh-24}
Let $a,b$ be two measurable
functions on $\ii$ such that $a$ is right continuous, $b$ is
Lebesgue-integrable and
    $
    \Limh \int_{[0,T-h]}\left|\frac1h\{a(t+h)-a(t)\}-b(t)\right|\,dt=0.
    $
Then, $a$ is absolutely continuous and its distributional
derivative is $\dot a=b.$
\end{lemma}

\begin{proof}
Remark first that $t\mapsto\1_{\{0\le t\le
T-h\}}\frac1h\{a(t+h)-a(t)\}$ is integrable for any small enough 
$0<h\le T.$ Take any $0\le r\le s<T.$ On one hand, we have
    $\Limh \II rs  \frac1h\{a(t+h)-a(t)\}\,dt=\II rs b(t)\,dt$
and on the other one:
    $\II rs  \frac1h\{a(t+h)-a(t)\}\,dt=\frac1h \II s{s+h}a(t)\,dt -\frac1h \II r{r+h}a(t)\,dt,$
so that with the assumed right continuity of $a$ the integrals $ \II s{s+h}a(t)\,dt$ and  $\II r{r+h}a(t)\,dt$ are well defined for any small enough $h$ and we have
     $\Limh \II rs  \frac1h\{a(t+h)-a(t)\}\,dt=a(s)-a(r).$
Therefore $a(s)-a(r)=\II rs b(t)\,dt$ which is the claimed
property.
\end{proof}

\begin{proposition}\label{resh-03}
Let  $u$ be  a  measurable real function on $\XXb$, and  $v$ be an adapted process
such that $u(\Xb)$ and $v$ are $\Qb$-integrable,  $t\mapsto u(\Xb_t)$ is right continuous (for instance $u$ might be continuous) and
 \begin{equation}\label{eqh-03}
    \Limh E_Q\II0{T-h}\left|\frac1h E_Q\left. [u(\Xb_{t+h})-u(\Xb_{t})\,\right|\,  X _{ [0,t]}]-v_t\right|\,dt =0.
\end{equation}
Then, $u$ belongs to $\dom\LQ$ and $\dom L^Q$,  and
    $
      \LQ u=L^Q u=v,\ \Qb\ae
    $
\end{proposition}

\begin{proof}

 We
write $E=E_Q$ and $u_t=u(\Xb_t)$   to simplify the notation.
Fix $0\le
r<T.$ We have
\begin{multline*}
    \left|\left. E\left[\II r{T-h}\left(\frac1h\{u_{t+h}-u_t \}-v_t\right)\,dt\,\right|\,  X _{[0,r]}\right]\right| \\
   \le  E\left.\left[\II r{T-h}E\left.\left(\big|\frac1h\{u_{t+h}-u_t \}-v_t\big|\,\right|\,  X _{[0,t]}\right)\,dt\,\right|\, 
   X _{ [0,r]}\right].
\end{multline*}
With \eqref{eqh-03} and Fatou's lemma, we obtain
\begin{multline*}
  E\left(\liminf_{h\downarrow0} \left|E\left.\left[\II r{T-h}\left(\frac1h\{u_{t+h}-u_t \}-v_t\right)\,dt\,\right|\,  X _{[0,r]}\right]\right|\right) \\
   \le \Limh  E \II r{T-h}E\left(\left.\big|\frac1h\{u_{t+h}-u_t \}-v_t\big|\,\right|\,  X _{[0,t]}\right)\,dt=0.
\end{multline*}
Hence, there exists a sequence $\seq hn$ of positive numbers such
that $\Lim n h_n=0$ and
\begin{equation*}
    \Lim n \II r{T-h_n} \Big|E\left[\left.\left(\frac1{h_n}\{u_{t+h_n}-u_t \}-v_t\right)\,\right|\, 
   X_r\right]\Big|\,dt=0,\quad Q\ae
\end{equation*}
It remains to apply Lemma \ref{resh-24} with
    $a(t)=E\left[u_t\mid X _{[0,r]}\right]$ and
    $b(t)=E\left[v_t\mid X _{[0,r]}\right]$
to see that for all $0\le r\le s<T,$
    $E\left[u_s-u_r-\II rs v_t\,dt\mid X _{ [0,r]}\right]=0.$
As $u_t$ is $Q$-integrable since it is  assumed that $u(\Xb)$ is   $\Qb$-integrable and $t\mapsto u(\Xb_t)$ is right-continuous, this proves that $M$ is a $Q$-martingale where
$M_s:=u(\Xb_s)-u(\Xb_0)-\II0s v_t\,dt.$
Therefore, $u$
belongs to $\dom\LQ$ and $\LQ u=v$. To obtain the remaining identity $\LQ u=L^Q u$, apply  Proposition \ref{resh-01} with $p=1.$
\end{proof}

\begin{corollary}\label{res-19}
Let  $ (\mathcal{G}_t) _{ 0\le t\le T}$ be the generator of the semigroup $(T_s^t) _{ 0\le s\le t\le T}$  associated to the Markov measure $Q$, and let  $u:\XXb\to\RR$ be an $x$-continuous and $t$-differentiable function  such that for each $t,$ $u(t,\sbt)\in\dom \mathcal{G}_t,$ and 
$ 
\Limh \Iii \sup _{\XX}| h ^{ -1}(T _t^{ t+h}-\Id)u- ( \partial_t+ \mathcal{G}_t)u|\,dt
=0.
$ 
Then $u$ belongs to $\dom\LQ$ and \  
$
 \LQ u=(\partial_t+ \mathcal{G}_t)u.
$
\end{corollary}

\begin{proof}
Immediate consequence of Proposition \ref{resh-03}.
\end{proof}

\section{Stochastic derivatives. Extensions} \label{sec-sd-ext}

The results of previous section are extended by means of the notions of $P$-locality and integration times. The main result of this section is Proposition \ref{res-H02c}.

Before this, we start introducing the backward in time analogues of the already defined (forward) extended generators and stochastic derivatives.

\subsection*{Reversing time}

Let $Q\in\MO$ be any path measure. Its time reversal is 
\begin{align*}
Q^*:=(X^*)\pf Q\in\MO,
\end{align*}
where 
\begin{align*}
\left\{ \begin{array}{ll}
X^*_t:= \Limh X _{ T-t+h},\quad & 0\le t< T,\\
 X^*_T:=X_0,& t=T,
\end{array}\right.
\end{align*}
is the reversed canonical process. We assume that $Q$ is such that $Q(X _{ T^-}\neq X_T)=0,$ i.e.\ its sample paths are left-continuous at $t=T.$ This implies that the time reversal mapping $X^*$ is (almost surely) one-one on $\OO.$ Similarly, we define
\begin{align*}
\Xb^*(t, \omega):=(T-t, X^*_t( \omega)), \quad 0\le t\le T.
\end{align*}
 We introduce the backward extended generator and the backward stochastic derivative
\begin{align}\label{eq-F05}
\begin{split}
\LLb^Qu(t, X _{ [t,T]})&:=\LLf ^{ Q^*} u^*(t^*,X^* _{ [0,t^*]}),\\
\Lb^Qu(t, X _{ [t,T]})&:=\Lf ^{ Q^*} u^*(t^*,X^* _{ [0,t^*]}),
\end{split}
\end{align}
where $u^*:\XXb\to\RR$ is defined by $ u^*(t^*, \omega^* _{ [0,t^*]}):=u(t, \omega _{ [t,T]}),$  with $t^*:=(T-t)^+$, $ \omega^*(t):= \omega(t^*)$, and $ \LLf ^{ Q^*}$ and  $\Lf ^{ Q^*}$ stand respectively for the standard (forward) generator and derivative of $Q^*$ as introduced at Definitions \ref{defh-05} and \ref{defh-04}. Definitions \eqref{eq-F05} match with the following ones.

As a notation, the $ \sigma$-field generated by $ X _{ [t^-,T]}$ is $ \sigma(X _{ [t^-,T]}):= \cap _{ h>0} \sigma(X _{ [t-h,T]})= \sigma(X _{ t^-})\vee \sigma(X _{ [t,T]}).$

\begin{definition}[Extended backward generator]
Let $Q$ be a conditionable path measure. A process $u$ adapted to the predictable backward filtration $( \sigma(X _{ [t^-,T]}); 0\le t\le T)$ is said to be in
the domain of the extended backward generator of $Q$ if there exists a
process $v$ also  adapted to the predictable backward filtration such that
 $\Iii|v(t,X _{ [t^-,T]})|\,dt<\infty,$ $Q\ae$ and the process 
 $$
u(t, X _{ [t^-,T]})-u(T,X_T)-\II tT v(s, X _{ [s^-,T]})\,ds,\quad 0\le t\le T,
 $$ 
 is a local backward $Q$-martingale.
We denote
$$
   \LLb^Q u(t, X _{ [t^-,T]}):=v(t, X _{ [t^-,T]})
$$
and call $\LLb^Q $ the extended backward generator of $Q.$ Its domain
  is denoted by $\dom
\LLb^Q .$
\end{definition}

Remark that denoting  $ \LLb^Q=v$ is consistent with \eqref{eq-F05}.

\begin{definition}[Stochastic backward derivative]\label{defh-04b}
Let $Q$ be a conditionable path measure and a measurable function $u$ on $\iX$  such that $E_Q|u(s, X _s)|<\infty$ for all $0\le s\le T.$ 
\begin{enumerate}
\item
We say
that $u$ admits a stochastic backward derivative  under $Q$ at time
$t\in(0,T]$ if   the following limit
\begin{align*}
 \Lb^Qu(t,X _{ [t^-,T]}):=\Limh E_Q\left.\left(\frac1h [u(\Xb_{t-h})-u(\Xb_t)]
    \,\right|\,   X _{ [t^-,T]}  \right)
\end{align*}
exists in $L^1(Q).$  
\\
In this case,  $\Lb^Qu(t,\sbt)$ is called  the stochastic backward derivative of $u$ at time $t$.
\item
If $u$ admits a stochastic backward derivative for  almost all
$t,$ we say that $u$ belongs to the domain $\dom \Lb^Q $ of
the stochastic backward derivative $\Lb^Q$ of $Q.$
\end{enumerate}
\end{definition}

Remark that this definition is consistent with \eqref{eq-F05}.

Mimicking almost verbatim the proofs of Propositions \ref{resh-01} and \ref{resh-03}, (consider the convolution kernel $k ^{ -h}:=\frac 1h \1_{[0,h]}$ instead of $k^h=\frac 1h \1_{[-h,0]}$), we arrive at

\begin{proposition}\label{resh-01b} Let $Q$ be a conditionable path measure. 
\begin{enumerate}[(a)]
\item
If  $u$ is in $\dom\LLb^Q$ and  such that
$E_Q\Iii
 \big|\LLb^Q u(t,X _{ [t^-,T]})\big|^p\,dt<\infty$ for some  $p \ge 1,$ then
\begin{equation*}
    \Limh E_Q\II h{T} \left|
    \frac 1hE_Q\left.\Big[u(\Xb_{t-h})-u(\Xb_t)\,\right|\,  X _{ [t^-,T]} \Big]-\LLb^Q u(t,X _{ [t^-,T]})\right|^p\,dt=0.
\end{equation*}
In particular, this implies that
$u\in\dom \Lb^Q,$ and
the limit 
$$\LLb^Q u(t, X _{ [t^-,T]})=\Lb^Q u(t, X _{ [t^-,T]}):= \Limh \frac 1hE_Q\big[u(\Xb _{ t-h})-u(\Xb_t)\mid X _{ [t^-,T]} \big]$$ takes place in $L ^p(\Qb).$

\item
Let $u$ be a  measurable real function on $\XXb$ and $v$ a process adapted to the predictable backward filtration,
such that $u(\Xb), v$ are $\Qb$-integrable,  $t\mapsto u(\Xb^*_t)$ is right continuous (for instance $u$ might be continuous) and
 \begin{equation*}
    \Limh E_Q\II h{T}\left|\frac1h E_Q[u(\Xb_{t-h})-u(\Xb_{t})\mid X _{ [t^-,T]}]-v(t, X _{ [t^-,T]})\right|\,dt =0.
\end{equation*}
Then, $u$ belongs to $\dom\LLb^Q$ and $\dom \Lb^Q$,  and
    $
      \LLb^Q u=\Lb^Q u=v,\ \Qb\ae
    $
\end{enumerate}
\end{proposition}

When working with forward generators and derivatives in contexts where time reversal is not mentioned, we often omit the forward arrow,  denoting: $\LLf= \mathcal{L}$ and $\Lf=L.$

\subsection*{$P$-locality}
We  present some notions which will be useful when looking at the extended HJB equation at Section \ref{sec-HJB}.
Let us recall  a slight modification of the standard definitions of stochastic integral and  local martingale, introduced in \cite{Leo11a}.

\begin{definition}[$P$-locality] \label{def-z01}
Let $P\in\PO,$ $Q\in\MO$ such that $P\ll Q.$
\begin{enumerate}
\item
A process $M$ is said to be a $P$-local $Q$-martingale if there exists an increasing sequence of $\ii \cup \{\infty\}$-valued  stopping times $\seq \tau k$ such that $\Lim k \tau_k =\infty,$ $P\ae$ (and not necessarily $Q\ae$) such that the stopped processes $M^{\tau_k}$ are $Q$-martingales, for all $k\ge1.$
\item
A process $Y$ is said to be a $P$-local $Q$-stochastic integral if there exists an increasing sequence of $\ii \cup \{\infty\}$-valued  stopping times $\seq \tau k$ such that $\Lim k \tau_k =\infty,$ $P\ae$ (and not necessarily  $Q\ae$) such that the stopped processes $Y^{\tau_k}$ are $L^2(Q)$-stochastic integrals, for all $k\ge1.$
\end{enumerate}
\end{definition}

The filtration is the $Q$-completion of the canonical filtration. Since any local $Q$-martingale admits a  càdlàg $Q$-version and $P\ll Q,$ any $P$-local $Q$-martingale admits a  càdlàg $P$-version.

In connection with the notion of $P$-locality, we introduce a modification of the notion of extended generator.

\begin{definition}[$P$-local extended generator of the path measure $Q$]\ \label{defh-05b}
Let $P\in\PO,$ $Q\in\MO$ such that $P\ll Q.$
\\
A measurable function $u$ on $\XXb$ is said to be in
the domain of the $P$-local extended forward generator of the path measure $Q$ if there exists a $\Pb$-almost everywhere defined 
adapted process $v$    such that
 $\Iii|v_t|\,dt<\infty,$ $P\ae$ and the process $u(\Xb_t)-u(\Xb_0)-\II0t v_s\,ds,$\ $0\le t\le T,$ is a $P$-local $Q$-martingale.
We denote
$$
   \LLLf QP u:=v,\quad 
$$
and call $ \LLLf QP  $ the $P$-local extended  forward generator of $Q.$ Its domain is denoted by $\dom
\LLLf QP .$

A similar definition holds for the $P$-local extended \emph{backward} generator: $\LLLb QP$, of $Q$.

\end{definition}

The $P$-local extended generator of $Q$  only consists in requiring that the local $Q$-martingales   entering the Definitions \ref{defh-05} and \ref{defh-05b} of $ \LLf^Q$ and $\LLb^Q$ are replaced by {$P$-local} $Q$-martingales.
This notion will allow us to extend in a natural way Itô formula  in a diffusion setting at Proposition \ref{res-14}.

\subsection*{Integration times}
 
In order to motivate Definition \ref{def-H07} below, let us start with a remark.
Since  $u(\Xb)$ might not be $Q$-integrable,  the conditional increment $E_Q[u(\Xb_{t+h})-u(\Xb_t)\mid X _{ [0,t]}]$  appearing in the expression \eqref{eq-F03} of the stochastic derivative, might not be meaningful. To take this trouble into account, let us introduce the following notion.

\begin{definition}[Integration time] \label{def-H07}
Let $u$ be in $\dom\LQ$ and $\tau$ be a stopping
time. We say that $\tau$ is a $Q$-integration time of $u$ if the stopped local martingale $M^{u, \tau}$ (recall Definition \ref{defh-05}) and the stopped process  $u(\Xb ^{ \tau})$ are  $\Qb$-integrable.
\end{definition}

Clearly, for any $Q$-integration time $ \tau$ of $u$,   
$E_Q[u(\Xb^\tau_{t+h})-u(\Xb^\tau_t)\mid X _{ [0,t]}]$ is well defined. 
\\
Of course, for the stopping time $ \tau$ to be a $Q$-integration time of $u$ it is necessary and sufficient that two of the following properties hold:
\begin{enumerate}[\qquad -]
\item
$u(\Xb ^{ \tau})$ is  $\Qb$-integrable;
\item
$M^{u, \tau}$  is  $\Qb$-integrable;
\item
$(t, \omega)\mapsto\1 _{\{ \tau( \omega)>t\}}\LQ u(t, \omega _{ [0,t]})$ is  $\Qb$-integrable.
\end{enumerate}
In which case the three properties are satisfied.

\begin{lemma}\label{res-H22}
 Let us take a conditionable path measure $Q,$  a function $u$ in $\dom\LQ $ and $\tau$ any $Q$-integration time of  $u$.
  Then,
\begin{equation}\label{eq-H43b}
    \Limh E_Q\II0{T-h} \1 _{ \left\{\tau>t\right\} }  \left|
    \frac 1hE_Q\left.[u(\Xb^\tau_{t+h})-u(\Xb^\tau_t)\,\right|\,  X _{ [0,t]}]-\LQ u(t, X _{ [0,t]})\right|\,dt=0.
\end{equation}
\end{lemma}

\begin{proof}
Denoting $Q^\tau:=(X^\tau)\pf Q$ the law of the stopped canonical process $X^\tau$ at the random time $ \tau$, for any $0\le t\le T$ and any bounded  measurable function $U:X^\tau(\OO)\to\RR,$  as a consequence of general considerations in measure theory, we have
\begin{align*}
E _{ Q^\tau}(U\mid X _{ [0,t]})\circ X^\tau
	=E_Q(U(X^\tau)\mid X^\tau _{ [0,t]}),
	\qquad Q\ae
\end{align*}
 It follows that
\begin{align*}
&E _{ Q^\tau}\II0{T-h}  \left|
    \frac 1hE _{ Q^\tau}[u(\Xb_{t+h})-u(\Xb_t)\mid X  _{ [0,t]}]- \mathcal{L} ^{ Q^\tau} u(t, X _{ [0,t]})\right|\,dt\\
    =\ & E _{ Q}\II0{T-h}  \left|
    \frac 1hE _{ Q^\tau}[u(\Xb_{t+h})-u(\Xb_t)\mid X  _{ [0,t]}]\circ X^\tau- \mathcal{L} ^{ Q^\tau} u(t, X _{ [0,t]})\circ X^\tau\right|\,dt\\
     =\ & E _{ Q}\II0{T-h}  \left|
    \frac 1hE _{ Q}[u(\Xb^\tau_{t+h})-u(\Xb^\tau_t)\mid X^\tau  _{ [0,t]}]- \mathcal{L} ^{ Q^\tau} u(t, X^\tau_{ [0,t]})\right|\,dt\\
    =\ & E _{ Q}\II0{T-h}  \left|
    \frac 1hE _{ Q}[u(\Xb^\tau_{t+h})-u(\Xb^\tau_t)\mid X^\tau  _{ [0,t]}]-  \1 _{ \left\{ \tau>t\right\} }\ \LQ u(t, X _{ [0,t]})\right|\,dt,
\end{align*}
where at last equality we use
\begin{align*}
\mathcal{L} ^{ Q^\tau}u(t,X^\tau _{ [0,t]})
	= \1 _{ \left\{ \tau>t\right\} }\ \LQ u(t, X _{ [0,t]}),\quad \textrm{for almost all }t,
\end{align*}
which is a direct consequence of  the very definitions of $ \mathcal{L} ^Qu$ and $ \mathcal{L} ^{ Q^\tau}u.$ 
\\
Applying Proposition \ref{resh-01} to $Q^\tau$ with $p=1,$ this identity gives us 
\begin{equation*}
    \Limh E_Q\II0{T-h}  \left|
    \frac 1hE_Q[u(\Xb^\tau_{t+h})-u(\Xb^\tau_t)\mid X^\tau  _{ [0,t]}]-\1 _{ \left\{\tau>t\right\} }\LQ u(t, X _{ [0,t]})\right|\,dt=0,
\end{equation*}
which is equivalent to \eqref{eq-H43b}.
\end{proof}

\begin{remarks}\label{rem-H04}\ \begin{enumerate}[(a)]
\item
Writing $ \Limh E_Q\left(\frac1h [u(\Xb ^{ \tau}_{t+h})-u(\Xb_t ^\tau)]
    \mid  X_{[0,t]}^{\tau} \right)$ is a slight abuse of notation. It should be written $\Limh E_Q\left(\frac1h [u(\Xb ^{ \tau}_{(t+h)\wedge T})-u(\Xb_t^\tau)]  \mid  X_{[0,t]}^{\tau}\right),$ with $(t+h)\wedge T$ instead of $t+h.$
 This simplification will be kept at several places in the sequel.
 \item
It is necessary that $Q(\tau>t\mid X_{[0,t]})>0$ for \eqref{eq-H43b} to be a nontrivial assertion. { Lemma \ref{res-H26} below  tells us that such stopping times $\tau$ exist. }
\end{enumerate}\end{remarks}

\begin{remark}[The Markov case]\label{rem-H05}
Denoting $DU_t^\tau:=u(\Xb^\tau_{t+h})-u(\Xb^\tau_t)$,
it is tempting to infer from Lemma \ref{res-H22} that for any Markov measure $Q$,
\begin{equation*}
E_Q[DU_t^\tau\mid X_{[0,t]}]=\1 _{ \left\{\tau>t\right\}}  E_Q[DU_t^\tau\mid X_t].
\end{equation*}
But \emph{this is false in general}, unless $ \tau$ is a ``Markov stopping time'', i.e.\ unless $X^\tau\pf Q$ is a Markov measure. 
More precisely,   in general  last equality in:
$E_Q(DU_t^\tau\mid X_{[0,t]})
= \1 _{ \left\{\tau>t\right\}} E_Q(DU_t^\tau\mid X_{[0,t]})
= \1 _{ \left\{\tau>t\right\}} E_Q(DU_t^\tau\mid X_t),
$
fails.
\\
Note in passing that \eqref{eq-H43b} is a statement about the stopped path measure $X^\tau\pf Q.$
\end{remark}

Lemma \ref{res-H22} suggests a way to compute $\LQ u$ when integrability is lacking.  It is stated below at Proposition \ref{res-H02c}. Let us start with a simple remark stated as a lemma.

\begin{lemma}\label{res-H26}
Any function in $\dom\LQ$ admits a sequence of $Q$-integration times tending $Q$-almost everywhere to infinity.
\end{lemma}

\begin{proof}
For any $u\in \dom\LQ,$ we denote
$
M^u_t:=u(\Xb_t)-u(\Xb_0)-\int _{ [0,t]}\LQ u(s, X _{ [0,s]})\,ds,
$
the local $Q$-martingale which appears at Definition \ref{defh-05}. By the very definition of a local $Q$-martingale, there exists a sequence $( \sigma_k)$ of stopping times such that $\Lim k\sigma_k= \infty,$ $Q\ae$ and the stopped local martingales $M ^{ u, \sigma_k}$ are genuine uniformly integrable $Q$-martingales. Furthermore, since it is assumed that $\Iii |\LQ u(t, X _{ [0,t]})|\,dt< \infty$, $Q\ae,$ the stopping times $ \theta_k:= \inf_{}\{t: \int _{ [0,t]}|\LQ u(s, X _{ [0,s]})|\,ds\ge k\}$ satisfy $\Lim k \theta_k= \infty, Q\ae$ Therefore, the sequence $(\tau_k:=\sigma_k\wedge \theta_k)$ tends almost surely to infinity. We also see that for any $k$,  $M ^{u, \tau_k }$ and $u(\Xb ^{ \tau_k})$ are uniformly $Q$-integrable. 
\end{proof}

\begin{remark}[Integration time trick]\label{rem-H08}
Note that if $( \tau_k)$ is a sequence of $Q$-integration times tending $Q\ae$  to infinity, and $( \tau'_k)$  is any sequence of stopping times tending $Q\ae$ to infinity, then as a consequence of the stopped martingale theorem, $ \tau''_k:= \tau_k\wedge \tau'_k,$  defines a sequence of $Q$-integration times tending $Q\ae$  to infinity.  To keep notation easy, in this situation we still write $ \tau_k$ instead of $ \tau''_k.$ Let us call this operation the \emph{integration time trick}.
\end{remark}

\begin{proposition}\label{res-H02c}
Let $Q$ be a conditionable measure, take  a function $u$  in
$\dom\LQ$ and fix $t \in[0,T).$
There exist an increasing sequence $(\tau_k)$ of $Q$-integration times of $u$  and a sequence $(h_n)$ of positive numbers such that $\Lim k \tau_k=\infty,$ $Q\ae$, $\Lim n h_n=0$ and for each $k$ we have 
\begin{equation}\label{eq-H46}
\LQ u(t, X _{ [0,t]})=\Lim n \frac{1}{h_n}E_Q\left.\Big[u(\Xb ^{ \tau_k}_{t+h_n})-u(\Xb ^{ \tau_k}_t)\,\right|\,  X_{[0,t]} \Big],
\quad (\1 _{ \left\{\tau^k>t\right\}} \Qb)\ae
\end{equation}
Moreover, by taking a further subsequence
\begin{equation*}
 \LQ u(t, X _{ [0,t]})= \Lim k \Lim n \frac 1{h_n}E_Q\left. \Big[u(\Xb ^{ \tau_k}_{t+h_n})-u(\Xb ^{ \tau_k}_t)\,\right|\,  X_{[0,t]}  \Big],\quad \Qb\ae
\end{equation*}
\end{proposition}

\begin{proof}
We have already seen at Lemma \ref{res-H26} that there exists a sequence $(\tau'_k)$ of $Q$-integration times which tends almost surely to infinity. By means of the integration time trick (see Remark \ref{rem-H08}), it can be chosen as an increasing sequence by considering $\tau_k=\max _{ 1\le i\le k}\tau' _i.$  The almost everywhere convergence  \eqref{eq-H46} is a direct consequence of the $L^1$-limit \eqref{eq-H43b}. The remaining statement  is an easy consequence of \eqref{eq-H46}, $\Lim k Q(\tau_k>t)=1$, and a diagonal subsequence argument to extract a  sequence $(h_n)$ from the array $(h _{ k,n})$. 
\end{proof}

Next result is  intuitively obvious. Nevertheless, its proof is not as  direct as one could wish. As a corollary of Proposition \ref{res-H02c},  its proof relies on Proposition \ref{resh-01}, again.

\begin{corollary}\label{res-18}
For any Markov measure $Q\in\MO$ and any $u\in \dom\LQ$ 
\begin{align*}
\LQ u(t, X _{ [0,t]})=\LQ u(t,X_t),\qquad \Qb\ae
\end{align*}
(with some obvious abuse of notation).
\end{corollary}

\begin{proof}
Assuming that  for some $t$ the limit in the first term of next series of identities exists $Q\ae$, we have
\begin{align*}
&\Lim n h_n ^{ -1}E_Q\left. \Big[u(\Xb ^{ \tau}_{t+h_n})-u(\Xb ^{ \tau}_t)\,\right|\,  X_{[0,t]}  \Big]\\
	=\ &\Lim n h_n ^{ -1} \Big(E_Q\left. \Big[\1 _{  \left\{t< \tau\right\} }\{u(\Xb ^{ \tau}_{t+h_n})-u(\Xb ^{ \tau}_t)\}\,\right|\,  X_{[0,t]}  \Big]\\
	&\hskip 6cm +E_Q\left. \Big[\1 _{  \left\{t\ge  \tau\right\} }\{u(\Xb ^{ \tau}_{t+h_n})-u(\Xb ^{ \tau}_t)\}\,\right|\,  X_{[0,t]}  \Big]\Big)\\
	=\ &\1 _{  \left\{t< \tau\right\} } \Lim n h_n ^{ -1} E_Q\left. \Big[u(\Xb ^{ \tau}_{t+h_n})-u(\Xb ^{ \tau}_t)\,\right|\,  X_{[0,t]}  \Big] \\
	=\ &\1 _{  \left\{t< \tau\right\} } \Lim n h_n ^{ -1} E_Q\left. \Big[u(\Xb _{t+h_n})-u(\Xb _t)\,\right|\,  X_{[0,t]}  \Big] \\
	=\ &\1 _{  \left\{t< \tau\right\} } \Lim n h_n ^{ -1} E_Q\left. \Big[u(\Xb _{t+h_n})-u(\Xb _t)\,\right|\,  X_t  \Big],
\end{align*}
where all these equalities hold $Q\ae$ The first equality is an obvious decomposition. The second one follows from $t\ge \tau\implies \Xb^\tau _{ t+h_n}=\Xb ^{ \tau}_t=( \tau,X_\tau)\implies u(\Xb ^{ \tau}_{t+h_n})-u(\Xb ^{ \tau}_t)=0$
 for the rightmost term, and the fact that the event $ \left\{t< \tau\right\} $ is $ \sigma(X _{ [0,t]})$-measurable because $ \tau$ is a stopping time, for the leftmost term. The third equality is obvious because the limit is pointwise.
Last equality follows from the assumed Markov property of $Q.$
We conclude with Proposition \ref{res-H02c}.
\end{proof}

Note that although in presence of a stopping time, this proof does not rely on a strong Markov property argument.

\begin{remarks}\ 
\begin{enumerate}[(a)]
\item \emph{About càdlàg versions of local martingales.} 
In the previous two sections, we took some care making precise the conditions for the filtration to fulfill the "usual hypotheses" to imply that martingales admit càdlàg versions. This will not be used in what follows because it will be assumed later that the diffusion field $\aa$ is invertible to assure the Brownian martingale representation   (hence any local martingale is continuous as a Brownian stochastic integral). These precautions about the filtration are   necessary if one wishes to extend our results   to a more general setting, including jumps for instance, as a careful reading of the proof of Theorem \ref{res-13} indicates.
\item \emph{About time reversal.} 
Similarly, time reversal will not play any role in the present article. However, we decided to include considerations about it because it plays a major role in entropic optimal transport theory. In particular,  our previous results about backward generators and derivatives are used in the recent article  \cite{CCGL20} and will be utilized in  future research of the author.
\item
\emph{About $P$-locality.} Of course, if $P$ and $R$ are equivalent path measures, then the notion of $P$-locality is useless. However, there are situations where the scalar potential $V$ is irregular enough for $R$ not to be absolutely continuous with respect to $P$. For instance this happens when $V$ is strong enough for the sample paths of $P$ not to reach the nodal set $ \left\{(t,x): dP_t/dx=0\right\} $, see \cite[Remark.\,4.1]{CWZ} and   \cite{Z85} in a similar but different context.

\end{enumerate}
\end{remarks}

\section{Some more preliminary results}
\label{sec-EDM}

\subsection*{Relative entropy}
Let us start with the definition. Let $A$ be any measurable space. The relative entropy of $\pp \in \mathrm{P}(A)$ with respect to the reference  measure $\rr \in \mathrm{M}(A)$ is 
\begin{align*}
H(\pp |\rr ):=\int_A  \log(d\pp /d\rr )\,d\pp  \in( - \infty, \infty]
\end{align*}
if $\pp $ is absolutely continuous with respect to $\rr $  and  $\log_-(d\pp/d\rr)$ is $\pp $-integrable. We set $H(\pp| \rr)= \infty$ otherwise. 
\\
We set $\log_- x:=\max(-\log x,0) $ and $\log_+ x:=\max(\log x,0). $

\begin{remark}\label{rem-02}
As an immediate consequence of this definition,  any $\pp\in \mathrm{P}(A)$ such that  $H(\pp|\rr)< \infty$ verifies $\log(d\pp/d\rr)\in L^1(\pp).$  However it is useful to define $H(\sbt\mid\rr)$ on the following  subset of $\mathrm{P}(A):$
\begin{align*}
D _{ H(\cdot\mid\rr)}
	:= \left\{\pp\in \mathrm{P}(A); \int_A W\,d\pp< \infty, \textrm{ for some }W:A\to[0, \infty) \textrm{ such that }\int_A e ^{ -W}\,d\rr< \infty\right\} .
\end{align*}
Note that when $\rr$ is a bounded measure, then $W=0$ does the work and $D _{ H(\cdot\mid\rr)}= \mathrm{P}(A).$
\end{remark}

The reason for this choice is guided by the next 
 result which provides us with a criterion for $\log_-(d\pp/d\rr)\in L^1(\pp).$ 

\begin{proposition}\label{res-24}
Let $\rr$ be a $ \sigma$-finite measure on $A.$ 
\begin{enumerate}[(a)]
\item
\begin{enumerate}[(i)]
\item
There exists a measurable function $W:A\to [0, \infty)$ such that 
$	
\int_A e ^{ -W}\,d\rr< \infty.
$	

\item
For any such function $W$ and any $\pp\in \mathrm{P}(A)$ verifying $\pp\ll\rr$ and  $\int_A W\, d\pp< \infty,$ we have: $\log_-(d\pp/d\rr)\in L^1(\pp).$

\end{enumerate}
\item
Moreover, $D _{ H(\cdot\mid\rr)}$ is a convex subset of the vector space of all bounded signed measures and $H(\sbt|\rr)$ is a $(- \infty, \infty]$-valued convex function on this set.
\end{enumerate}
\end{proposition}

\begin{proof}
Take $\pp\in \mathrm{P}(A)$  such that $\pp\ll\rr$ and denote $\rho :=d\pp/d\rr$ for simplicity. 
\\
Remark that when $\rr(A)< \infty$,  because $\sup _{ 0\le z\le 1}z|\log z|=e ^{ -1},$ we always have 
\begin{align*}
\int_A \log_- \rho \,d\pp=
\int _{ \left\{\rho \le 1\right\} } |\log \rho |\, d\pp
	= \int _{ \left\{\rho \le 1\right\} } \rho  |\log \rho |\, d\rr
	\le e ^{ -1}\rr(A)< \infty.
\end{align*}
This corresponds to $W=0.$ 

\Boulette{(a)}
Now, we only assume that $\rr$ is $ \sigma$-finite. Statement (i) is a direct consequence of this hypothesis. Therefore, let us prove (ii). We set $\rr_W:=  e ^{ -W}\,\rr\in \mathrm{M}(A)$ and denote $\rho _W:=d\pp/d\rr_W= e^W \rho .$ For any $\pp\in  \mathrm{P}(A)$, we see  that

\begin{align*}
\int _{ \{\rho \le 1\}}|\log \rho |\, d\pp
	&=\int _{ \{\rho \le 1\}}|\log \rho _W-W|\, d\pp\\
	&\le \int _{ \{\log \rho _W\le W\}}|\log \rho _W|\, d\pp+\int_A W\,d\pp\\
	&\le \int _A \log_- \rho _W\, d\pp+2\int_A W\,d\pp.
\end{align*}
But we know by our preliminary remark that $\int _A \log_- \rho _W\, d\pp$ is finite because $\rr_W(A)=\int_A e ^{ -W}\,d\rr< \infty.$

\Boulette{(b)}
Finally,   $D _{ H(\cdot\mid\rr)}$  is a convex set because for any $\pp,\pp'$ in the convex set $ \mathrm{P}(A)$ such that $\int_A W\, d\pp, \int_A W'\,d\pp'<\infty$ with $W,W':A\to[0, \infty),$ $\int_A (e^W+ e ^{ W'})\,d\rr< \infty,$ we have $\int_A W''\,d( \theta\pp+(1- \theta)\pp')< \infty$ where $0\le  \theta\le 1$ and  $W'':=\min(W,W')$ satisfies $\int_A e ^{ W''}\,d\rr< \infty$. To prove that $H(\sbt|\rr)$ is convex, write it as $H(\pp|\rr)=\int_A [\log\rho_W-W]\,d\pp=\int_A h(\rho_W)\, d\rr_W-\int_A W\,d\pp$ with $h( \rho):= \rho\log  \rho $ a convex function.
\end{proof}

\subsection*{Reference diffusion measure}

Let $\aa$ be a diffusion matrix field on $\iZ$, $\cc$ be some $\ZZ$-valued predictable  process and $ \nu\in\ \mathrm{M}(\ZZ).$ We  say that the path measure $Q \in\MO$ solves the martingale problem with initial measure $ \nu$ and characteristics $(\cc,\aa)$   if $Q _0= \nu$ and
\begin{align*}
\begin{array}{l}
dX_t= \cc_t\,dt+ dM^Q _t,\\
d \langle X\rangle _t=d \langle M^Q \rangle _t=\aa(\Xb_t)\,dt,
\end{array}
\quad \Qb \ae,
\end{align*}
where $M^Q $ is some local $Q $-martingale.    We denote this property by
\begin{align}\label{eq-85}
Q \in\MP(\aa,\cc; \nu).
\end{align}
This is equivalent to: $Q$ solves $\MP( \mathcal{C}, \mathcal{L};Q_0)$ using Definition   \ref{defh-02} with $ \mathcal{C}= \mathcal{C} ^{ 1,2}(\iZ)$ and $ \mathcal{L}=\cc\cdot\nabla + \Delta_\aa/2$.  
We also write shortly $Q\in \MP(\aa,\cc)$ instead of $Q\in\MP(\aa,\cc;Q_0).$
\\
It is implicitly assumed  that $\Iii |\cc_t|\,dt< \infty$ and $\Iii \|\aa (\Xb_t)\|\, dt< \infty,\ Q\ae$ This last property is satisfied for instance if 
$\aa$ is locally bounded:
\begin{align}\label{eq-50}
\sup _{ t\in\ii, x\in K}\|\aa(t,x)\|< \infty,\quad \textrm{for any compact set $K\subset\ZZ$.}
\end{align}

\begin{definition}[The reference Markov measure $R$]
Let $R\in\MO$ be a Markov measure solution to 
\begin{align}\label{eq-z03}
R\in \MP(\aa,\BBf),
\end{align}
for some vector field $\BBf:\iZ\to\ZZ$. In addition to \eqref{eq-50}, we assume that $R$ fulfils the uniqueness condition:
\begin{align}\label{eq-H17a}
\forall R'\in\MO,\ [R'\in\MP(\aa,\BBf;R_0) \textrm{ and }R'\ll R]\implies R'=R.
\end{align}
\end{definition}

Assumption \eqref{eq-H17a}  is necessary to write explicit formulas for   relative entropies and Radon-Nikodym derivatives  of path measures with respect to $R.$
It is proved in \cite[Thm.\,12.21]{Jac79} that \eqref{eq-H17a} holds if and only if $R$ is an extremal point of the (convex) set of solutions to its own martingale problem: $\MP(\aa,\BBf;R_0)$. To emphasize this uniqueness property, we shall sometimes write
\begin{align*}
R=\MP(\aa,\BBf)
\end{align*}
with an equality. \\
Clearly 
$R(\sbt)=\IZ R ^{ x_o}(\sbt)\,R_0(dx_o),
$ with   $R ^{ x_o}$ the law of the Markov process with initial position $x_o$ and  generator 
\begin{align*}
\partial_t u+ \BBf\scal\nabla u + \Delta _{ \aa}u/2,
\qquad u\in C_c ^{ 1,2}(\iZ),
\end{align*}
see Lemma \ref{res-20} below for a precise statement. 
The measure $R$ is  possibly an unbounded $ \sigma$-finite positive measure. This occurs for instance when $R$ is the  reversible Wiener measure: its reversing measure is  Lebesgue measure, and we take $\BBf=0,$ $\aa=\Id$ and $R_0=\Leb$.
\\
We note for future reference that it is assumed implicitly that 
\begin{align}\label{eq-51}
\Iii|\BBf(\Xb_t)|\,dt< \infty,\quad R\ae
\end{align}

\subsection*{Girsanov theory}
Take $Q\in\PO$ such that
\begin{align}\label{eqd-04}
H(Q|R)< \infty.
\end{align}
We know by the Girsanov theory under a finite entropy  condition  \cite{Leo11a} that there exists some $\ZZ$-valued predictable   process $ \beta ^{ Q|R} $ which is defined $\Qb \ae$    such that  $   \beta ^{ Q|R}  \in \range[\aa(\Xb)],$ $\Qb \ae$,  $Q$ solves the martingale problem
\begin{align}\label{eq-H04}
Q=\MP(\aa,\BBf+\aa \beta ^{ Q|R} ),
\end{align}
and $Q$ inherits the uniqueness property \eqref{eq-H17a} from $R.$
Furthermore, because of  \eqref{eq-H17a}, we know that
\begin{align}\label{eq-H31}
\frac{dQ}{dR}
	&= 	\1 _{ \left\{ dQ/dR>0\right\} }\ 
		\frac{dQ_0}{dR_0}(X_0)
		\exp \left( \Iii  \beta ^{ Q|R} _t\scal dM^R_t
		-\Iii | \beta ^{ Q|R}   _t|^2 _{ \aa(\Xb_t)} /2\ dt\right)\nonumber\\
	&= 	\1 _{ \left\{ dQ/dR>0\right\} }\ 
		\frac{d Q_0}{dR_0}(X_0)
		\exp \left( \Iii  \beta ^{ Q|R} _t\cdot dM^Q_t
		+\Iii | \beta ^{ Q|R}   _t|^2 _{ \aa(\Xb_t)} /2\ dt\right) ,
\end{align}
where 
we denote 
 $ |\beta|^2_\aa:= \beta\cdot \aa\beta,$  
\[
dM^R_t=dX_t-\BBf_t\,dt
\quad \textrm{and}\quad dM^Q=dX_t-(\BBf_t+\aa(\Xb_t)  \beta ^{ Q|R} _t)\,dt.
\]
Here $M^Q$ is a local $Q$-martingale and the local $R$-martingale $M^R$ is seen as a $Q$-local $R$-martingale to define $\Iii  \beta ^{ Q|R}_t\cdot dM^R_t$ as a $Q$-local $R$-stochastic integral, see Definitions \ref{def-z01}.
Moreover, 
\begin{align}\label{eq-H05}
H( Q|R)
	= H( Q_0|R_0)+E_ Q  \Iii  | \beta ^{ Q|R}   _t|^2 _{ \aa(\Xb_t)}/2\ dt.
\end{align}
Of course, the requirement $H( Q |R )< \infty$ implies that 
\begin{align}\label{eq-z01}
E_ Q  \Iii   | \beta ^{ Q|R}   _t|^2 _{ \aa(\Xb_t)}\, dt< \infty.
\end{align} 

\begin{lemma}\label{res-20}
Under the assumptions \eqref{eq-50}, \eqref{eq-H17a} and \eqref{eqd-04}, $C_c ^{ 1,2}(\iZ)$ is included in $\dom\LQ,$ and for any $u\in C_c ^{ 1,2}(\iZ),$
\begin{align*}
\LQ u= (\partial_t +\vv^Q\cdot\nabla + \Delta_\aa/2)u,
\end{align*}
with $\vv^Q:=\BBf+ \aa \beta ^{ Q|R}.$
\\
If in addition $ Q $ is Markov, then the process $  \beta ^{ Q|R}   $ turns out to be a vector field: 
\begin{align}
\label{eq-H25} \beta ^{ Q|R}   _t=\bbf QR  (\Xb_t), \ \bar  Q \ae,
\end{align} 
for some measurable   $\bf ^{ Q|R}:\iZ\to\ZZ$ which is defined $\qb\ae$
\end{lemma}

Note that unlike Lemma \ref{res-23} below, next assumption \eqref{eq-54} is not required for this lemma to hold. 

\begin{proof}
For any  $u\in C_c ^{ 1,2}(\iZ),$ identifying Itô formula and the basic identity attached to $Q\in\MP(\aa,\vv^Q)$  
\begin{align*}
du(\Xb_t)&= \partial_tu(\Xb_t)\,dt+ \nabla u(\Xb_t)\cdot dX_t
	+ \Delta_\aa u(\Xb_t)/2\ dt\\
	&= (\underbrace{ \partial_tu(\Xb_t)+\vv^Q(t,X _{ [0,t]})\cdot \nabla u(\Xb_t)+ \Delta_\aa u(\Xb_t)/2}_{ \textrm{candidate to be } \mathcal{L}^Qu(t,X _{ [0,t]})}) \,dt+dM ^{ u,Q}_t,
\end{align*}
(rely on the uniqueness of the Doob-Meyer decomposition), we see that the increment of the local $Q$-martingale $M ^{ u,Q}$ is 
\begin{align}\label{eq-z10}
dM ^{ u,Q}_t= \nabla u(\Xb_t)\cdot dM^Q_t 
\end{align}
with 
\begin{align*}
dM^Q_t=dX_t-\vv^Q(t,X _{ [0,t]})\,dt,
\end{align*}
the increment of the canonical local $Q$-martingale. By assumption \eqref{eq-50},  $M ^{ u,Q}$ is a square integrable martingale because   $E_Q|M ^{ u,Q}_T|^2=E_Q\Iii |\nabla u|_\aa^2(\Xb_t)\,dt\le T \sup|\nabla u|^2_\aa< \infty.$ 
To prove the first part of the lemma, it remains to verify 
$$
\Iii|( \partial_t+[\BBf(\Xb_t)+\aa(\Xb_t) \beta ^{ Q|R}(t,X _{ [0,t]})]\scal \nabla + \Delta_{\aa(\Xb_t)}/2)u(\Xb _t)|\,dt< \infty,\quad Q\ae
$$
We already know with \eqref{eq-50} that $\sup _{ t,x} |(\partial_t + \Delta_\aa/2)u(t,x)|< \infty,$ and with \eqref{eq-51} and $Q\ll R$ that $\Iii|\BBf(\Xb_t)|\,dt< \infty,\ Q\ae$ On the other hand, by Cauchy-Schwarz inequality, 
\begin{align*}
E_Q\Iii |(\aa \beta ^{ Q|R})\scal \nabla u|^2(t,X _{ [0,t]})\,dt
	&\le E_Q\Iii |\beta ^{ Q|R}|^2_\aa\ | \nabla u|^2_\aa(t,X _{ [0,t]})\,dt\\
	&\le \sup | \nabla u|^2_\aa\  E_Q\Iii |\beta ^{ Q|R}|^2_\aa\ (t,X _{ [0,t]})\,dt
	< \infty,
\end{align*}
where we use assumption \eqref{eq-50} to control $ \sup | \nabla u|_\aa$, and assumption \eqref{eqd-04} to obtain \eqref{eq-z01}.
Therefore, $\Iii |(\aa \beta ^{ Q|R})\scal \nabla u|(t,X _{ [0,t]})\,dt< \infty,\ Q\ae,$
and the proof of the first statement is done.

Now, suppose that in addition $Q$ is Markov. By Corollary \ref{res-18},  $\LQ u(t,X _{ [0,t]})=v(\Xb_t),\ \Qb\ae$ for some function $v:\iZ\to\RR.$ Hence, for any $u\in C_c ^{ 1,2}(\iZ),$ $\aa(\Xb_t) \beta ^{ Q|R}_t\scal \nabla u(\Xb_t)=v(\Xb_t)- \partial_tu(\Xb_t)-\BBf\scal\nabla u(\Xb_t)- \Delta _{ \aa(\Xb_t)}u(\Xb_t)/2,$ $\Qb\ae$ This implies that  $\aa_{\Xb} \beta ^{ Q|R}= \vv(\Xb),\ \Qb\ae,$ for some vector field $\vv:\iZ\to \ZZ.$ Of course, $\vv\in\range(\aa), $ $\qb\ae$ As $ \beta ^{ Q|R}$ also satisfies $ \beta ^{ Q|R}\in\range(\aa _{ \Xb})$, without assuming that $\aa$ is invertible, one can solve uniquely  this  equation to obtain $ \beta ^{ Q|R}=\aa ^{ -1}\vv(\Xb),$ $\Qb\ae$ where $\aa ^{ -1}$ is the generalized inverse of $\aa$ defined by $\Lim n (\aa + \Id/n) ^{ -1}.$ We conclude  taking $\bbf QR=\aa ^{ -1}\vv.$
\end{proof}

\subsection*{Kinetic action}

We observe that 
\begin{align}\label{eq-46}
H(Q|R)-H(Q_0|R_0)
	=H(Q|R ^{ Q_0})
	=E_Q\Iii\ud |\vv^{ Q|R}_t|^2_{\gg (\Xb_t)}\, dt
\end{align}
is an average kinetic action, 
where $$\vv^{ Q|R}:=\aa  \beta ^{ Q|R} $$ should be interpreted as a stochastic  relative velocity between $Q$ and $R,$ and 
\begin{align*}
|\vv|^2_\gg  := \left\{ \begin{array}{ll}
\vv\cdot \aa ^{ -1}\vv,\quad & \textrm{if }\vv\in \mathrm{range}(\aa),\\
+ \infty,& \textrm{otherwise.}
\end{array}\right. 
\end{align*}

\subsection*{Martingale representation}

The martingale representation theorem will play an important role. 
We assume that $R$ is the law of a diffusion process solution to the stochastic differential equation
\begin{align*}
d\yy_t=\BBf(t,\yy_t)\, dt + \sigma(t,\yy_t)\, d\ww_t
\end{align*}
where $\ww$ is an $n$-dimensional Brownian motion built on some unspecified filtered triple    $(\Xi, ( \mathcal{F}_t) _{ 0\le t\le T}, \mathcal{F}, \PP) $ where the $n\times n$-matrix field  $ \sigma$ satisfies $ \sigma \sigma ^{ \mathsf{t}}=\aa.$ In this setting, the martingale representation theorem states that if 
\begin{align}\label{eq-54}
\aa(t,x) \textrm{ is invertible for all } (t,x)\in\iZ,
\end{align}
for any local $ (\mathcal{F} ^{ \yy},\PP)$-martingale $ \mathsf{M}$,  where $ \mathcal{F} ^{ \yy}:= ( \sigma(\yy _{ [0,t]})) _{ 0\le t\le T}$ is the $\PP$-complete natural filtration of the process $\yy,$  there exists a predictable process $\Phi ^{ \mathsf{M}}$ such that $ \sigma \Phi ^{ \mathsf{M}}$ is locally square integrable  and
\begin{align*}
 \mathsf{M}_t= \mathsf{M}_0+\int_0^t \Phi ^{ \mathsf{M}}_s\cdot \sigma(s,\yy_s)\,d\ww_s,\quad 0\le t\le T,\quad \PP\ae
\end{align*}
Now, let us go back to the canonical setting by taking the image by 
$\yy:\Xi\to \OO$ of $(\Xi, \mathcal{F} ^{ \yy}, \PP)$ to obtain $R=\yy\pf\PP$ and some sub-filtration of the canonical one. After completing it by the $R$-negligible sets, this filtration coincides with the $R$-completion of the canonical filtration. We also see that under the law $R,$ $dM^R_t:=dX_t-\BBf(\Xb_t)\,dt\eqlaw \sigma(t,\yy_t)\, d\ww_t.$
More generally, Girsanov's theory also tells us that for any $Q$ such that $H(Q|R)< \infty,$ under the law $Q$ we have $dM^Q_t:=dX_t-\vv^Q(t,X _{ [0,t]})\,dt\eqlaw \sigma(t,\yy_t)\, d\ww_t.$ Otherwise stated,
\begin{align}\label{eq-57}
dM^Q_t:=dX_t-\vv^Q(t,X _{ [0,t]})\,dt
	= \sigma(\Xb_t)\,dW^Q_t,\quad \Qb\ae,
\end{align}
where $W^Q$ is a $Q$-Brownian motion on the canonical space equipped with the $Q$-completion of the canonical filtration. Moreover, for any function $u$ in $\dom \LLF^Q$, by the martingale representation theorem which holds because of  the assumed invertibility of $ \aa$,
\begin{align}\label{eq-52}
dM ^{ u,Q}_t:=du(\Xb_t)-\LLF^Qu(t,X _{ [0,t]})\,dt
	= \alpha ^{ u,Q}_t \cdot  dM ^{ Q}_t,\quad \Qb\ae,
\end{align}
for some predictable process $ \alpha ^{ u,Q}.$ This implies in particular
\begin{align}\label{eq-52b}
d\langle u(\Xb),v(\Xb)\rangle^Q_t
	= \alpha ^{ u,Q}_t\scal \aa(\Xb_t) \alpha ^{ v,Q}_t\, dt.
\end{align}

\begin{lemma}\label{res-23}
We assume  \eqref{eq-54}, i.e.\  $\aa$ is invertible.
\begin{enumerate}[(a)]
\item
Then $\dom\LLF^Q$  is an algebra,
meaning that for any $u,v\in \dom\LLF^Q,$ the product $uv$ is still in $\dom\LLF^Q.$

\item
For any $u\in\dom\LLF^Q$ and any  function $F\in C ^{2}(\RR),$ $F(u)$ is also in $\dom\LLF^Q.$
\end{enumerate}

\end{lemma}

\begin{proof}
Under the  ellipticity assumption \eqref{eq-54}, we have  \eqref{eq-52b}. This is the key  of the proof.\\ 
\boulette{(a)}
As a definition of the forward generator
\begin{align*}
du(\Xb_t)=\LLF^Q_tu(t, X _{ [0,t]})\,dt + dM^u_t,
\qquad dv(\Xb_t)=\LLF^Q_t v(t, X _{ [0,t]})\, dt+dM^v_t,
\end{align*}
and  applying Itô formula in the forward sense of time \begin{align*}
d(uv)(\Xb_t)&=u(\Xb_t)dv(\Xb_t)+v(\Xb_t)du(\Xb_t)+d[u(\Xb),v(\Xb)]_t\\
	&=u(\Xb_t)dv(\Xb_t)+v(\Xb_t)du(\Xb_t)+d\langle u(\Xb),v(\Xb)\rangle_t\\
	&=[u\LLF^Q_tv(t, X _{ [0,t]})+v\LLF^Q _tu(t, X _{ [0,t]})]\,dt+d\langle u(\Xb),v(\Xb)\rangle_t\\
	&\hskip 4cm +u(\Xb_t)dM^v_t+v(\Xb_t)dM^u_t.
\end{align*}
The bounded variation part of this semimartingale is
\begin{align}\label{eq-56}
[u\LLF^Q_tv(t, X _{ [0,t]})+v\LLF^Q _tu(t, X _{ [0,t]})]\,dt+d\langle u(\Xb),v(\Xb)\rangle^Q_t.
\end{align}
With \eqref{eq-52b},  this shows that $uv(\Xb)$ is a \emph{nice} $Q$-semimartingale, which is the announced result.

\Boulette{(b)}
By Itô formula and  \eqref{eq-52b}:
\begin{align*}
dF(u)(\Xb_t)
	&= F' (u(\Xb_t)) du(\Xb_t)+  F''(u(\Xb_t))d\langle u(X),u(X)\rangle_t/2\\
	&= [F' (u(\Xb_t)) \LLF^Q_tu(t, X _{ [0,t]})+F''(u(\Xb_t)) \alpha ^{ u,Q}_t\scal \aa(\Xb_t) \alpha ^{ u,Q}_t/2]\, dt+ F' (u(\Xb_t)) \, dM ^u_t,
\end{align*}
we see that the bounded variation part of $F(u(\Xb))$ is absolutely continuous, which is the announced result.
\end{proof}

This  simple consequence of the martingale representation theorem  is  crucial for the remainder of the article.

\subsection*{Extended Itô formula}

We already saw at	\eqref{eq-z10} while proving Lemma \ref{res-20} that for any regular function $u\in C ^{ 1,2}_c(\iZ),$
$
dM ^{ u,Q}_t= \nabla u(\Xb_t)\cdot dM^Q_t .
$
Let us look at the extension of this identity to the case where $u$ belongs to $\dom \mathcal{L}^Q$.

\begin{proposition}[Extended Itô formula]\ \label{res-14}
We assume that $R$ satisfies \eqref{eq-50}, \eqref{eq-H17a},  \eqref{eq-54}, and that $Q$ is Markov and satisfies  \eqref{eqd-04}. \\ Remember that by Lemma \ref{res-20}, $Q=\MP(\aa,\vf^Q)$ with
\begin{align}\label{eq-44}
\vf^Q=\BBf+\aa\bbf QR.
\end{align}
 The following statements are verified.
\begin{enumerate}[(a)]
\item
For any $u$ in $\dom \mathcal{L}^Q,$ there exists a $\qb\ae$ defined vector field $\gradt^Qu$ such that
\begin{align}\label{eq-z04}
\begin{split}
du(\Xb_t)
	= \mathcal{L}^Qu(\Xb_t)\, dt +&\gradt^Qu(\Xb_t)\cdot dM^Q_t\\
	&=[\mathcal{L}^Q-\vf^Q\cdot \gradt^Q]u(\Xb_t)\, dt
			+\gradt^Qu(\Xb_t)\cdot dX_t,
\quad \Qb\ae
\end{split}
\end{align}

\item
Let $P\ll Q$. For any $u$ in $\dom \LLL QP$ (recall Definition \ref{defh-05b}), there exists a $\pb\ae$-defined vector field $\Gt QP u$ such that
\begin{align}\label{eq-z04b}
\begin{split}
du(\Xb_t)
	&= \LLL QPu(\Xb_t)\, dt +\Gt QPu(\Xb_t)\cdot dM ^{ Q}_t\\
	&=[\LLL QP-\vf^Q\cdot \Gt QP]u(\Xb_t)\, dt
			+\Gt QP u(\Xb_t)\cdot dX_t,
\quad \Pb\ae,
\end{split}
\end{align}
where   $\Gt QPu(\Xb_t)\cdot dM ^{ Q}_t$ is the increment of a $P$-local $Q$-stochastic integral (recall Definition \ref{def-z01}).
\end{enumerate}
\end{proposition}

\begin{remark}
Note that $\Gt QP u$ is only defined $\pb\ae$ (and not $\qb\ae$), and unlike  \eqref{eq-z04}, the identity  \eqref{eq-z04b} only holds $\Pb\ae$, but is meaningless $\Qb\ae$ in general.
\end{remark}

\begin{proof}
The proof of statement (b) is an almost verbatim modification of the proof of statement (a). It is left to  the reader.\\
 In what follows, the identities hold $\Qb\ae$
Let $u\in\dom \mathcal{L}^Q.$
By the very  definition of the extended generator:
\begin{align}\label{eq-z06}
du(\Xb_t)= \mathcal{L}^Qu(\Xb_t)\,dt+dM ^{ u,Q}_t,
\end{align}
with $M ^{ u,Q}$ a local $Q$-martingale. And by \eqref{eq-52}: 
\[dM ^{ u,Q}_t=\alpha ^{ u,Q}_t\cdot dM^Q_t\]
for some predictable process $\alpha ^{ u,Q}$.
Taking next Lemma \ref{res-21} for granted, it follows that
$	
dM ^{ u,Q}_t=\gradt^Q u(\Xb_t)\cdot dM^Q_t,
$	
and we complete the proof   with \eqref{eq-z06}.
\end{proof}

It remains to prove
\begin{lemma}\label{res-21}
Under the ellipticity condition \eqref{eq-54}, 
the  Markov measure  $Q$ is  such that the process $ \alpha ^{ u,Q}$ only depends on the current position: $ \alpha ^{ u,Q}_t=\gradt^Q u(\Xb_t),$ $0\le t\le T,$ for some $\qb$-almost everywhere defined vector field $\gradt^Q u:\iZ\to\ZZ.$ 
\end{lemma}

\begin{proof}
We have  just seen that 
\begin{align*}
dM ^{ u,Q}_t=\alpha ^{ u,Q}_t\cdot dM^Q_t
	=du(\Xb_t)- \mathcal{L}^Qu(\Xb_t)\,dt,
\end{align*}
and for any  function $v$ in $C_c ^{ 1,2}(\iZ)$, Itô formula is
\begin{align*}
dv(\Xb_t)= \partial_tv(\Xb_t)\,dt+ \nabla v(\Xb_t)\cdot dX_t
	+ \Delta_\aa  v(\Xb_t)/2\ dt.
\end{align*}
Hence, the quadratic covariation of $u(\Xb)$ and $v(\Xb)$ satisfies
\begin{align*}
d \left\langle u(\Xb),v(\Xb) \right\rangle _t
	= \alpha ^{ u,Q}_t\cdot \aa \nabla v(\Xb_t)\, dt.
\end{align*}
Therefore, with intuitive arguments we see that
\begin{align*}
\alpha ^{ u,Q}_t\cdot \aa \nabla v(\Xb_t)\, dt
	&= E_Q\left.\big(du(\Xb_t)dv(\Xb_t)\,\right|\,  X _{[0,t] }\big)
	=  E_Q\left.\big(du(\Xb_t)dv(\Xb_t)\,\right|\,  X _{t }\big)\\
	&= E_Q\big(\alpha ^{ u,Q}_t\mid  X _{t}\big)\cdot \aa \nabla v(\Xb_t)\,dt,
\end{align*}
where we used the Markov property of $Q$ at second equality. Since $v$ is arbitrary and $\aa$ is invertible, we see that $ \alpha_t ^{ u,Q}=E_Q\big( \alpha ^{ u,Q}_t\mid X _{t}\big),$ proving the lemma.

However, it is necessary to justify rigorously  the above string of identities. Firstly, one must  localize (by means of stopping times) to give some meaning at the conditional expectation $E_Q\big( \alpha ^{ u,Q}_t\mid X _{t}\big),$ and  the above identity
\begin{align}\label{eq-53}
d \left\langle u(\Xb),v(\Xb) \right\rangle _t= E_Q\big(du(\Xb_t)dv(\Xb_t)\mid X _{[0,t] }\big),\quad \Qb\ae,
\end{align}
still needs to be carefully established, again with some localization argument.
\\
The first localization argument is standard, so we leave it to the reader. 
\\
Finally,    \eqref{eq-53} follows from Lemma \ref{res-22}  in the appendix. To verify that the hypotheses of this lemma are satisfied note that in the present Brownian setting  under the hypothesis \eqref{eq-54}, by Lemma \ref{res-23} we know that $\dom\LLF^Q$ is an algebra. On the other hand the martingale representation theorem implies that   any martingale is continuous. Hence the  assumption about $M ^{ Q,[u,v]}$ in Lemma \ref{res-22} is trivially satisfied because $M ^{ Q,[u,v]}=0.$
\end{proof}

\begin{remarks}\ \begin{enumerate}[(i)] \label{rem-03}
\item
By \eqref{eq-z10}, if $u$ is $C ^{ 1,2}$-regular, then $\gradt^Qu=\nabla u,\ \qb\ae$

\item
It is proved by Cont and Fournié in \cite[Thm.\,5.9]{CF13} that under   hypotheses slightly more general than those of Proposition \ref{res-14} (no entropy appears), any square integrable $Q$-martingale $N$ is represented as the stochastic integral
\begin{align*}
N_t=N_0+\int_0^t \nabla _{ M^Q} N_s\cdot dM^Q_s
\end{align*}
where $M^Q$ appears at \eqref{eq-57} and the \emph{Cont-Fournié non-anticipative derivative} $ \nabla _{ M^Q} N$ is   introduced in \cite{CF13}.  Therefore,
\begin{align*}
\gradt^Qu(\Xb)=\nabla _{ M^Q}\,M ^{ Q,u},\qquad \Qb\ae
\end{align*}
 The Cont-Fournié derivative is the predictable projection of the Malliavin (anticipative) derivative.
 
 \item
 A remarkable extension of Itô formula is also obtained in \cite{CF13}, which goes in another direction than Proposition \ref{res-14}: stochastic differentials of regular non-anticipative functionals are considered in \cite{CF13}, while Proposition \ref{res-14} gives a result for possibly not regular functions only depending on the current position.
\end{enumerate}\end{remarks}

\section{Feynman-Kac formula}
\label{sec-HJB}

The main character of this section is the Feynman-Kac measure $P$ already encountered in the introduction at \eqref{eq-77}. It will allow us to derive pathwise properties of an extended \eqref{eq-HJB} equation at Theorem \ref{res-13} and an extended \eqref{eq-FK} equation at Theorem \ref{res-31}.

\subsection*{The Feynman-Kac  measure}
It  is defined by
\begin{align}\label{eq-H03}
P:=f_0(X_0)\exp \left(\Iii V(\Xb_t)\,dt\right) g_T(X_T)\ R\in\PO,
\end{align}
where as in \eqref{eq-z03}  $$R=\MP(\aa,\BBf)$$ is a reference diffusion measure satisfying \eqref{eq-H17a}, $f_0$ and $g_T$ are nonnegative measurable functions on $\ZZ,$ $V:\iZ\to \RR\cup \left\{- \infty\right\} ,$ and we  use the convention $e ^{ - \infty}=0.$ It is also assumed that all these quantities are such that $P$ is a probability measure and
\begin{align*}
H(P|R)< \infty.
\end{align*}
For $P$ to be a probability measure it is necessary that the measurable subsets
\begin{align*}
D_+&:= \left\{ (x,y)\in\ZZZ; R\Big(\Iii V_+(\Xb_t)\,dt<\infty\mid X_0=x;X_T=y\Big)>0\right\} ,\\
 D_-&:= \left\{ (x,y)\in\ZZZ; R\Big(\Iii V_-(\Xb_t)\,dt< \infty\mid X_0=x;X_T=y\Big)>0\right\} ,\\
 D&:= \left\{ (x,y)\in\ZZZ; f_0(x)>0, g_T(y)>0\right\} ,
\end{align*}
verify $R _{ 01}\ae$
\begin{align}\label{eq-87}
(a):\ D_+\cup D_-=\ZZZ,\qquad
(b):\ D\subset D_+,\qquad 
(c):\ D \cap D_-\neq\emptyset.
\end{align}
Indeed, (a), (b) and (c) are implied respectively by the properties of  $dP/dR$ of being well defined,   finite, and  not identically equal to zero, up to some $R _{ 0T}$-negligible set. We denote $R _{ 0T}:=(X_0,X_T)\pf R$ the law of the endpoint position under $R.$

Sufficient conditions for \eqref{eq-87} and $H(P|R)< \infty$ will be given at Section \ref{sec-GC}. From now on, when assuming that $P$ is in $\PO$ it is supposed implicitly that \eqref{eq-87} holds.

The measure $P$  will bring us valuable informations about equation \eqref{eq-FK}.

\subsection*{Positive integration}

For all $0\le t\le T$ and $ P_t$-almost every $x\in\ZZ,$ define the measurable functions $f,g: \ii\times\ZZ\to\RR$ by
\begin{align}\label{eq-z14}
\begin{split}
f(t,x):=f_t(x)&:=E_R\left.\left(f_0(X_0)\exp \Big(\int _{ [0,t]} V(\Xb_s)\,ds\Big) \,\right|\, X_t=x\right),\\
g(t,x):=g_t(x)&:=E_R\left.\left(\exp \Big(\II tT V(\Xb_s)\,ds\Big)  g_T(X_T)\,\right|\,  X_t=x\right).
\end{split}
\end{align}
One must be careful with these conditional expectations because it is not assumed that the integrands are $R$-integrable, neither that $R$ is bounded. However,  they are well defined $\pb$-almost everywhere (but not $\rb\ae$) as conditional expectations of nonnegative functions with respect to a possibly unbounded Markov measure. This is warranted by next

\begin{lemma}[{\cite[\S 4]{Leo12b}}] \label{res-15}
Let   $R\in\MO$  be a Markov measure and 
     $P\in\PO$  a  probability measure such that $P\ll R$ and  $\displaystyle{\frac{dP}{dR}=\alpha \zeta 
\beta }$ with $\alpha,\zeta, \beta$ \emph{nonnegative} functions such that $\alpha \in\sigma(X_{[0,s]})$, $\zeta\in \sigma(X_{[s,t]})$ and $\beta \in\sigma(X_{[t, T]})$ for some
$0\le s\le t\le T$. Then,
$$
\left\{
\begin{array}{l}
  E_R(\alpha \mid X_s),E_R(\beta \mid X_t)\in (0,\infty) \\
  E_R(\alpha \beta \mid X_{[s,t]})=E_R(\alpha \mid X_s)E_R(\beta \mid
    X_t)\in (0,\infty)\\
\end{array}
    \right.\quad P\ae
$$
(and not  $R\ae$ in general). In addition,
\begin{align*}
 \frac{dP _{ [s,t]}}{dR _{ [s,t]}}(X _{ [s,t]})
    = E_R(\alpha \mid X_s)\zeta E_R(\beta \mid X_t)\in (0,\infty),\qquad P\ae
\end{align*}
 \end{lemma}

\begin{remark}
Even if the product $\alpha
\beta $ is integrable, it is not true in general that the nonnegative factors $\alpha $ and $\beta $ are integrable. Therefore, a priori the conditional expectations  $E_R(\alpha \mid
X_s)$ and $E_R(\beta \mid X_t)$ might have been infinite.
\end{remark}

\begin{proposition}\label{res-16}
The  Feynman-Kac measure $P$ given at \eqref{eq-H03} is Markov and
\begin{align}\label{eq-H34}
 \frac{dP_t}{dR_t}=f_tg_t,
 \quad P_t\ae,\quad \forall 0\le t\le T.
\end{align}
\end{proposition}

\begin{proof}
For any  $0\le t\le T$,   the Radon-Nikodym derivative  $Z:= dP/dR$ equals $Z=Z _{ [0,t]} Z _{ [t,T]}$ with $Z _{ [0,t]}:=f_0(X_0)\exp \Big(\int _{ [0,t]} V(\Xb_r)\,dr\Big),$ $Z _{ [t,T]}:=\exp \Big(\II tT V(\Xb_r)\,dr\Big)  g_T(X_T)$.

The identity \eqref{eq-H34} is a direct consequence of Lemma \ref{res-15} applied with $s=t,$ $\zeta=1,$ $\alpha=Z _{ [0,t]}$ and $ \beta=Z _{ [t,T]}$.

The Markov property is proved  by showing that  for any bounded measurable functions $u(X _{ [0,t]})$ and $v(X _{[t,T]})$, we obtain
\begin{align*}
E_P[u(X _{ [0,t]})v(X _{[t,T]})\mid X_t]
	&= \frac{E_R[Z\,u(X _{ [0,t]})v(X _{[t,T]})\mid X_t]}{E_R[Z\mid X_t]}\\
	&=\frac{E_R[Z _{ [0,t]}u(X _{ [0,t]})Z _{ [t,T]}v(X _{[t,T]})\mid X_t]}{E_R[Z _{ [0,t]}Z _{ [t,T]}\mid X_t]}\\
	&=\frac{E_R[Z _{ [0,t]}u(X _{ [0,t]})\mid X_t]}{E_R[Z _{ [0,t]}\mid X_t]}\ 
	\frac{E_R[Z _{ [t,T]}v(X _{[t,T]})\mid X_t]}{E_R[Z _{ [t,T]}\mid X_t]}\\
	&= E_P[u(X _{ [0,t]})\mid X_t]\, E_P[v(X _{[t,T]})\mid X_t].
\end{align*}
We invoked the Markov property of $R$ and Lemma \ref{res-15} at the last but one equality. \end{proof}

\subsection*{Extended HJB equation}

We also introduce  the measure $\pb (dtdx)=dt P_t(dx)$ on $\iZ$ and  the notation
\begin{align}\label{eq-47}
 \varphi:=\log f,\qquad \psi:=\log g
\end{align}
where $f$ and $g$ are defined at \eqref{eq-z14}. 
 We are  ready to state the  first main result of this section.
 
\begin{theorem}\label{res-13}
Let $R\in\MO$ be a Markov measure  satisfying \eqref{eq-z03},  \eqref{eq-H17a} with $\aa$ verifying \eqref{eq-50}, \eqref{eq-54}, and let $P\in\PO$ be given by \eqref{eq-H03}.   Then, under the assumption that \[- \infty<H(P|R)< \infty,\]  the following statements hold.

\begin{enumerate}[(a)]
\item
The $\pb\ae$\,defined function 
\begin{align*}
 \psi(t,x):=\log g_t(x)=\log E_R\left.\Big(\exp \Big(\int _{ [t,T]} V(\Xb_s)\,ds\Big)  g_T(X_T)\,\right|\,  X_t=x\Big)\in\RR,
 \end{align*}
 is in $\dom \LLL RP$ and it satisfies the extended HJB equation
\begin{align}\label{eq-H38}
\big(\LLL RP \psi+| \Gt RP \psi|^2_\aa /2 +V\big)(\Xb)=0,
\qquad \Pb\ae
\end{align}
with $ \psi_T=\log g_T,$  $P_T\ae$ Remember that the existence of $\Gt RP \psi$ and its definition are stated  at Proposition \ref{res-14} (extended Itô formula).

\item
The Feynman-Kac measure $P$ solves
$\MP(\aa,\vf^P)$ where
\begin{align}\label{eq-z09}
\vf^P=\BBf+\aa \Gt RP \psi.
\end{align}
\item
In addition, $ \psi\in\dom \mathcal{L}^P,$ 
\begin{align}\label{eq-z08}
\mathcal{L}^P \psi(\Xb)=\big(| \gradt^P \psi|^2_\aa /2 -V\big)(\Xb),\quad \Pb\ae
\end{align}
and
\begin{align}\label{eq-z07}
\begin{split}
\LLL RP\psi(\Xb)&= \big(\mathcal{L}^P \psi - |\gradt^P\psi|^2_\aa \big)(\Xb),
\qquad\Pb\ae
\\
 \Gt RP \psi(\Xb)&=\gradt^P \psi(\Xb),
\qquad\Pb\ae
\end{split}
\end{align}
\end{enumerate}
\end{theorem}
Before proving this theorem, we make some remarks and establish  preliminary estimates at Lemma \ref{res-H02} under a finite entropy condition.

\begin{remarks}\ \begin{enumerate}[(a)]
\item
It  follows  from \eqref{eq-H38} that
\begin{align}\label{eq-H39}
\LLL RP \psi+| \Gt RP \psi|^2_\aa /2 +V=0,
\qquad \pb\ae
\end{align}
and it turns out that when  $ \psi$ is a finite $C ^{ 1,2}$ function, \eqref{eq-H39} is the standard HJB equation \eqref{eq-HJB}:
\begin{align*}
(\partial_t+ \BBf\cdot \nabla+ \Delta _{ \aa}/2) \psi+ |\nabla \psi|^2_\aa  /2+V=0.
\end{align*}
\item
By \eqref{eq-H34}, $f_tg_t>0,\ P_t\ae$ Consequently, for all $t$, $ \psi_t:=\log g_t$ is well defined $ P_t\ae$ as a real valued function.
\end{enumerate}\end{remarks}

For any $0\le s<t\le T,$ the restriction $P _{ [s,t]}$ of $P$ to  $ \sigma(X _{ [s,t]})$ satisfies 
\begin{align}\label{eq-H35}
\frac{dP _{ [s,t]}}{dR _{ [s,t]}}
	= \frac{dP_s}{dR_s}(X_s)  \exp \left( \psi(\Xb_t)- \psi(\Xb_s)+\int_{[s,t]} V(\Xb_r)\,dr\right),
 \qquad P\ae
\end{align}
This formula will be used in a while. 
\\
 Let us prove it. The following identities hold $P\ae$ (and not necessarily $R\ae$) without further mention:
\begin{align*}
&\frac{dP _{ [s,t]}}{dR _{ [s,t]}}
=E_R \left.\left[ \frac{dP}{dR}\,\right|\,  X _{ [s,t]}\right]
	 \overset{(i)}= E_R \left.\left[f_0(X_0)\exp \left(\Iii V(\Xb_r)\,dr\right) g_T(X_T)\,\right|\,  X _{ [s,t]}\right]\\
	=\ &\exp \left( \int _{ [s,t]} V(\Xb_r)\,dr\right) 
	 E_R \Big[f_0(X_0)\exp \left(\int _{ [0,s]} V(\Xb_r)\,dr\right)\\
	&\hskip 8cm \times  \exp \left(\int _{ (t,T]} V(\Xb_r)\,dr\right) g_T(X_T)\,\Big|\,  X _{ [s,t]}\Big]\\
	 \overset{(ii)}=\ &\exp \left( \int _{ [s,t]} V(\Xb_r)\,dr\right) 
	 E_R \Big[f_0(X_0)\exp \left(\int _{ [0,s]} V(\Xb_r)\,dr\right)\,\Big|\,  X _{ [s,t]}\Big]\\
	&\hskip 7cm \times   E_R \Big[\exp \left(\int _{ (t,T]} V(\Xb_r)\,dr\right) g_T(X_T)\,\Big|\,  X _{ [s,t]}\Big]\\
			 \overset{(iii)}=\ & f_s(X_s) \exp \left( \int _{ [s,t]} V(\Xb_r)\,dr\right) g_t(X_t)
	=f_s(X_s) g_s(X_s)\ \frac{ g_t(X_t)}{g_s(X_s)} \exp \left( \int _{ [s,t]} V(\Xb_r)\,dr\right)\\
	 \overset{(iv)}=\ &\frac{dP_s}{dR_s}(X_s)  \exp \left( \psi(\Xb_t)- \psi(\Xb_s)+\int_{[s,t]} V(\Xb_r)\,dr\right)
\end{align*}
where we used \eqref{eq-H03} at (i),  the Markov property of $R$  and  Lemma \ref{res-15} at (ii), the definitions \eqref{eq-z14} of $f_s$ and $g_t$ at (iii), and \eqref{eq-H34} and the definition $ \psi:=\log g$ at (iv). Recall  the time-symmetric formulation of the Markov property:  a path measure $R$ is Markov  if and only if  for any $0\le s\le t\le T,$ $X _{ [0,s]}$ and $X _{ [t,T]}$ are $R(\sbt\mid X _{ [s,t]}) $-independent.

Remark that no division by zero occurs $P\ae$ because $g_s(X_s)$ is positive $P\ae$ by \eqref{eq-H34}.

\begin{lemma}[Finite entropy estimates]\label{res-H02}
Suppose  that  $- \infty<H(P|R)< \infty$.  Then,  
\begin{enumerate}[\quad (i)]
\item
$ \varphi_t, \psi_t\in L^1( P_t),$ for all $0\le t\le T$;
\item
$V(\Xb)\in L^1(\Pb);$\quad
$V_t\in L^1( P_t)$ for almost every $0\le t\le T;$
\item
 $t\mapsto \langle V_t,P_t\rangle :=\IZ V_t\, dP_t $ is $dt$-integrable;
\item
$\psi(\Xb_t)-\psi(\Xb_s)+\int_{[s,t]} V(\Xb_r)\,dr\in L^1(P )$ for all $0\le s\le t\le T$;
\end{enumerate}
 and
\begin{align*}
H(P|R)
	= \langle \varphi_0,P_0\rangle +\Iii \langle V_t, P_t\rangle\,dt+\langle\psi_T,P_T\rangle .
\end{align*}
\end{lemma}

\begin{proof}
For any $0\le s\le t\le T,$ $H(P_s|R_s)\le H(P _{ [s,t]}| R _{ [s,t]})\le H(P|R)< \infty$.
On the other hand, we obtain with \eqref{eq-H34} and \eqref{eq-H35} that 
\begin{align*}
H(P_s|R_s)
	&= \IZ\log(f_sg_s)\, dP_s=\IZ( \varphi_s+ \psi_s)\,dP_s,\\
H(P _{ [s,t]}| R_{[s,t] })	
	&= H(P_s|R_s)
		+ E_P \Big(\psi(\Xb_t)-\psi(\Xb_s)+\int_{[s,t]} V(\Xb_r)\,dr\Big).\end{align*}
Hence (see Remark \ref{rem-02} and Proposition \ref{res-24}),
$ \log(f_s g_s)\in L^1(P_s),$ and  
$\psi(\Xb_t)-\psi(\Xb_s)+\int_{[s,t]} V(\Xb_r)\,dr\in L^1(P)$.
We conclude with Fubini-Lebesgue theorem applied with the product measure $( \delta_s+\Leb _{ [s,t]}+ \delta_t)\otimes P _{ [s,t]}$ on $[s,t]\times\OO _{ [s,t]}.$ 
\end{proof}

\begin{proof}[Proof of Theorem \ref{res-13}]

The boundary condition  $ \psi_T=\log g_T$ is an obvious outcome of the expression of $ \psi.$
\\
By  \eqref{eq-H31}, \eqref{eq-H05} and the Markov property of $P$ proved at Proposition \ref{res-16} and Lemma \ref{res-20}, there exists a vector field $ \gamma: \iZ\to \ZZ$  such that
\begin{align}\label{eq-48}
\IiZ | \gamma|^2_\aa\,d\pb< \infty
\end{align}
and for any $0\le s<t\le T,$
\begin{align}\label{eq-H32}
\frac{dP _{ [s,t]}}{dR _{ [s,t]}}
	=\frac{dP_s}{dR_s}(X_s)
		\exp \left( \int_{[s,t]} \gamma(\Xb_r)\cdot dM^P_r
			+\int_{[s,t]} | \gamma|^2_\aa(\Xb_r)/2\ dr\right),	
		\quad P\ae,
\end{align}
where   $dM^P_r=dX_r-(\BBf+ \aa \gamma)(\Xb_r)\,dr$.
\\
By Lemma \ref{res-H02}, we have: $V(\Xb)\in L^1(\Pb)$, and because of \eqref{eq-48} we have also: $| \gamma|^2_\aa(\Xb)\in L^1(\Pb).$ Hence, applying  Lemma \ref{resh-21} with $k^h=\frac 1h \1_{[-h,0]}$, we obtain that
\begin{align*}
\lim _{ h\to 0^+}\Big[t\mapsto h ^{ -1}\int_{ [t,t+h]}  ( | \gamma|^2_\aa /2- V)(\Xb_r)\,dr \Big]
	=(| \gamma|^2_\aa /2- V)(\Xb)
\quad \textrm{in }L^1(\Pb).
\end{align*}
Identifying \eqref{eq-H35} with \eqref{eq-H32} gives us
 for all $0\le s<t\le T,$
\begin{equation}\label{eq-H33}
\psi(\Xb_t)- \psi(\Xb_s)
	= \int_{[s,t]} ( | \gamma|^2_\aa /2- V)(\Xb_r)\,dr 
		+ \int_{[s,t]}  \gamma(\Xb_r)\cdot dM^P_r,
	\quad P\ae,
\end{equation}
which implies that for any $0\le t< t+h\le T,$
\begin{align*}
E_P\left. \left[h ^{ -1}\big( \psi(\Xb _{ t+h})- \psi(\Xb_t)\big) \,\right|\,  X_t \right] 
	= E_P\left. \left[h ^{ -1}\int_{[t,t+h]} ( | \gamma|^2_\aa /2- V)(\Xb_r)\,dr  \,\right|\,  X_t \right].
\end{align*}
By \eqref{eq-H33} again, $ \psi(\Xb)$ admits a continuous version because 
$t\mapsto \int_{[0,t]} ( | \gamma|^2_\aa /2-V)(\Xb_r)\,dr$ is continuous $P\ae$, since $( | \gamma|^2_\aa /2-V)(\Xb)$ is $\Pb$-integrable, and $t\mapsto  \int_{[0,t]}  \gamma(\Xb_r)\cdot dM^P_r$ admits a continuous version as a Brownian  stochastic integral. Therefore, we are in position to apply
Proposition \ref{resh-03} which ensures that $ \psi$ belongs to  $\dom \mathcal{L}^P$, and 
\begin{align}\label{eq-H36}
\mathcal{L}^P \psi(\Xb)=[| \gamma|^2_\aa /2 -V](\Xb),\quad \Pb\ae
\end{align}
As $ \psi$ is in $\dom \mathcal{L}^P$, the extended Itô formula \eqref{eq-z04} is valid and writes as
 \begin{align}
d\psi(\Xb_t)\label{eq-H33b}
	=\mathcal{L}^P \psi(\Xb_t)\, dt
		+ \gradt^P \psi(\Xb_t)\cdot dM^P_t,
		\qquad \Pb\ae,
\end{align}
for some  vector field $\gradt^P \psi$.
Note that identifying the bounded variation parts of \eqref{eq-H33} and \eqref{eq-H33b} also leaves us with \eqref{eq-H36}, while the identification of the martingale parts  yields $| \gamma-\gradt^P \psi|^2_\aa(\Xb)=0,\ \Pb\ae$, or equivalently since $\aa$ is invertible,
\begin{align}\label{eq-H37}
\gamma(\Xb)=\gradt^P \psi(\Xb),\quad \Pb\ae
\end{align}
With \eqref{eq-H36}, this gives us \eqref{eq-z08}.
\\
Because of
\begin{align*}
dM^P_t=dM^R_t-\aa \gamma(\Xb_t)\,dt,
\qquad \Pb\ae,
\end{align*}
and \eqref{eq-H37}, 
equation \eqref{eq-H33} becomes
 \begin{align*}
d\psi(\Xb_t)
	=[\mathcal{L}^P \psi -|\gradt^P\psi|^2_\aa ](\Xb_t)\, dt
		+ \gradt^P \psi(\Xb_t)\cdot dM^R_t,
		\qquad \Pb\ae
\end{align*}
As   $\Iii | \gradt^P\psi|^2_\aa(\Xb_t)\,dt=\Iii | \gamma|^2_\aa(\Xb_t)\,dt< \infty,$ $P\ae$, this proves that  $ \psi$ is in  $\dom \LLL RP$ and that \eqref{eq-z07} holds.
Remark that $M^R$ is a $P$-local $R$-martingale.
\\
Finally, \eqref{eq-z08} and \eqref{eq-z07} directly imply \eqref{eq-H38}, while \eqref{eq-z09} follows from \eqref{eq-H37} and \eqref{eq-z07}.
\end{proof}

\subsection*{Feynman-Kac formula}

Let us go back to the function 
\begin{align*}
g_t(x):=E_R\left.\left(\exp \Big(\II tT V(\Xb_s)\,ds\Big)  g_T(X_T)\,\right|\,  X_t=x\right),
\quad (t,x)\in\iZ.
\end{align*}
introduced at \eqref{eq-z14}. 

\begin{theorem}\label{res-31}
Under the hypotheses of  Theorem \ref{res-13},
 $g$ belongs to $\dom \LLL RP$ and 
\begin{align*}
[\LLL RP  +V]g (\Xb)=0,\quad \Pb\ae
\end{align*}
\end{theorem}

\begin{proof}
We know by Theorem \ref{res-13} that  $ \psi(\Xb)$ is a $P$-local $R$-semimartingale. Applying  Itô formula to $g= e ^{ \psi}$, we immediately see that  $g(\Xb)$ is also  a $P$-local $R$-semimartingale. More precisely, we have 
\begin{align*}
d g(\Xb_t)
	&= g(\Xb_t) \Big[ d \psi(\Xb_t) + \ud d \left\langle \psi(\Xb), \psi(\Xb) 			\right\rangle _t\Big] \\
	&= g(\Xb_t) \Big[ \LLL RP \psi(\Xb_t)\,dt +dM ^{ \psi}_t+ \ud | \Gt RP \psi| 			^2_\aa (\Xb_t)\,dt \Big] ,\quad \Pb\ae,
\end{align*}
where we used  $d \psi(\Xb_t)= \LLL RP \psi(\Xb_t)\,dt +dM ^{ \psi}_t$ with $M ^{  \psi}$ a $P$-local $R$-martingale, which is an alternate statement for $ \psi\in\dom \LLL RP,$ and we also wrote $ d \left\langle \psi(\Xb), \psi(\Xb) 	\right\rangle _t=| \Gt RP \psi| ^2_\aa (\Xb_t)\,dt$ which is a direct consequence of \eqref{eq-H33} and \eqref{eq-H37}. Finally, with \eqref{eq-H38}:
$\big(\LLL RP \psi+| \Gt RP \psi|^2_\aa /2\big)(\Xb)=-V(\Xb),$ we arrive at
\begin{align*}
dg(\Xb_t)=-V \,g(\Xb_t) \,dt + dM^g_t,
\quad \Pb\ae,
\end{align*}
where $dM^g_t=g(\Xb_t)\, dM ^{ \psi}_t$ is the infinitesimal increment of a $P$-local $R$-martingale. This completes the proof of the theorem.
\end{proof}

\begin{remark}\label{rem-05}
This proof  relies on Theorem \ref{res-13} about the HJB equation. We didn't find a direct proof keeping the same minimal hypotheses. The main advantage of starting from HJB is the identification based on Girsanov theory which  led us to \eqref{eq-H36} and \eqref{eq-H37}.
\end{remark}

\section{Growth conditions}
\label{sec-GC}

We present at Theorem \ref{res-30} and Theorem \ref{res-32} some growth conditions on the coefficients $\aa,\BBf, V, g_T$ and $f_0$ which are sufficient for  $H(P|R)< \infty.$
\\
Consider the  Markov diffusion generator 
\begin{align*}
 A ^{ U} := \partial_t+\vU\scal \nabla + \Delta_\aa/2,
\end{align*}
where the velocity field is of the gradient form
\begin{align*}
\vU(t,x):=-\aa(t,x)\nabla U(t,x),\quad (t,x)\in\iZ, 
\end{align*} with $U$  some differentiable numerical function. The reference measure in next Theorem \ref{res-30} solves 
\begin{align*}
R=\MP(\aa,\vU+\vs)
\end{align*}
with $\vs$  a $\range(\aa)$-valued bounded vector field without any  regularity. \\
In this section, it is not necessary  to assume that $\aa$ is invertible.

\begin{hypotheses}\label{ass-01}\ 
\begin{enumerate}[(i)]
\item
$U\in C ^{ 1,2}(\iZ)$, 
\item
$\aa= \sigma \sigma ^{ \mathsf{t}}$ for some $ \sigma$ which is locally Lipschitz:\\ $ \displaystyle{\sup _{ x,y \in B, x\neq y, 0\le t\le T} \frac{| \sigma(t,y)- \sigma(t,x)|}{|y-x|}  < \infty}$, for any bounded subset $B\subset \ZZ$,
\item
for  some $K\ge 0,$\quad  $\sup _{ 0\le t\le T}\{x\scal \vU(t,x)+\tr \aa(t,x)\}\le K(1+|x|^2)$ for all $x\in\ZZ$.
\item
The measurable  vector field  $ \vs :\iZ\to\ZZ$  is bounded in the sense that
\begin{align}\label{eq-65}
\sup _{ \iZ}| \vs |_{ \mathsf{g}}< \infty.
\end{align}
\end{enumerate}
\end{hypotheses}
We denote $ \mathsf{g}:=\aa ^{ -1}$ the generalized inverse of $\aa$ and $|\vv|^2 _{ \mathsf{g}}:=\vv\scal \mathsf{g}\vv$ with $|\vv| _{ \mathsf{g}}= \infty$ when $\vv$ is outside the range of $\aa.$

 It is a standard result that under the hypotheses (i), (ii) and (iii), 
any martingale problem  associated with the characteristics $\vU$ and $\aa$ admits a unique solution. The uniqueness is a consequence of the Lipschitz hypothesis \ref{ass-01}-(ii). It implies  in particular \eqref{eq-H17a}. Let us denote by
\begin{align*}
\RU  =\MP(\aa,\vU;\mU),
\end{align*}  
the solution of this martingale problem with  some initial marginal $\mU.$
It is also true that adding the hypothesis (iv),
for any initial marginal $\ms $ the following martingale problem   admits a unique  (in the sense of \eqref{eq-H17a}) solution
\begin{align}\label{eq-79}
R=\MP(\aa,\vU+\vs;\ms).
\end{align}
 Let us give some details about this last assertion, when $\ms\ll\mU$ to keep  minimal notation.
We denote $$ \beta:= \mathsf{g}\vs,$$ recall that $\vs$ is assumed to  live in $\range(\aa)$. 
 By Novikov's criterion, Girsanov's formula  is valid under the finite energy estimate
\begin{align}\label{eq-70}
\sup\Iii | \beta_t|^2_\aa\,dt\le T\sup |\vs |^2 _{ \mathsf{g}}< \infty,
\end{align}
implied by the  hypothesis \eqref{eq-65}. This formula is:
\begin{align}\label{eq-72}
\frac{dR}{d\RU}
	= \frac{d\ms}{d\mU}(X_0) Z ^{ ( \beta)}_T,
\end{align}
where we set 
\begin{align}\label{eq-71}
Z_t ^{( \xi)}:=\exp \Big( \int _{ [0,t]} \xi_s\scal dM ^{ \RU}_s-\int _{ [0,t]}| \xi_s|^2_\aa/2\, ds\Big),
\quad 0\le t\le T,
\end{align}
with $dM ^{ \RU}_t=dX_t-\vU(X_t)\,dt$ the increment of a local $\RU$-martingale.
\\
The notation $Z ^{ (\xi)}$ will be used again later in the proof of Theorem \ref{res-30}. For any adapted process  $\xi$, $Z ^{ (\xi)}$ an $\RU$-supermartingale. When $\xi=\beta$ and  under \eqref{eq-70}, it is a genuine martingale. 

Now, we establish some preliminaries for the proof of Theorem \ref{res-30} at Lemma \ref{res-27} and Lemma \ref{res-29}.

\begin{lemma}\label{res-27}
For any function $H$   in $C ^{ 1,2}(\iZ)$, we set 
\begin{align*}
 \mathcal{H}:= e ^{ H}A^U e ^{ -H} 
 	 =- \partial_t H + \aa\nabla U\scal \nabla H- \Delta_\aa H/2+|\nabla H|_\aa^2/2 .
\end{align*}
The process 
\begin{align}\label{eq-62}
Z_t=\exp \left(H(\Xb_0)-H(\Xb_t)-\int _{ [0,t]} \mathcal{H}(\Xb_s)\,ds\right),
\quad 0\le t\le T,
\end{align}
is a $\RU$-supermartingale. In particular, if $\IZ e ^{-2H_0}\,d\mU< \infty,$   then
\begin{align*}
E _{ \RU } \exp \left(-H(\Xb_0)-H(\Xb_T)-\Iii \mathcal{H}(\Xb_t)\,dt\right) < \infty.
\end{align*}
\end{lemma}

\begin{proof}
Because it is a nonnegative local $\RU$-martingale, the process
\begin{align*}
Z_t:=\exp \left( - \Iii \nabla H(\Xb_t)\scal (dX_t-\vU(\Xb_t)\,dt)
		- \Iii|\nabla H|^2_\aa(\Xb_t)/2\,dt\right),
\qquad 0\le t\le T,
\end{align*}
is a $\RU$-supermartingale. With Itô formula 
\begin{align*}
-\Iii \nabla H(\Xb_t)\scal dX_t
	=H(\Xb_0)-H(\Xb_T)
		+\Iii \{ \partial_t+\Delta_\aa/2 \}H(\Xb_t)\,dt,
	\quad R^0\ae,
\end{align*}
we see that $Z$ is expressed by \eqref{eq-62}.
This implies
\begin{multline*}
E _{ \RU } \exp \left(-H(\Xb_0)-H(\Xb_T)-\Iii \mathcal{H}(\Xb_t)\,dt\right)
	= E_\RU [ e ^{ -2H(\Xb_0)}Z_T]\\
	\le E_\RU [ e ^{ -2H(\Xb_0)}Z_0]
	= E_\RU  e ^{ -2H(\Xb_0)}=\IZ e ^{ -2H_0}\, d\mU< \infty,
\end{multline*}
where the inequality is due to the supermartingale property of $Z$.
\end{proof}

\begin{lemma}\label{res-29}
For all measures $\rr,\qq$ and any probability measure $\pp$ on the same measurable space, such that $\pp\ll\qq\ll\rr,$ $E_\pp \log_+(d\pp/d\qq)< \infty$ and $E_\qq \max (1,d\qq/d\rr)< \infty,$ we have: 
$$H(\pp|\rr)\le 2H(\pp|\qq)+E_\qq \left(d\qq/d\rr\right)\in [- \infty, \infty),$$
where we set $H(\pp|\rr)=- \infty$ when $E_\pp\log_-(d\pp/d\rr)= \infty$ and $E_\pp\log_+(d\pp/d\rr)< \infty.$
\end{lemma}

\begin{proof}
We start considering positive parts of integrands to manipulate well-defined integrals:
\begin{align*}
E_\pp\log_+(d\pp/d\rr)
	=E_\pp \log_+(d\pp/d\qq\times d\qq/d\rr)
	\le E_\pp \log_+(d\pp/d\qq)+E_\pp \log_+(d\qq/d\rr).
\end{align*}
With the convex inequality: 
\[ab\le a\log a+ e ^{ b-1}\le a\log a+ e^b\le a\log_+ a+ e^b,\quad  \textrm{for all }a\ge 0,b\in\RR,
\]
 we see that
\begin{align*}
E_\pp \log_+(d\qq/d\rr)
	=E_\qq\big( d\pp/d\qq\times \log_+(d\qq/d\rr)\big)
	\le E_\pp \log_+(d\pp/d\qq)+E_\qq \max (1,d\qq/d\rr)< \infty
\end{align*}
which is finite by hypothesis. It follows that the integrals $H(\pp|\rr):=E_\pp\log(d\pp/d\rr),$  $H(\pp|\qq):=E_\pp\log(d\pp/d\qq)$ and $E_\pp\log(d\qq/d\rr)$ are well-defined in $[- \infty, \infty)$ and we are allowed to write
\begin{align*}
H(\pp|\rr)=E_\pp\log(d\pp/d\rr)
	=E_\pp \log(d\pp/d\qq)+E_\pp\log( d\qq/d\rr)
	=H(\pp|\qq)+E_\pp\log( d\qq/d\rr).
\end{align*}
Using the  convex inequality again, we obtain
\begin{align*}
E_\pp \log(d\qq/d\rr)
	=E_\qq\big( d\pp/d\qq\times \log(d\qq/d\rr)\big)
	\le H(\pp|\qq)+E_\qq (d\qq/d\rr)
\end{align*}
and conclude plugging this estimate into the above identity.
\end{proof}

Besides $U$, let us introduce another function $\Uo\in  C ^{ 1,2}(\iZ)$ and denote
\begin{align*}
\mathcal{U}&:= - e ^{ U}A^U e ^{ -U}=|\nabla U|_\aa^2/2- \partial_t U - \Delta_\aa U/2,\\ 
\UUo &:= - e ^{ \Uo}A ^{ \Uo} e ^{ -\Uo}=|\nabla \Uo|_\aa^2/2- \partial_t \Uo - \Delta_\aa \Uo/2.
\end{align*}

\begin{theorem}\label{res-30}
Let $\aa$, $U$ and $\vs$ entering the definition of $R$ at \eqref{eq-79} satisfy the Hypotheses \ref{ass-01}. 
Take   a function  $\Uo$ in $C ^{ 1,2}(\iZ)$ and a nonnegative measurable function $\Us :\ZZ\to [0, \infty)$ verifying
\begin{equation}\label{eq-61}
\begin{split}
(a)\quad \IZ  e ^{2U_0 -2\Uo_0 }&\,d\mU < \infty,\qquad
(b)\  \IZ \Us  e ^{ 2U_0-2\Uo_0 }\,d\mU< \infty,\\
&(c)\quad \IZ e ^{ -\Us}\,d\mU< \infty,
\end{split}
\end{equation}
recall the notation $\mU:=R^U_0.$
We assume  that  $R_0:=\ms$ is absolutely continuous with respect to $\mU$ and satisfies
\begin{align}\label{eq-67}
\sup \frac{|\log (d\ms /d\mU)|}{1+\Us}< \infty.
\end{align}
Let $f_0, g_T$ and $V$ verify \eqref{eq-87}-(c) and 
suppose that there exist   $c,\kappa\ge 0$ such that
\begin{align}\label{eq-59}
\begin{split}
 f_0\log_+f_0
	&\le \kappa\exp(U_0-\Uo_0),\\
 g_T\log_+g_T
	&\le \kappa\exp(U_T-\Uo_T),\\
V+ \log(1+ V_+)
	&\le  \mathcal{U} - \UUo +c.
\end{split}
\end{align}
Then $P$ defined by \eqref{eq-H03} can be normalized as a probability measure  and the relative entropy $H(P|R)$ is   finite.
\end{theorem}

\begin{proof}
We shall prove in a moment that $E_P\log_+(dP/d\RU )< \infty$. Since its proof does not require a priori that $P$ is a finite measure, this estimate implies that $dP/d\RU$ is finite. It follows with \eqref{eq-66} that $dP/dR$ is also finite, implying that \eqref{eq-87}-(a,b) is satisfied.

We divide  the rest of the proof into five steps. First, we consider the reference measure $R=\RU$ and show that under the assumption \eqref{eq-59} where $\Uo$ and $\Us$ satisfy \eqref{eq-61}, we have  $H(P|\RU)< \infty,$ requiring two steps: (i): $E_P\log_+(dP/d\RU )< \infty,$ 
and 
(ii): $E_P \log_-(dP/d\RU )< \infty.$
Then, we introduce $\vs$ and consider a specific $\ms$, given at \eqref{eq-68}. We first show that (iii): $H(P|R) < \infty$ , then (iv):   $H(P|R) >- \infty$ with this $\ms$. Finally, (v):  we extend the result to the case where $\ms$ satisfies \eqref{eq-67}.

\Boulette{$E_P\log_+(dP/d\RU )< \infty$}  By  definition \eqref{eq-H03} of $P$, with notation \eqref{eq-47}
\begin{multline} \label{eq-84}
\log_+  \Big(\frac{dP}{d\RU }\Big) 
	= \Big(\varphi_0(X_0)+ \psi_T(X_T) +\Iii V(\Xb_t)\, dt\Big)_+\\
	\le A:= [\varphi _{ 0}]_+(X_0)+ [\psi _{ T}]_+(X_T) +\big[\Iii V(\Xb_t)\, dt\big]_+.
\end{multline}
As
\begin{align*}
E_P \log_+  \Big(\frac{dP}{d\RU }\Big) \le E_P A
	=E_\RU \exp \left(\log A + \varphi_0(X_0)+ \psi_T(X_T)+\Iii V(\Xb_t)\,dt\right) ,
\end{align*}
by Lemma \ref{res-27} for the estimate ($+$) to hold it suffices that
\begin{multline*}
\log A + \varphi_0(X_0)+ \psi_T(X_T)+\Iii V(\Xb_t)\,dt
	\le -H_0(X_0)-H_T(X_T)-\Iii \mathcal{H}(\Xb_t)\,dt+c, 
\end{multline*}
 for some real $c\ge 0$ and some function $H$ in $C ^{ 1,2}(\iZ)$ verifying
\begin{align}\label{eq-60}
\IZ e ^{-2H_0}\,d\mU< \infty.
\end{align}
Since
\begin{align*}
\log A\le \log (1+ [\varphi _{ 0}]_+(X_0))
		+ \log( 1+ [\psi _{ T}]_+(X_T))
		+\log\Big( 1+\big[\Iii V(\Xb_t)\, dt\big]_+\Big),
\end{align*}
this is implied by
\begin{align*}
  \varphi_0(X_0) + \log (1+ [\varphi _{ 0}]_+(X_0))
	&\le  -H_0(X_0)+c,\\
  \psi_T(X_T)+\log( 1+ [\psi _{ T}]_+(X_T))
	&\le -H_T(X_T)+c,\\
 \Iii V(\Xb_t)\,dt +\log\Big( 1+\big[\Iii V(\Xb_t)\, dt\big]_+\Big)
	&\le -\Iii \mathcal{H}(\Xb_t)\,dt+c.
\end{align*}
Because  $[\Iii V(\Xb_t)\, dt]_+\le \Iii V_+(\Xb_t)\, dt,$
 for the above inequalities to be fulfilled, it is enough that there exists    $c\ge 0$ such that 
\begin{equation}\label{eq-63}
\begin{split}
 \varphi_0+ \log(1+ [\varphi _{ 0}]_+)
	&\le -H_0+c,\\
 \psi_T+ \log(1+[\psi _{ T}]_+)
	&\le - H_T+c,\\
V+ \log(1+ V_+)
	&\le - \mathcal{H} +c.
\end{split}
\end{equation}
Writing  $$H=-U+\Uo ,$$ these inequalities are \eqref{eq-59}  because
\begin{multline*}
- \mathcal{H}
	=\partial_t (-U+\Uo)- \aa\nabla U\scal \nabla (-U+\Uo)+ \Delta_\aa (-U+\Uo)/2-|\nabla (-U+\Uo)|_\aa^2/2\\
	= - \partial_t U- \Delta_\aa U/2+|\nabla U|_\aa^2/2
		+ \partial_t \Uo+ \Delta_\aa \Uo/2-|\nabla \Uo |^2/2
	= \mathcal{U}- \UUo .
\end{multline*}
Finally  \eqref{eq-60} becomes
$	
\IZ e ^{ -2\Uo_0 +2U_0}\,d\mU  < \infty,
$	
which is \eqref{eq-61}-(a).
These last considerations prove that \eqref{eq-61}-(a) and \eqref{eq-59} imply $E_P \log_+(dP/d\RU )< \infty.$

\Boulette{$E_P \log_-(dP/d\RU )< \infty$} 
By Proposition \ref{res-24} it is enough to obtain
\begin{align}\label{eq-69}
E_P \Us(X_0)< \infty
\end{align}
because it is assumed at \eqref{eq-61}-(c) that $E _{ \RU }e ^{ -\Us(X_0)}< \infty.$ We have
\begin{align*}
E_P \Us(X_0)
	&=E _{ \RU} \left[ \Us (X_0)\exp\Big( \varphi_0(X_0)+ \psi_T(X_T)+\Iii V(\Xb_t)\,dt\Big)\right] \\
	&\overset{(a)}{\le} E _{ \RU} \left[ \Us (X_0)\exp\Big( -H_0(X_0) -H_T(X_T)-\Iii \mathcal{H}(\Xb_t)\,dt+3c\Big)\right]\\
	&\overset{(b)}{=} e ^{ 3c} E _{ \RU} \left[ \Us (X_0)\exp( -2H_0(X_0)) Z_T\right] \\
	&\overset{(c)}{\le} e^{ 3c} E _{ \RU} \left[ \Us (X_0)\exp( -2H_0(X_0)) \right] \\
	&\overset{(d)}{=}e ^{ 3c}\IZ \Us  e ^{ -2H_0}\, d\mU
	=e ^{ 3c}\IZ \Us  e ^{ 2U_0-2\Uo_0 }\,d\mU.
\end{align*}
We used \eqref{eq-63} at (a),   the definition \eqref{eq-62} of the process $Z$ at (b), the fact that $Z=Z ^{ (-\nabla H)}$ is an $\RU$-supermartingale (see the proof of Lemma \ref{res-27}) at (c), and   $R^U_0=\mU$ at (d).
\\
Finally, our   assumption \eqref{eq-61}-(b)  implies that   $E_P \Us (X_0)$ is  finite, completing the proof of $E_P \log_-(dP/d\RU )< \infty$ and $- \infty<H(P|\RU)< \infty.$

Now, we consider the "bounded perturbation" $R$ of $\RU$ under the assumption \eqref{eq-65}. 
We choose for a while 
\begin{align}\label{eq-68}
\ms = e ^\Us\,\mU.
\end{align}
This will be relaxed later.

\Boulette{$H(P|R)< \infty$ under \eqref{eq-68}}
We already know that   $H(P|\RU)< \infty.$ Hence, Lemma \ref{res-29} gives us
\begin{align*}
- \infty\le H(P|R)\le 2H(P|\RU)+E _{ \RU}(d\RU/dR)<  \infty,
\end{align*}
provided that
\begin{align}\label{eq-66}
E _{ \RU}(d\RU/dR)<  \infty.
\end{align}
 Let us show that our assumptions \eqref{eq-65} and \eqref{eq-68} imply \eqref{eq-66}. 
 \\
The identities $dX_t=\vv(\Xb_t)\,dt+dM^R_t,$ $R\ae$, $dX_t=\vU(X_t)\,dt+dM ^{ \RU}_t,$ $\RU\ae,$ with $\vv=\vU+\vs =\vU+\aa \beta,$ imply that $dM^R=dM ^{ \RU}-\aa \beta\, dt$.  With Girsanov's formula  \eqref{eq-72}, we have
\begin{align*}
E _{ \RU}&(d\RU/dR)\\
	&= E _{ \RU}\bigg[\frac{d\mU}{d\ms }(X_0)\exp \Big( -\Iii \beta_t\scal dM^R_t-\Iii| \beta|^2_\aa/2\, dt\Big)\bigg] \\
	&=E _{ \RU}\bigg[\frac{d\mU}{d\ms }(X_0)\exp \Big( -\Iii \beta_t\scal dM ^{ \RU}_t+\Iii| \beta|^2_\aa/2\, dt\Big)\bigg]  \\
	&=E _{ \RU}\bigg[\frac{d\mU}{d\ms }(X_0)Z ^{( - \beta)}_T \exp \Big( \Iii| \beta|^2_\aa\, dt\Big)\bigg] 
	\le  \exp (T\sup| \beta|^2_\aa) E _{ \RU}\bigg[\frac{d\mU}{d\ms }(X_0)Z ^{( - \beta)}_T \bigg] \\
	&\le  \exp (T\sup| \beta|^2_\aa) E _{ \RU}\Big(\frac{d\mU}{d\ms }(X_0)\Big) 
	= \exp (T\sup| \vs |^2 _{ \mathsf{g}}) \IZ \frac{d\mU}{d\ms }\,d\mU,
\end{align*}
where last inequality holds because $Z ^{ (- \beta)}$
is a $\RU$-supermartingale,  recall notation \eqref{eq-71}. 
\\
It is assumed that
$
\IZ e ^{ -\Us}\,d\mU< \infty 
$
and we know that 
$
\IZ \Us\,dP_0< \infty,
$
see \eqref{eq-69}.
Hence $\IZ {d\mU}/{d\ms }\ d\mU=\IZ e ^{ -\Us}\,d\mU< \infty,$ proving \eqref{eq-66} when \eqref{eq-68} holds.

\Boulette{$H(P|R)>- \infty$  under \eqref{eq-68}}
By Proposition \ref{res-24}, this will be done by showing that  
$$E_R e ^{ -2\Us(X_0)}< \infty,$$  
because we already know that $E_P\Us(X_0)< \infty$, see \eqref{eq-69}.
   Let us control
\begin{align*}
E_R e ^{ -2\Us(X_0)}
	=E _{ \RU} \left( \frac{dR}{d\RU} e ^{ -2\Us(X_0)}\right) 
	\le E _{ \RU} \left[ \left( \frac{dR}{d\RU}\right) ^2 e ^{ -3\Us(X_0)}\right] ^{ 1/2}
		\left( E _{ \RU} e ^{ -\Us(X_0)}\right) ^{ 1/2}
\end{align*}
by Cauchy-Schwarz inequality. As
\begin{multline*}
 \left( \frac{dR}{d\RU}\right) ^2
 	= \left( \frac{d\ms }{d\mU}(X_0)\ Z ^{ ( \beta)}\right) ^2
	= \left( \frac{d\ms }{d\mU}(X_0)\right) ^2
		Z ^{ (2 \beta)}_T\exp \left(\Iii | \beta_t|^2_\aa\, dt\right)\\
	\le \exp \left(T \sup | \beta|^2_\aa\right) \left( \frac{d\ms }{d\mU}(X_0)\right) ^2
		Z ^{ (2 \beta)}_T 
\end{multline*}
where $Z ^{ (2 \beta)}$ is a $\RU$-supermartingale, we obtain
\begin{multline*}
E _{ \RU} \left[ \left( \frac{dR}{d\RU}\right) ^2 e ^{ -3\Us(X_0)}\right]
	\le \exp \left(T \sup | \beta|^2_\aa\right)
	\IZ \left( \frac{d\ms }{d\mU}\right) ^2 e ^{ -3\Us}\, d\mU\\
	= \exp \left(T  \sup | \beta|^2_\aa\right)
	\IZ e ^{ -\Us}\, d\mU ,
\end{multline*}
showing that
\begin{align*}
E_R e ^{ -2\Us(X_0)}\le  \exp \left(T\sup | \vs |^2 _{ \mathsf{g}}/2 \right)
	 \IZ e ^{ -\Us}\, d\mU < \infty
\end{align*}
and completing the proof of $$- \infty<H(P|R)< \infty $$ when $\ms $ is given by \eqref{eq-68}.

\par\medskip\noindent $\bullet$\ Relaxation of \eqref{eq-68}.\ 
Let us extend this property to any  $\ms $ satisfying \eqref{eq-67}. As for any positive function $r$ on $\ZZ$ we have
\begin{align*}
H(P| r(X_0)\,R)
	=H(P|R)-E_P\log r(X_0),
\end{align*}
one can extend our  previous result from $R$ to $r(X_0)\,R$ provided that $E_P|\log r (X_0)|< \infty.$  Simply requiring that $|\log r|\le \kappa (1+\Us )$ for  some $\kappa \ge 0$ so that $E_P |\log r (X_0)|\le \kappa (1+ E_P \Us(X_0)) < \infty,$ permits us to extend our result from $\ms = e^\Us  \,\mU$ to $\ms = re^\Us \, \mU$ with $\sup \{| \log r|/(1+\Us )\}< \infty,$ which is the assumption \eqref{eq-67}.
\end{proof}

\begin{remarks} 
 \begin{enumerate}[(a)]
\item
In the case where $\RU$, hence $\mU$, is a bounded measure, one can choose $\Us =0$ because it verifies
$	
\IZ e ^{ -\Us}\,d\mU< \infty.
$
Moreover, the assumption  \eqref{eq-61} simplifies as $\IZ  e ^{ 2U_0-2\Uo_0 }\,d\mU < \infty$ and \eqref{eq-67} becomes $\sup |\log (d\ms /d\mU)|< \infty.$

\item
One can always choose $\ms =\mU$ in \eqref{eq-67}. 

\item
There is some freedom in the choice of $\mU$ because the initial measure $R_0$ does not appear explicitly in the Feynman-Kac equation \eqref{eq-FK} and it is erased by the conditional expectation in the Feynman-Kac formula \eqref{eq-74}. A possible choice is
\begin{align*}
\mU= e ^{ -2U_0}\,\Leb.
\end{align*}
In this case, the hypotheses \eqref{eq-61} become
\begin{align*}
\IZ  e ^{-2\Uo_0 }\,d\Leb < \infty,\qquad
 \IZ \Us  e ^{ -2\Uo_0 }\,d\Leb< \infty,\qquad
 \IZ e ^{-2U_0 -\Us}\,d\Leb< \infty,
\end{align*}
and when $\IZ e ^{ -2U_0}\,d\Leb < \infty$, taking $\Us=0,$ a possible choice for $\Uo$ is
\begin{align*}
\Uo(t,x)=\Uo_0 (x)= \gamma \log \sqrt{1+|x|^2},
\quad \textrm{ with }  \gamma>n/2
\end{align*} 
to ensure $\IZ  e ^{-2\Uo_0 }\,d\Leb < \infty.$
It gives
\begin{align*}
\UUo (t,x)
	= \frac{ \gamma (2+ \gamma)|x|^2 _{ \aa(t,x)}- \gamma \tr\aa(t,x)(1+|x|^2)}{2(1+|x|^2)^2} .
\end{align*}
\item
Because for any $ \varepsilon>0$ and $0<q\le 1$, there is some $c( \varepsilon)\ge 0$ such that  
\[
a+\log(1+a_+)
	\le a+ \varepsilon a_+^q+ c( \varepsilon),
	\quad \forall a\in\RR,
\]
for the upper bounds \eqref{eq-59} to hold, it suffices that
 there exist $ \varepsilon>0$, $0<q\le 1$ and $c\ge 0$ such that
\begin{align*}
\begin{split}
 \varphi_0+ \varepsilon [\varphi _{ 0}]_+^q
	&\le U_0-\Uo_0 +c,\\
 \psi_T+ \varepsilon[\psi _{ T}]_+^q
	&\le U_T-\Uo_T +c,\\
V+ \varepsilon V_+^q
	&\le  \mathcal{U} - \UUo +c.
\end{split}
\end{align*}
 \end{enumerate}\end{remarks}

Roughly speaking, the hypothesis \eqref{eq-59}  imposes that $V$ should not grow faster than the opposite of  the confinement potential $ \mathcal{U}$ associated to $\RU.$ It implies in particular that $V$ is locally upper bounded because $ \mathcal{U}$ is continuous. Next result presents a set of hypotheses where this is relaxed.

\subsection*{Kato class}

It is  known since the article \cite{Kh59} by Khas'minskii  that if for some $ \tau>0,$
\begin{align*}
 \sup _{ x\in \ZZ} \sup _{0\le s<s+ \tau\le T}
	E _{ Q} \left.\left( \int _{ [s,s+ \tau]}|W|( \Xb_t)\,dt \,\right|\,  X_s=x\right) =: \alpha< 1,
\end{align*}
then
\begin{align*}
\exp( \Lambda( \tau)):=\sup _{ x\in\ZZ}\sup _{0\le s<s+ \tau\le T}
	E _{ Q}\left. \left[\exp\left( \int _{ [s,s+ \tau]} e^{|W|( \Xb_t)}\,dt\right) \,\right|\,  X_s=x\right]
	\le \frac{1}{1- \alpha}.
\end{align*}
The Markov property of $Q$ is essential to prove this result.
On the other hand, the Markov property also implies that $ \tau\mapsto \Lambda( \tau)$ is subadditive. Consequently, the finiteness of $ \Lambda( \tau)$ for a small enough positive $ \tau$  implies that there exists $c>0$ such that $ \Lambda(t)\le c+ct$ for all $0\le t\le T,$ proving that
\begin{align}\label{eq-80}
\sup _{ x\in\ZZ}
	E _{ Q}\left. \left[\exp\left( \Iii e^{|W|( \Xb_t)}\,dt\right) \,\right|\,  X_0=x\right]
	\le e ^{c+cT }.
\end{align}

\begin{definition}[Kato class] \label{def-kato}
A measurable function $W:\iZ\to (- \infty, \infty)$ belongs to the Kato class $J(Q)$ of the Markov measure $Q\in\MO$  if
\begin{align*}
 \sup _{ x\in \ZZ} 
	E _{ Q}\left. \left( \Iii |W|( \Xb_t)\,dt \,\right|\,  X_0=x\right) < \infty
\end{align*}
 and
\begin{align*}
\lim _{ h\downarrow 0} \sup _{ x\in \ZZ} \sup _{0\le s<s+h\le T}
	E _{ Q} \left.\left( \int _{ [s,s+h]}|W|( \Xb_t)\,dt \,\right|\,  X_s=x\right) =0.
\end{align*}
We say that $W$ is in $J^*(Q)$ if the function $W^*$ defined by $W^*(t,x):=W(T-t,x),$ $(t,x)\in\iZ,$ is in $J(Q^*)$ where $Q^*$ is the time reversal of $Q.$
\end{definition}
It follows from the above considerations that any $W$  in $J(Q)$ verifies \eqref{eq-80} and we see  immediately  that $J(Q)$ contains all the bounded functions.
\\
Let us introduce the Feynman-Kac operators
\begin{align}\label{eq-82}
Su(x)
	:=E _{ Q}\left. \left[\exp\left( \Iii e^{|W|( \Xb_t)}\,dt\right) u(X_T) \,\right|\,  X_0=x\right],\quad x\in\ZZ,
\end{align}
and
\begin{align*}
S^*u(y)
	:=E _{ Q} \left.\left[u(X_0)\, \exp\left( \Iii e^{|W|( \Xb_t)}\,dt\right)  \,\right|\,  X_T=y\right],\quad y\in\ZZ,
\end{align*}
defined for any measurable function $u:\ZZ\to\RR$ such that the above conditional expectations are meaningful. Note that $Su$ is defined $Q_0\ae$ and $S^*u$ is defined $Q_T\ae$ Of course, if $Q$ is reversible and $W$ does not depend explicitly on the time variable, i.e.\, $W(t,x)=W(x)$, then $S^*=S.$

\begin{lemma}\label{res-33}
Let $Q\in\MO$ be a Markov measure and  $W:\iZ\to (- \infty, \infty)$ a measurable function in the Kato class $J(Q)\cap J^*(Q)$. Then, for any $1\le p\le \infty,$ the linear operators  $S: L^p(Q_T)\to L^p(Q_0)$ and $S^*: L^p(Q_0)\to L^p(Q_T)$ are bounded.\end{lemma}

\begin{proof}
By \eqref{eq-80},  $S: L^\infty(Q_T)\to L^\infty(Q_0)$ is a bounded operator with  $\|S \| _{\infty, \infty }\le e ^{ c+cT}.$
Similarly, as  $W\in J^*(Q),$  we see  by time-reversal that    $S^*: L^\infty(Q_0)\to L^\infty(Q_T)$ is a bounded operator with  $\|S^* \| _{\infty, \infty }\le e ^{ c^*+c^*T}.$
For any $g\in L^1(Q_T),$
\begin{align*}
\IZ |Sg|\,dQ_0
	\le \IZ S|g|\,dQ_0=\IZ |g|S^*\1\,dQ_T\le \|S^* \| _{\infty, \infty }\,\|g\| _{ L^1(Q_T)},
\end{align*}
proving that $S: L^1(Q_T)\to L^1(Q_0)$ is a bounded operator with  $\|S \| _{1, 1 }\le \|S^* \| _{\infty, \infty }.$  The term $Sg$ in the first integral is justified a posteriori by the finiteness of the second integral. The equality follows from the properties of the conditional expectations of nonnegative functions.
\\
By  the Riesz-Thorin interpolation theorem it is also true that  for any $1\le p\le \infty,$ $S: L^p(Q_T)\to L^p(Q_0)$ is a bounded operator, and a similar proof works with $S^*.$
\end{proof}

It is uneasy to obtain practical sufficient conditions for a function to belong to $J(Q)$ except when some upper bound  for the transition  kernel is known. Of course, the Gaussian case is well understood, corresponding to $Q$ being the law of a Brownian or an Ornstein-Uhlenbeck process. 
 For a clear exposition of the properties of the Kato class  and the semigroup of  Feynman-Kac operators  of the Wiener measure, see \cite[Ch.\,3]{CZ95}. In particular, when $Q$ is the Wiener measure $  \mathcal{W}^{ \epsilon}$  of the reversible Brownian motion with diffusion coefficient $ \epsilon>0,$ that is
\begin{align*}
 \mathcal{W} ^{ \epsilon}=\MP( \epsilon\Id,0;\Leb),
\end{align*}
with Markov generator $ \epsilon \Delta/2$ and the Lebesgue measure as its initial marginal, we know that $ \mathcal{W} ^{ \epsilon} = ( \mathcal{W}^ \epsilon)^*$ and $W:\ZZ\to\RR$ is in $J( \mathcal{W}^ \epsilon)$   if and only if
\begin{align}\label{eq-86}
\lim _{ \alpha\downarrow 0} \sup _{ x\in\ZZ}
	 \int _{ |y-x|\le \alpha} | \mathsf{g}(y-x)W(y)|\,dy=0
\end{align}
where $ \mathsf{g}$ is related to  the Green potential  and defined by
\begin{align*}
 \mathsf{g}(z):= \left\{ \begin{array}{ll}
 |z| ^{ 2-n}, & n\neq 2;\\
 \log(1/|z|), &n=2.\\
 \end{array}\right.
\end{align*}

\subsection*{A variant of Theorem \ref{res-30}}
We are now ready to prove another criterion for $H(P|\RU)< \infty$ with a possibly locally unbounded potential $V.$

Taking $U=0$ in the definition of $\RU$, we obtain a law 
$
R ^{ 0,\aa}=\MP( \Delta_\aa/2)=\MP(\aa,0)
$
of a Markov Brownian martingale with diffusion matrix $\aa,$ whose generator is $A^0= \partial_t+ \Delta_\aa/2.$ In this section we prefer looking at   the law 
\begin{align*}
R^\aa=\MP(\nabla\scal (\aa\nabla \sbt)/2;\Leb)
\end{align*}
of the Markov process with initial measure $R^\aa_0=\Leb$ and generator
$
\nabla\scal (\aa\nabla \sbt)/2
$
in divergence form.  When $\aa$ does not depend on time explicitly, it is reversible with Lebesgue measure as reversing measure; in particular $J^*(R^\aa)=J(R^\aa)$. In the more general case where $\aa$ depends on $t$, $R^\aa$ is not reversible anymore but its marginal flow remains constantly equal to $\Leb.$

\begin{theorem}\label{res-32}
Let
\begin{align*}
R=\MP\Big((\vU+\vs)\scal\nabla +\nabla\scal (\aa\nabla \sbt)/2\Big)
\end{align*}
with  $\aa$, $U$ and $\vs$ 
satisfying the Hypotheses \ref{ass-01}. 
Let $h_0, h_T$ and $\Us$ be nonnegative  measurable functions on $\ZZ$ such that
\begin{align}\label{eq-81}
(1+\Us)h_0\in L^p(\Leb),\qquad
h_T\in L ^{ p'}(\Leb),
\end{align}
with $1\le p\le \infty$, $1/p+1/p'=1,$ and
\begin{equation}\label{eq-81b}
\IZ e ^{ -\Us-2U_0}\,d\Leb< \infty.
\end{equation}
We assume  that  $R_0$ is absolutely continuous and satisfies
\begin{align}\label{eq-67b}
\sup \frac{|\log (dR_0 /d\Leb)+2U_0|}{1+\Us}< \infty.
\end{align}
Let $f_0, g_T$ and $V$ verify \eqref{eq-87}-(c) and 
suppose that there exists a measurable function   $$W\in J(R^\aa)\cap J^*(R^\aa)$$ 
 such that
\begin{align}\label{eq-59b}
 f_0
	\le  e ^{ U_0}\, h_0,\qquad
 g_T
	\le  e ^{ U_T}\,h_T,\qquad
V+ \log(1+ V_+)
	\le  \mathcal{U} +W.
\end{align}
Then $P$ defined by \eqref{eq-H03} can be normalized as a probability measure  and the relative entropy $H(P|R)$ is   finite.
\end{theorem}

\begin{proof}
As in the proof of Theorem \ref{res-30}, the announced estimate $E_P\log_+(dP/dR )< \infty$  implies that $dP/dR$ is finite, hence \eqref{eq-87}-(a,b) is satisfied.

\Boulette{$E_P\log_+ \frac{dP}{d\RU}< \infty$}
In this proof we change a little bit the path measure 
\begin{align}\label{eq-79c}
\RU=\MP(\vU\scal\nabla +\nabla\scal (\aa\nabla \sbt)/2;\mU),\qquad \mU=e ^{ -2U_0}\,\Leb.
\end{align}
  We replace $ \Delta_\aa$ by $\nabla\scal(\aa\nabla \sbt)$ but keep the same notation $\RU.$ 
 \\
Clearly \eqref{eq-59b} implies 
 \begin{align}\label{eq-59c}
 \begin{split}
 &f_0\log_+f_0+f_0
	\le  3e ^{ U_0}\, h_0,\quad
 g_T\log_+g_T+g_T
	\le  3e ^{ U_T}\,h_T,\\
&V+ \log(1+ V_+)
	\le  \mathcal{U} +W.
 \end{split}
\end{align}
 As in the proof of Theorem \ref{res-30} we obtain
\begin{align*}
 \frac{d\RU}{dR^\aa}
 	= \frac{d\mU}{d\Leb}(X_0)\exp \left( U(\Xb_0)-U(\Xb_T)-\Iii \mathcal{U}(\Xb_t)\,dt\right) 
\end{align*}
and with \eqref{eq-59c} we see that
\begin{align*}
E_P\log_+ \frac{dP}{d\RU}
	&\le E _{ \RU} \exp \left( \log \frac{dP}{d\RU}+\log \left( 1+\log_+ \frac{dP}{d\RU}\right) \right) \\
	&\le E _{ \RU} \exp \left( \log A+ \varphi_0(X_0)+ \psi_T(X_T)+\Iii V(\Xb_t)\,dt\right) \\
	&\le 9\, E _{ \RU} \left\{ \left( \frac{d\RU}{dR^\aa}\right)  ^{ -1} h_0(X_0)h_T(X_T) \exp \left( \Iii W(\Xb_t)\,dt\right) \right\} \\
	&= 9\,E _{ R^\aa} \left\{ h_0(X_0) \exp \left( \Iii W(\Xb_t)\,dt\right) \,h_T(X_T)\right\} \\
	&=9\,E _{ R^\aa} \left\{ h_0(X_0) [Sh_T](X_0)\right\} ,
\end{align*}
where the notation $A$ is used at \eqref{eq-84} and  $S$ is the Feynman-Kac operator associated with $R^\aa$ and $W\in J(R^\aa)\cap J^*(R^\aa).$ We conclude   with Lemma \ref{res-33} under the assumption \eqref{eq-81} that $E_P\log_+ \frac{dP}{d\RU}< \infty.$

\Boulette{$E_P\log_- \frac{dP}{d\RU}< \infty$}
By Proposition \ref{res-24}, under the assumption \eqref{eq-81} it suffices to prove that $E_P\Us(X_0)< \infty.$ Proceeding as above
\begin{align*}
E_P\Us(X_0)
	&= E _{ \RU} \left[\Us(X_0)\exp \left( \varphi_0(X_0)+\psi_T(X_T)
		+\Iii V(\Xb_t)\,dt\right) \right] \\
	&\le 9\, E _{ R^\aa} \left\{\Us h_0(X_0) [Sh_T](X_0)\right\}
\end{align*}
and we conclude   with Lemma \ref{res-33} under the assumption \eqref{eq-81} that $E_P\Us(X_0)< \infty.$

\Boulette{ the rest}
Already done at Theorem \ref{res-30}.
\end{proof}

\begin{remarks}\ \begin{enumerate}[(a)]
\item
Since any Kato class contains all the bounded functions, in \eqref{eq-59b} it is valid   to choose $$W=c.$$  
\item
Again, if $\IZ e ^{ -2U_0}\,d\Leb< \infty,$ taking $\Us=0$  in \eqref{eq-81} and \eqref{eq-81b} is all right.
\end{enumerate}\end{remarks}

\begin{corollary}[Classical cases]\ 
\begin{enumerate}[(a)]
\item \emph{Brownian motion}. 
When $R=\MP(\epsilon \Delta/2;\Leb)$ with $ \epsilon>0$, for  the relative entropy $H(P|R)$ to be finite, it suffices that for some $1\le p\le \infty,$
\begin{align*}
&(1+\log_+|\sbt|) f_0\in L^p(\Leb),\qquad
g_T\in L ^{ p'}(\Leb),\\
& \sup _{ t\in\ii}V_-(t,\sbt)\in L^1 _{ \mathrm{loc}}(\Leb),
\qquad \sup _{ t\in\ii}V_+(t,\sbt) \textrm{ verifies } \eqref{eq-86}.
\end{align*}

\item \emph{Ornstein-Uhlenbeck process}. 
When $R=\MP\Big(-k x\scal\nabla+\epsilon \Delta/2;  \mathcal{N}_{\epsilon/(2k)} \Big)$ with $ \epsilon,k>0$, for  the relative entropy $H(P|R)$ to be finite, it suffices that for some $1< p< \infty$ and some $W$ verifying \eqref{eq-86}
\begin{align*}
&f_0\in L^p( \mathcal{N}_{\epsilon/(pk)}), \qquad
g_T\in L ^{ p'}( \mathcal{N}_{\epsilon/(p'k)}),
\qquad \quad \sup _{ t\in\ii}V_-(t,\sbt)\in L^1 _{ \mathrm{loc}}(\Leb)\\
&\sup _{ t\in\ii}V_+(t,x)
	\le \frac{k^2}{2 \epsilon}\,|x|^2-2\log_+|x|+W(x),\quad \forall x\in\ZZ,
\end{align*}
where $ \mathcal{N}_a$ stands 
 for the normal distribution with zero mean and variance  $a\Id.$
  \end{enumerate}
\end{corollary}

\begin{proof}
\boulette{(a)}
Apply Theorem \ref{res-32} with $\vs=0,$ $U=0$,  $\Us=(n+1)\log_+|\sbt|,$ and remark that a Kato class is a vector space which is stable by the lattice operations.  The assumption about $V_-$ implies \eqref{eq-87}-(c).
 
\Boulette{(b)} Apply Theorem \ref{res-32} with $\vs=0,$ $U(x)=k|x|^2/(2 \epsilon)$, $ \mathcal{U}(x)=k^2|x|^2/(2 \epsilon)-nk/2$ and $\Us=0.$ Remember that a Kato class is a lattice vector space which contains the constants, and take advantage of: $\forall a,b\ge 0, a\le b-\log(1+b)\implies a+\log(1+a)\le b,$ applied to $a=\sup_t V_+.$ The assumption about $V_-$ implies \eqref{eq-87}-(c).
\end{proof}

\begin{remarks}\ \begin{enumerate}[(a)]
\item
In both cases  $R$ is chosen to be reversible.
\item
By time symmetry, (a) also holds if the hypothesis on $f_0,g_T$ is replaced by
\begin{align*}
 f_0\in L^p(\Leb),\qquad
(1+\log_+|\sbt|)g_T\in L ^{ p'}(\Leb).
\end{align*}
\end{enumerate}\end{remarks}

\appendix

\section{Carré du champ}

Lemma \ref{res-22} below is a simplified version of  \cite[Lemma 3.9]{CCGL20}, which was used during the proof of Lemma \ref{res-21}. For the confort of the reader, we give its detailed proof which is slightly simpler, but essentially the same as  \cite{CCGL20}'s one.

Let $Q\in\MO$ be a conditionable path measure.
Its forward  carré du champ  is defined by
\begin{align*}
\GaF^Q(u,v):= \LLF^Q(uv)-u\LLF^Q v-v\LLF^Q u,
\quad \ 0\le t\le T,
\end{align*}
 for any functions $u,v$ in $\dom\LLF^Q$ such that their product $uv$ also belongs to  $\dom\LLF^Q$.

The quadratic covariation $[u(X),v(X)]$ is a $Q$-semimartingale.
We denote by $\langle u(X),v(X)\rangle^Q$ its bounded variation part,  i.e.
\begin{align*}
d[u(X),v(X)]_t=d\langle u(X),v(X)\rangle^Q_t+ dM_t ^{Q, [u,v]},
\qquad \Qb\ae
\end{align*}
where, here and below, $M^{Q,\sbt}$ or $M ^{ \sbt}$ stands for any local $Q$-martingale. As next lemma indicates, we are interested in situations where the bounded variation process $\langle u(X),v(X)\rangle^Q$ is predictable (as a continuous process). Therefore, in the whole article $\langle u(X),v(X)\rangle^Q$ is the usual sharp bracket (sometimes called conditional quadratic variation) of stochastic process theory.

\begin{lemma}\label{res-13b}
For any $u,v\in\dom\LLF^Q$ such that $uv\in\dom\LLF^Q,$ the process $ \langle u(X),v(X)\rangle^Q $ is absolutely continuous $Q\ae$ and 
\begin{align*}
 d\langle u(\Xb),v(\Xb)\rangle^Q_t=\GaF^Q(u,v)(t, X _{ [0,t]})\,dt, \qquad Q\ae
\end{align*}
\end{lemma}

\begin{proof}
As a definition of the forward generator
$
d(uv)(\Xb_t)=\LLF^Q_t(uv)(t, X _{ [0,t]})\,dt+dM ^{ uv}_t.
$ 
Comparing this expression with \eqref{eq-56}, the Doob-Meyer decomposition theorem  gives the announced result.
\end{proof}

We say that a process $Y$ \emph{can be localized as} a bounded (resp.\ integrable)  process if there exists a sequence of stopping times $( \sigma_k)$ tending almost surely to infinity and such that for each $k$, the stopped process $Y^{ \sigma_k}$ is bounded almost surely (resp.\ integrable).

\begin{lemma}\label{res-22}
For any conditionable path measure $Q\in\MO$, almost every $t\in\ii,$ and  any locally bounded functions $u,v\in\dom\LLF^Q$ such that $uv\in\dom\LLF^Q$, and    $M ^{ Q,[u,v]}$ as defined at Lemma \ref{res-13b} can be localized as an integrable $Q$-martingale,    there exist an increasing sequence $(\tau_k)$ of $Q$-integration times of $u$ and $v,$ and a sequence $(h_n)$ of positive numbers such that $\Lim k \tau_k=\infty,$ $Q\ae$, $\Lim n h_n=0$ and for each $k$ we have 
\begin{equation*}
\begin{split}
 \GaF^Q (u,v)&(t, X _{ [0,t]})\\&= \Lim k \Lim n \frac 1{h_n}E_Q \left.\Big[\{u( \Xb ^{ \tau_k}_{t+h_n})-u( \Xb ^{ \tau_k}_t)\}\{v( \Xb ^{ \tau_k}_{t+h_n})-v( \Xb ^{ \tau_k}_t)\}\,\right|\,  X_{[0,t]}  \Big],\quad Q\ae
\end{split}
\end{equation*}
\end{lemma}

\begin{proof}
Since $u$ and $v$ are assumed to be  locally bounded,  $u(\Xb)$ and $v(\Xb)$  can be localized as bounded processes. Furthermore, the processes $\int_0^{\sbt}|\LLF^Qu(t,X _{ [0,t]})|\,dt$ and  $\int_0^{\sbt}|\LLF^Qv(t,X _{ [0,t]})|\,dt$ can also be localized as bounded processes. It follows that the local $Q$-martingales $M^u, M^v$,   (where $M^u_t:=u(\Xb_t)-\int_0^t\LLF^Qu(s,X _{ [0,s]})\,ds$) can also be localized as bounded processes.
Localizing as in the proof of Proposition \ref{res-H02c}, it is enough to show that
 \begin{align}\label{eq-55}
 \begin{split}
\Limh E_Q\int_0 ^{ T-h} \Big|E_Q \big[h ^{ -1}\{u( X_{t+h})-u( X_t)\}\{v( X_{t+h})-&v( X_t)\}\,\big|\, X _{[0,t]}\big]\\
	&-   \GaF^Q(u,v)(t,X _{ [0,t]})\,\Big|\,dt =0,
 \end{split}
\end{align} 
and we can assume that  all the above mentioned  processes are bounded.
\\
For each $0\le t\le T-h$ with $0<h\le T,$
\begin{align*}
[u( X_{ t+h})-&u( X_t)][v( X _{ t+h})-v( X_t)]\\
	= &\Big[ \int_t ^{ t+h} dM^u_s
		+ \int_t ^{ t+h} \LLF^Qu(\Xb_s)\,ds\Big]
		 \Big[ \int_t ^{ t+h} dM^v_s
		+ \int_t ^{ t+h} \LLF^Qv(\Xb_s)\,ds\Big]\\
	= & A_t^h+B_t^h+C_t^h+D_t^h,
	\qquad Q\ae,
\end{align*}
where
\begin{align*}
A_t^h&= \int_t ^{ t+h} dM^u_s\ \int_t ^{ t+h} dM^v_s,\hskip 2,5cm
B_t^h=  \int_t ^{ t+h} \LLF^Qu(\Xb_s)\,ds\  \int_t ^{ t+h} dM^v_s,\\
C_t^h&=  \int_t ^{ t+h} \LLF^Q v(\Xb_s)\,ds\  \int_t ^{ t+h} dM^u_s,\hskip 1,1cm 
D_t^h=  \int_t ^{ t+h}\LLF^Qu(\Xb_s)\,ds\  \int_t ^{ t+h} \LLF^Qv(\Xb_s)\,ds.
\end{align*}
Let us control $A^h_t.$ Denoting $U_{t,s}:=M^u_s-M^u_t$ and 
$V_{t,s}:=M^v_s-M^v_t,$
\begin{align*}
A_t^h&=
	\int_t ^{ t+h} d(U_{t,s}V_{t,s})\\
&= \int_t ^{ t+h} U_{t,s}  dM^v_s
		+ \int_t ^{ t+h} V_{t,s}  dM^u_s
		+ \int_t ^{ t+h} dM ^{Q, [u,v]}_s
		+ \int_t ^{ t+h} d \langle M^u,M^v\rangle ^Q _s,
\end{align*}
and  with Lemma \ref{res-13b}
\begin{align}\label{eq-02}
h ^{ -1}E_Q\left.[A^h_t\,\right|\,  X _{ [0,t]}]= h ^{ -1}\int_t ^{ t+h} E_Q[\GaF^Q(u,v)(\Xb_s)\mid X _{ [0,t]}]\,ds.
\end{align}
Remark that the boundedness  properties obtained above by localization, together with the extra assumption that $M ^{Q, [u,v]}$ is integrable, justify the cancelation of the expectations of the martingale terms.
\\
Let us control $B^h$:
\begin{align*}
h ^{ -1}E_Q \int_0 ^{ T-h}|B_t^h|\,dt
	&\le  E_Q \int_0 ^{ T-h} h ^{ -1}\Big|\int_t ^{ t+h} \LLF^Qu(\Xb_s)\,ds\Big|\ |M^v_{ t+h}-M^v_t|\, dt
	\\
	&=  E_Q \int_0 ^{ T-h} \Big| k^h\ast (\LLF^Qu)(\Xb_t)\Big|\ |M^v_{ t+h}-M^v_t|\, dt\\
	&= o _{ h\to 0^+}(1) ,
\end{align*}
where we took $k^h:= h ^{ -1} \1 _{ [-h,0]}$ as our convolution kernel. The last identity is a consequence of  Lemma \ref{resh-21}  under the assumption  $\LLF^Qu(\Xb)\in L^1(\Qb)$ (because $\int_0^{\sbt}|\LLF^Qu(t,X _{ [0,t]})|\,dt$ is bounded), the uniform boundedness and right-continuity of $M^v$ and the dominated convergence theorem. 
\\
Similarly, $\Limh h ^{ -1}E_Q \int_0 ^{ T-h}|C_t^h|\,dt=0.$
\\
The control of $D^h$  is analogous:
\begin{align*}
h ^{ -1}E_Q \int_0 ^{ T-h}|D_t^h|\,dt
	&\le  E_Q \int_0 ^{ T-h} \Big| k^h\ast (\LLF^Qu)(\Xb_t)\Big|\ \Big|\int_t ^{ t+h} \LLF^Qv(\Xb_s)\,ds\Big|\, dt\\
	&= o _{ h\to 0^+}(1) ,
\end{align*}
thanks to the uniform boundedness of $\Iii |\LLF^Qv(\Xb_s)|\,ds.$
\\
Putting everything together, we obtain
\begin{align*}
\Limh E_Q\int_0 ^{ T-h} \Big|E_Q \big[h ^{ -1}\{u( X_{t+h})-u( X_t)\}&\{v( X_{t+h})-v( X_t)\}\,\Big|\, X _{ [0,t]}\big]\\
	&-h ^{ -1}\int_t ^{ t+h} E_Q[\GaF^Q(u,v)(\Xb_s)\,\big|\, X _{ [0,t]}]\,ds\Big|\,dt =0.
\end{align*}
On the other hand, with Corollary  \ref{res-14b}  we obtain 
\begin{align*}
\Limh E_Q\int_0 ^{ T-h} \Big|h ^{ -1}\int_t ^{ t+h} E_Q[\GaF^Q(u,v)(\Xb_s)\mid X _{ [0,t]}]\,ds - \GaF^Q(u,v)(t,X _{ [0,t]})\Big|\,dt =0.
\end{align*}
The limit \eqref{eq-55} follows from these last two limits.
\end{proof}

\section{About Nelson velocities}
This section refers to the diffusion measure $Q$ of Section \ref{sec-EDM}. Its  content   is not used directly in this article. We propose it to the reader to stress the importance for our purpose  of considering the relative momentum field $ \beta ^{ Q|R}$  rather than the absolute velocity $\vv^Q.$
\\
Denoting by $\Id$ the identity mapping on $\ZZ,$ we see that the vector field $\vv^Q$ appearing in  the martingale problems satisfy $\vv^Q=\LLF^Q[\Id]$. 
Because of the identification $\LLF^Q= L^Q$  which was obtained at Section \ref{sec-SD}, one suspects that  $\vv^Q$    should satisfy
\begin{align*}
\vv^Q_t= L^Q_t[\Id]= \Limh E_Q\left. \Big( \frac{X _{ t+h}-X_t}{h}\,\right|\, X _{ [0,t]}\Big),
\end{align*} 
\emph{whenever this expression is meaningful}. 
The r.h.s.\ of this identity is the forward Nelson velocity of $Q.$ But in general it is  not  well defined, due to a possible lack of integrability.  
In order to give sense to limits of this type  in a general setting, one must introduce integration times and work as in Proposition \ref{res-H02c}.
Next result presents a situation where integration times can be avoided.

\begin{proposition}\label{res-17} Under the hypothesis  \eqref{eqd-04}, suppose that $\aa$ is bounded from above. 
Then, the  limit
\begin{align*}
\vv^{ Q|R}(t, \omega )= \Limh E_Q\left. \Big( \frac{X _{ t+h}-X_t}{h} - \frac{1}{h} \int _{ [t,t+h]} \BBf(\Xb_s)\,ds\,\right|\, X _{ [0,t]}=\omega _{ [0,t]}\Big),
\quad (t, \omega)\in\Ob,
\end{align*}
takes place in $L^2(\Qb).$ 

\end{proposition}

\begin{proof}
Under \eqref{eqd-04}, we know that $E _{ \Qb}|\beta ^{ Q|R}|^2_\aa < \infty.$ Because of the assumed upper boundedness of $\aa$, this implies that $E _{ \Qb}|\vv ^{ Q|R}|^2< \infty,$ where $\vv ^{ Q|R}:=\aa\beta ^{ Q|R}.$ Rewrite the assertion: $Q\in\MP(\aa,\BBf+\vv ^{ Q|R})$ as:
\begin{align*}
X _{ t+h}-X_t-\int _{ [t,t+h]} \BBf(\Xb_s)\,ds
	=\int _{ [t,t+h]}\vv ^{ Q|R}(\Xb_s)\,ds +M^Q _{ t+h}-M^Q_t,
\end{align*}
where $M^Q$ is a local $Q$-martingale. The   assumption $E _{ \Qb}|\vv ^{ Q|R}|^2< \infty,$ expressed with the Euclidean norm $|\sbt|$ rather than the Riemannian norm $|\sbt|_\gg ,$   permits us to { apply the convolution Lemma \ref{resh-21} to  $v=a \beta$ componentwise with $p=2$.} The critical step where this is used is Jensen's inequality right below \eqref{eqh-42}. With this at hand, proceeding as in the proof of Proposition \ref{resh-01} leads us to the  announced result. \end{proof}

\begin{remarks}\ \begin{enumerate}[(a)]
\item
In the setting of this proposition, if the Nelson velocity $ L^R[\Id]$ is  ill defined because $E_R\Iii |\BBf_t|\,dt= \infty$, it might happen that  $ L^Q[\Id]$ is also ill defined. Nevertheless, we have: $\Iii |\BBf_t|\,dt<\infty, R\ae,$ and $Q\in\MP(\aa,\vv^Q)$ where $\vv^Q=\BBf+\vv ^{ Q|R}$ satisfies $\Iii|\vv^Q_t|\,dt< \infty, Q\ae$ 

\item
Requiring that the diffusion matrix field $\aa$ is upper bounded is not a strong restriction for the applications, because in general temperature is upper bounded.

\item
If $\aa$ is only locally bounded, then 
there exists   a sequence $(h_n)$ of positive numbers such that   $\Lim n h_n=0$ and  the limit 
\begin{equation*}
\vv^{ Q|R}(t, \omega)
	= \Lim k \Lim n   E_Q\left. \Big( \frac{X  ^{ \tau_k}_{ t+h_n}-X^{\tau_k} _t}{h_n} -\frac 1{h_n}  \int _{ [t,t+h_n]}\1 _{ \left\{ s\le \tau_k\right\} } \BBf(\Xb_s)\,ds\,\right|\, X _{ [0,t]}=\omega _{ [0,t]}\Big)
\end{equation*}
holds $\Qb\ae,$ where for each integer $k\ge 1,$ $ \tau_k:= \inf\{t\in\ii: |X_t|\ge k\}.$
The proof of this statement  is similar to Proposition \ref{res-H02c}'s proof.
\end{enumerate}\end{remarks}


\begin{thebibliography}{10}

\bibitem{Alb00}
S.~Albeverio.
\newblock {\em Theory of {D}irichlet forms and applications, in Ecole d'Et{\'e}
  de Probabilit{\'e}s de Saint-Flour XXX-2000}, volume 1816 of {\em Lecture
  Notes in Mathematics}.
\newblock Springer, Berlin, 2003.

\bibitem{BGN21a}
E.~Bernton, P.~Ghosal, and M.~Nutz.
\newblock Entropic optimal transport: Geometry and large deviations.
\newblock Preprint arXiv:2102.04397.

\bibitem{Bo52}
D.~Bohm.
\newblock A suggested interpretation of the quantum theory in tems of "hidden"
  variables {I, II}.
\newblock {\em Physical Review}, 85:166--179, 180--193, 1952.

\bibitem{CCGL20}
P.~Cattiaux, G.~Conforti, I.~Gentil, and C.~L{\'e}onard.
\newblock Time reversal of diffusion processes under a finite entropy
  condition.
\newblock Preprint arXiv:2104.07708.

\bibitem{CL94}
P.~Cattiaux and C.~L{\'e}onard.
\newblock Minimization of the {K}ullback information of diffusion processes.
\newblock {\em Ann. Inst. H. Poincar{\'e}. Probab. Statist.}, 30:83--132, 1994.

\bibitem{CL96}
P.~Cattiaux and C.~L{\'e}onard.
\newblock Minimization of the {K}ullback information for some {M}arkov
  processes.
\newblock In {\em Seminar on {P}robability, {tome 30} ({U}niv. {S}trasbourg,
  {S}trasbourg, 1996)}, volume 1626 of {\em Lecture Notes in Math.}, pages
  288--311. Springer, Berlin, 1996.

\bibitem{ChZ03}
K.L. Chung and J.C. Zambrini.
\newblock {\em Introduction to random time and quantum randomness}.
\newblock World Scientific Publishing Co. Inc., 2003.

\bibitem{CZ95}
K.L. Chung and Z.~Zhao.
\newblock {\em From Brownian Motion to Schr\"odinger's Equation}, volume 312 of
  {\em Grundlehren der mathematischen Wissenschaften}.
\newblock Springer-Verlag, 1995.

\bibitem{Co18}
G.~Conforti.
\newblock A second order equation for {S}chr{\"o}dinger bridges with
  applications to the hot gas experiment and entropic transportation cost.
\newblock {\em Probability Theory and Related Fields}, 174(1):1--47, 2019.

\bibitem{CF13}
R.~Cont and D.-A. Fourni{\'e}.
\newblock Functional {I}t{\^o} calculus and stochastic integral representation
  of martingales.
\newblock {\em Ann. Probab.}, 41(1):109--133, 2013.

\bibitem{CWZ}
A.B. Cruzeiro, L.~Wu, and J.-C. Zambrini.
\newblock {B}ernstein processes associated with a {M}arkov process.
\newblock In R.~Rebolledo, editor, {\em Stochastic analysis and mathematical
  physics, ANESTOC'98. Proceedings of the Third International Workshop}, Trends
  in Mathematics, pages 41--71, Boston, 2000. Birkh\"auser.

\bibitem{DM4}
C.~Dellacherie and P.-A. Meyer.
\newblock {\em Probabilit{\'e}s et Potentiel. Ch. XII {\`a} XVI. Th{\'e}orie du
  potentiel associ{\'e}e {\`a} une r{\'e}solvante, th{\'e}orie des processus de
  Markov}.
\newblock Hermann. Paris, 1987.

\bibitem{Doob57}
J.L. Doob.
\newblock Conditional {B}rownian motion and the boundary limits of harmonic
  functions.
\newblock {\em Bull. Soc. Math. France}, 85:431--458, 1957.

\bibitem{Doob84}
J.L. Doob.
\newblock {\em Classical Potential Theory and Its Probabilistic Counterpart}.
\newblock Classics in Mathematics. Springer, 2nd edition, 2000.
\newblock (reprint of the 1984 first edition).

\bibitem{DT09}
D.~D{\"u}rr and S.~Teufel.
\newblock {\em Bohmian mechanics. The physics and mathematics of quantum
  theory}.
\newblock Springer-Verlag, Berlin, 2009.

\bibitem{evans-PDE}
L.C. Evans.
\newblock {\em Partial Differential Equations}, volume~19 of {\em Graduate
  Series in Mathematics}.
\newblock American Mathematical Society, 1998.

\bibitem{F94}
S.~Fang.
\newblock In{\'e}galit{\'e} du type de {P}oincar{\'e} sur l'espace des chemins
  riemanniens.
\newblock {\em C. R. Acad. Sci. Paris}, 318:257--260, 1994.

\bibitem{Fey48}
R.~Feynman.
\newblock Space-time approach to nonrelativistic quantum mechanics.
\newblock {\em Rev. Mod. Phys.}, 20:367--387, 1948.

\bibitem{FH65}
R.~Feynman and A.~Hibbs.
\newblock {\em Quantum Mechanics and Path Integrals}.
\newblock International Series in Pure and Applied Physics. McGraw-Hill, 1965.

\bibitem{FS06}
W.H. Fleming and H.M. Soner.
\newblock {\em Controlled Markov Processes and Viscosity Solutions}, volume~25
  of {\em Applications of Mathematics}.
\newblock Springer, second edition, 2006.

\bibitem{Foe85b}
H.~F{\"o}llmer.
\newblock An entropy approach to the time reversal of diffusion processes.
\newblock In {\em Stochastic Differential Systems - Filtering and Control},
  volume~69 of {\em Lecture Notes in Control and Information Sciences}, pages
  156--163. Springer, 1985.

\bibitem{Foe86}
H.~F{\"o}llmer.
\newblock Time reversal on {W}iener space.
\newblock In {\em Stochastic Processes - Mathematics and Physics}, volume 1158
  of {\em Lecture Notes in Math.}, pages 119--129. Springer, Berlin, 1986.

\bibitem{Foe85}
H.~F{\"o}llmer.
\newblock {\em Random fields and diffusion processes, in \'Ecole d'\'et{\'e} de
  Probabilit{\'e}s de Saint-Flour XV-XVII-1985-87}, volume 1362 of {\em Lecture
  Notes in Mathematics}.
\newblock Springer, Berlin, 1988.

\bibitem{Hsu97}
E.~Hsu.
\newblock Analysis on path and loop spaces.
\newblock volume~5 of {\em IAS/Park City Math. Series}. Amer. Math. Soc., 1997.

\bibitem{Jac79}
J.~Jacod.
\newblock {\em Calcul stochastique et probl{\`e}mes de martingales}, volume 714
  of {\em Lecture Notes in Mathematics}.
\newblock Springer, 1979.

\bibitem{Kac49}
M.~Kac.
\newblock On the distribution of certain {W}iener functionals.
\newblock {\em Trans. Amer. Math. Soc.}, 65:1--13, 1949.

\bibitem{Kac51}
M.~Kac.
\newblock On some connections between probability theory and differential and
  integral equations.
\newblock In J.~Neyman, editor, {\em Proc. Second Berkeley Symp. Math. Stat.
  Prob.}, pages 189--215. Univ. of California Press, 1951.

\bibitem{Kh59}
R.~Khas'minskii.
\newblock On positive solutions of the equation ${A}u+{V}u=0$.
\newblock {\em Theory Probab. Appl.}, 4:309--318, 1959.

\bibitem{Kry87}
N.V. Krylov.
\newblock {\em Nonlinear Elliptic and Parabolic Equations of the Second Order}.
\newblock Reidel, Dordrecht, 1987.

\bibitem{K69}
H.~Kunita.
\newblock Absolute continuity of {M}arkov processes and generators.
\newblock {\em Nagoya Mathematical Journal,}, 36:1--26, 1969.

\bibitem{Kun97}
H.~Kunita.
\newblock {\em Stochastic flows and stochastic differential equations},
  volume~24 of {\em Cambridge Studies in Advanced Mathematics}.
\newblock Cambridge University Press, 1997.

\bibitem{Leo11a}
C.~L{\'e}onard.
\newblock Girsanov theory under a finite entropy condition.
\newblock In {\em S{\'e}minaire de probabilit{\'e}s, vol. 44.}, pages 429--465.
  Lecture Notes in Mathematics 2046. Springer, 2012.

\bibitem{Leo12b}
C.~L{\'e}onard.
\newblock Some properties of path measures.
\newblock In {\em S\'eminaire de probabilit\'es, vol. 46.}, pages 207--230.
  Lecture Notes in Mathematics 2123. Springer, 2014.

\bibitem{Leo12e}
C.~L{\'e}onard.
\newblock A survey of the {S}chr{\"o}dinger problem and some of its connections
  with optimal transport.
\newblock {\em Discrete Contin. Dyn. Syst. A}, 34(4):1533--1574, 2014.

\bibitem{Lowther-cadlag}
G.~Lowther.
\newblock Cadlag modifications. \emph{{A}lmost {S}ure - {A} random mathematical
  blog}.
\newblock \\ {\tt \scriptsize
  https://almostsuremath.com/2009/12/18/cadlag-modifications/}.

\bibitem{MR92}
Z.-M. Ma and M.~R{\"o}ckner.
\newblock {\em Introduction to the theory of (non-symmetric) {D}irichlet
  forms}.
\newblock Universitext. Springer, 1992.

\bibitem{MZ84}
P.-A. Meyer and W.A. Zheng.
\newblock Tightness criteria for laws of semimartingales.
\newblock {\em Ann. Inst. H. Poincar{\'e}. Probab. Statist.}, 20(4):353--372,
  1984.

\bibitem{MZ85}
P.-A. Meyer and W.A. Zheng.
\newblock Construction de processus de {N}elson r{\'e}versibles.
\newblock In {\em S\'eminaire de probabilit\'es. Tome 19}, volume 1123 of {\em
  Lecture Notes in Mathematics}, pages 12--26. Springer, 1985.

\bibitem{Nel67}
E.~Nelson.
\newblock {\em Dynamical theories of {B}rownian motion}.
\newblock Princeton University Press, 1967.

\bibitem{O13}
Y.~Oshima.
\newblock {\em Semi-{D}irichlet forms and {M}arkov processes}, volume~48 of
  {\em de Gruyter Studies in Mathematics}.
\newblock Walter de Gruyter \&\ Co., 2013.

\bibitem{RY99}
D.~Revuz and M.~Yor.
\newblock {\em Continuous martingales and Brownian motion}, volume 293 of {\em
  Grundlehren der Mathematischen Wissenschaften}.
\newblock Springer, 3rd edition, 1999.

\bibitem{Sch31}
E.~Schr\"odinger.
\newblock {\"U}ber die {U}mkehrung der {N}aturgesetze.
\newblock {\em Sitzungsberichte Preuss. Akad. Wiss. Berlin. Phys. Math.},
  144:144--153, 1931.

\bibitem{Sch32}
E.~Schr{\"o}dinger.
\newblock Sur la th{\'e}orie relativiste de l'{\'e}lectron et
  l'interpr{\'e}tation de la m{\'e}canique quantique.
\newblock {\em Ann. Inst. H. Poincar\'e}, 2:269--310, 1932.
\newblock 

\bibitem{St99}
W.~Stannat.
\newblock {\em The theory of generalized {D}irichlet forms and its applications
  in analysis and stochastics}, volume 678 of {\em Mem. Amer. Math. Soc.}
\newblock American Mathematical Society, 1999.

\bibitem{SV69a}
D.W. Stroock and S.R.S. Varadhan.
\newblock Diffusion processes with continuous coefficients, {I}.
\newblock {\em Communications on Pure and Applied Mathematics}, 22(3):345--400,
  1969.

\bibitem{SV69b}
D.W. Stroock and S.R.S. Varadhan.
\newblock Diffusion processes with continuous coefficients, {II}.
\newblock {\em Communications on Pure and Applied Mathematics}, 22(4):479--530,
  1969.

\bibitem{SV79}
D.W. Stroock and S.R.S. Varadhan.
\newblock {\em Multidimensional diffusion processes}.
\newblock Number 233 in Grundlehren der mathematischen Wissenschaften. Springer
  Verlag, 1979.

\bibitem{vR11}
M.~von Renesse.
\newblock An optimal transport view on {S}chr{\"o}dinger's equation.
\newblock {\em Canad. Math. Bull.}, 55(4):858--869, 2011.

\bibitem{Zam86}
J.-C. Zambrini.
\newblock Variational processes and stochastic versions of mechanics.
\newblock {\em J. Math. Phys.}, 27:2307--2330, 1986.

\bibitem{Zam12}
J.-C. Zambrini.
\newblock The research program of stochastic deformation (with a view toward
  geometric mechanics).
\newblock In {\em Stochastic Analysis: A series of lectures}. Springer, 2015.
\newblock arXiv:1212.4186.

\bibitem{Z09}
X.~Zhang.
\newblock Clark-{O}cone formula and variational representation for {P}oisson
  functionals.
\newblock {\em Ann. Probab.}, 37(2):506--529, 2009.

\bibitem{Z85}
W.A. Zheng.
\newblock Tightness results for laws of diffusion processes. {A}pplication to
  stochastic mechanics.
\newblock {\em Ann. Inst. H. Poincar{\'e}. Probab. Statist.}, 21(2):103--124,
  1985.

\end{thebibliography}

\end{document}